\documentclass[a4paper,11pt]{article}
\usepackage{amsmath}
\usepackage{amssymb}
\usepackage{epsfig}
\usepackage{rotating} 
\usepackage{framed}
\usepackage{float}
\usepackage{subcaption}
\usepackage[dvipsnames]{xcolor}

\usepackage{multirow,longtable}
\usepackage[latin1]{inputenc}
\usepackage{amsfonts}
\usepackage{latexsym}
\usepackage{pstricks}
\usepackage{graphicx}
\usepackage{epstopdf}
\usepackage{chngcntr}
\counterwithin{figure}{section}
\counterwithin{table}{section}
\usepackage[comma]{natbib}
\newcommand{\dem}{\par \noindent{\bf Proof:} \par}
\newcommand{\fin}{\hfill $\square$  \par \bigskip}

\usepackage{verbatim}
\usepackage[ruled,vlined,linesnumberedhidden]{algorithm2e}
\usepackage[ruled]{algorithm2e}
\usepackage{caption}
\makeatletter
\newcounter{megaalgorithm}
\newenvironment{megaalgorithm}[1][htb]
{
	\let\c@algocf\c@megaalgorithm
	\begin{algorithm}[#1]%
	}{\end{algorithm}}
\makeatother

\usepackage{algorithm2e}
\usepackage{chngcntr}

\newtheorem{teor}{\bf Theorem}[section]

\newtheorem{coro}{\bf Corollary}[section]

\newtheorem{rem}{\bf Remark}[section]

\begin{document}

\title{\LARGE  Mixed Integer Linear Programming for Feature Selection in Support Vector Machine}

\author{Martine Labb\'e$^{1}$, Luisa I. Mart\'{\i}nez-Merino$^2$, Antonio M.\ Rodr\'{\i}guez-Ch\'{\i}a$^2$}

\date{\scriptsize
$^1${ Département d'Informatique, Université Libre de Bruxelles, Brussels, Belgium and INOCS, INRIA Lille Nord-Europe, Lille, France}\\
$^2$Departamento de Estad\'{\i}stica e Investigaci\'on Operativa, Universidad de C\'adiz, Spain\\\today}

\maketitle 

\begin{abstract}
This work focuses on support vector machine (SVM) with feature selection.
A MILP formulation is proposed for the problem.
The choice of suitable features to construct the separating hyperplanes has been modelled in this formulation by including a budget constraint 
that sets in advance a limit on the number of features to be used in the classification process. We propose both an exact and a heuristic procedure to solve this formulation in
an efficient way. Finally, the validation of the model is done by checking it with some well-known data sets and comparing it with classical classification methods. 
\end{abstract}
\section{Introduction}
{ 
In supervised classification, we are given a set of objects $\Omega$ partitioned into classes and the goal is to build 
a procedure for classifying new objects into these classes. This type of problem is faced in many fields, including insurance companies (to determine whether an applicant is a high insurance risk or not), banks (to decide whether an applicant is a credit 
risk or not), medicine (to determine whether a tumor is benign or malignant), etc. This wide field of application has attracted the attention of a number of 
researchers from different areas. Currently, these problems are analyzed from different perspectives: artificial intelligence,
machine learning, optimization or statistical pattern recognition among others. In this paper we analyze these problems from the point of view of optimization
and more specifically, from the perspective of Mathematical Programming \cite{Mang65,Mang68}.

In the partitioning process, the objects of $\Omega$ are considered as points in an $n$-dimensional feature space. However, the number of features
is often much larger than the number of objects in the population. Handling such a high number of features could therefore be a difficult task and in addition, 
the interpretation of the results  could also be impossible. In this sense, feature selection consists of eliminating as many features as possible in a given problem, keeping an acceptable accuracy where spurious random  variables are removed. Actually, having a minimal number of features often leads to better generalization and simpler models that can be interpreted more easily.

The support vector machine (SVM) is a type of mathematical programming approach developed 
by \cite{Vap98}, and \cite{cova95}. It has been widely studied and has become popular in many fields of application in recent
years; see the introductory
description of SVM by \cite{Burg98}. The SVM is based on margin maximization, which consists of finding the separating
hyperplane that is farthest from the closest object. SVM has proven to be a very powerful tool for supervised classification. Recently, \cite{Maldo14}
proposed two SVM-based models in which feature selection is taken into account by introducing a budget constraint in the formulation, limiting the
number of features used in the model, see also \cite{Ay15} and \cite{GAUDIOSO2017}. 

In this work, we are proposing a MILP formulation, based on \cite{Maldo14} idea to choose the best features and to obtain an adequate predictor. To find an efficient way for solving
this model, we exploit the tightening of the bounds on the separating hyperplane coefficients, enabling us to get better times. Exact and heuristic solution approaches are also presented using these improvements. Lastly, the model is validated by
comparing the proposed formulation with other models and with other feature selection techniques  known in the literature, such as  Recursive Feature Elimination (RFE) or the Fisher Criterion Score (F); see \cite{Guyon02} and \cite{Guyon06} respectively. 

The paper is organized into 8 sections. Section \ref{svm} is a revision of various formulations for SVM-based models that are analyzed in the 
literature. Section \ref{themodel} presents the model under study. In Section \ref{strategies}, several strategies to fix some of 
the ``big M'' parameters of the model are considered. Sections \ref{heuristic} and \ref{exact} develop heuristic and exact solution approaches, respectively.
Section \ref{computational} presents some of the computational results that illustrate the improvement analyzed in the paper. Section \ref{validation} 
 is devoted to analyzing the validation of the model presented in the paper. Finally some conclusions are addressed.

\section{Support Vector Machine}\label{svm}
{
Consider a training set $\Omega$ partitioned into two classes, each object $i\in\Omega$ is represented with a pair
$({\bf x}_i,y_i)\in\mathbb{R}^n\times\{-1,1\}$, where $n$ is the number of features analyzed over each element of $\Omega$, ${\bf x}_i$
contains the features' values and $y_i$ provides the labels, $1$ or $-1$, associated with the two classes in $\Omega$. 
The SVM determines a hyperplane $f(x)= w^T\cdot x +b$ that optimally separates the  training examples.
In the case of linearly separable data, this hyperplane maximizes the margin between the two data classes, i.e.,
it maximizes the distances between two parallel hyperplanes supporting some elements of the two classes. Even if the training
data is non-linearly separable, the constructed hyperplane also minimizes  classification errors. Thus, the classical SVM model minimizes an objective function that is a compromise between the structural risk, given by the inverse of the margin, $\|\omega\|$, and the empirical risk, given by the deviation of misclassified objects. Several SVM models have been proposed using different measures of margin and deviation. Among them, the standard $\ell_2$-SVM (\cite{BraMa98}) uses the following formulation:
\begin{eqnarray}
(\ell_2\textrm{-SVM})\quad &\displaystyle \min_{w,b,\xi}& \frac{1}{2}\|w\|_2^2 + C \sum_{i=1}^m \xi_i \nonumber\\
&s.t.& y_i (w^Tx_i+b) \ge 1-\xi_i, \quad \forall i=1,\ldots,m,\label{error} \\
& &\xi_i \ge 0,\quad \forall i=1,\ldots,m.  \label{bounds}
\end{eqnarray}
As can be observed, $\ell_2$-SVM  considers $\ell_2$-norm to measure the margin and it introduces the slack variables $\xi_i$ $i=1,\ldots,m$ to measure the deviations of misclassified elements. Additionally, a penalty parameter $C$ that regulates the trade-off between structural and empirical risk is added. Constraints \eqref{error} are the main restrictions appearing in classical SVM. In fact, constraints \eqref{error} determine whether or not the training data are separable by the classifier hyperplane.  

\cite{BraMa98} also presented the same model but considering $\ell_1$-norm instead of $\ell_2$-norm for the  margin. The resulting model is the following,
\begin{eqnarray}
&\displaystyle \min_{w,b,\xi}& \|w\|_1 + C \sum_{i=1}^m \xi_i \nonumber\\
&s.t.&\textrm{\eqref{error}-\eqref{bounds}}\nonumber.
\end{eqnarray}
An equivalent linear formulation of this problem is:
\begin{eqnarray}
(\ell_1\textrm{-SVM})\quad& \displaystyle \min_{w,b,\xi} & \sum_{j=1}^n  z_j + C \sum_{i=1}^m \xi_i\nonumber \\
&s.t.&\textrm{\eqref{error}-\eqref{bounds}}\nonumber,\nonumber\\
& & -z_j \le w_j \le z_j, \quad \forall j=1,\ldots,n, \label{abs}\\
& & z_j \ge 0, \quad \forall j=1,\ldots,n. \label{domz}
\end{eqnarray}
Due to constraints \eqref{abs}, the variables $z$ represent the absolute value of the hyperplane coefficients $w$.

Often, real data are composed of few sample elements ($m$), but each element has a large number of related features ($n$). Therefore, it is essential to select a suitable set of features to construct the classifier. From among the different techniques for feature selection, we focus on the embedded methods that perform feature selection at the same time as the classifier is constructed. Specifically, we focus on the SVM models that include feature selection constraints.

\cite{Maldo14} proposed a model inspired by the $\ell_1$-SVM in which the idea of feature selection was introduced through a budget constraint. Unlike $\ell_1$-SVM, the objective function of this model does not consider the margin, i.e., this model focuses on minimizing the sum of the deviations. Although structural risk is not explicitly included in the objective function it is, in some way, implicitly under control because the number of non-null $w$-variables is bounded by a budget constraint. This model is based on the use of a binary variable linked to each feature in order to restrict the
number of attributes used in the classifier via a budget constraint. 
 A cost vector $c\in\mathbb{R}^n$ is considered, where $c_j$ is the cost of acquiring 
attribute $j$, $j=1,\ldots,n$. The formulation is therefore given by
\begin{eqnarray}
\textrm{(MILP1)}\quad&\displaystyle \min_{{ v,w,b,\xi}} & \sum_{i=1}^m \xi_i\nonumber \\
& s.t. &\textrm{\eqref{error}-\eqref{bounds}},\nonumber\\ 
& & l_j v_j \le w_j \le u_j v_j, \quad \forall j=1,\ldots,n, \label{boundsw} \\
& & \sum_{j=1}^n c_j v_j \le B, \label{budget}\\
& & v_j \in \{0,1\}, \quad \forall j=1,\ldots,n. \label{binary}
\end{eqnarray}
 Constraints \eqref{boundsw} link the $v$- and $w$-variables and enable the identification of the $w$-variables which are non-null, i.e. $w_j$ will be non-null only if $v_j$ takes value $1$ for any $j=1,\ldots,n$. In fact, these are big M constraints.
 Values $l_j$ and $u_j$ correspond to the lower and upper bounds of the value of $w_j$, $j=1,\ldots,n$, respectively. As previously mentioned, constraint \eqref{budget} is the budget constraint that limits the number of non-null $w$-variables. Thus, an important issue for solving this model is the appropriate choice of these bounds because the efficiency of any enumeration solution approach will greatly depend on the tightness of the model's LP relaxation.
 
 A possible criticism to this model is that in the case in which, for a particular value of $B$, there are many optimal solutions with objective value $0$, this model does not provide a way to choose among those hyperplanes. This situation is very often in datasets with many features and very few objects.

\section{The model}\label{themodel}
Based on the idea introduced by \cite{Maldo14}, we propose the extension of $\ell_1$-SVM with a budget constraint, i.e., our model  takes into account the structural and empirical risk in the objective function with feature selection through a budget constraint, therefore we avoid the above mentioned criticism of MILP1.
\cite{Maldo14} were mainly focused on
validating their model by contrasting it against well-known classification methods from the literature, but little attention was paid to analyze how to solve their problem efficiently. However, our goal in this paper is to provide a deep analysis of the model, allowing us to produce  efficient  exact and heuristic solution approaches  in addition to validating the model by comparing it with classical classification methods. Hence the model that we propose is given by,
\begin{eqnarray}
&\displaystyle \min_{v,w,b,\xi,z} & 
 \sum_{j=1}^n z_j+C\sum_{i=1}^m \xi_i  \nonumber\\
& s.t. & \textrm{\eqref{error}-\eqref{binary}}.\nonumber
\end{eqnarray}
 In contrast to MILP1, the presented formulation considers the margin in the objective function. Thus, the problem looks for an optimal balance between deviations and the margin using $\ell_1$-norm.
In what follows and for the sake 
of clarity, analogously to \cite{Maldo14}, we will assume that $c_j$ in \eqref{budget} is set to $1$ for $j=1,\ldots,n$. Therefore, the budget $B$ will represent the maximum number of features that can be selected.  An equivalent formulation for this model is obtained by decomposing the unrestricted variables $w_j$ as two different non-negative variables, $w_j^+$ and $w_j^-$. In this reformulation, $w_j=w_j^+-w_j^-$, where $w_j^+,w_j^-\geq 0$ for $j=1,\ldots,n$. Thusly, by taking advantage of this definition, we have that $z_j=|w_j|=w_j^++w_j^-$ in any optimal solution since $w_j^-+w_j^+$ for $j=1,\ldots,n$ is part of the objective function to be minimized. This means that, at most, only one of the two variables is non-zero in the optimal solution.
Consequently, the following formulation is obtained,
\begin{eqnarray}
\textrm{(FS-SVM)}&\displaystyle \min_{v,w^+,w^-,b,\xi,z} & 
 \sum_{j=1}^n (w_j^++w_j^-)+C\sum_{i=1}^m \xi_i \nonumber\\
&s.t.& y_i \Big( \sum_{j=1}^n (w_j^+-w_j^-) x_{ij} +b\Big) \ge 1-\xi_i, \quad \forall i=1,\ldots,m,\label{planes} \\
&& \sum_{j=1}^n v_j\leq B,\label{numv}\\
& &  w_j^+ \le  u_j v_j, \quad \forall j=1,\ldots,n, \label{waconst} \\
& &  w_j^- \le -l_j v_j, \quad \forall j=1,\ldots,n, \label{wbconst} \\
& &   w_j^+\geq0,w_j^-\geq 0,\quad \forall j=1.\ldots,n,\label{wbounds}\\
& & \xi_i \ge 0, \quad \forall i=1,\ldots,m,\label{limxi}\\
& &  v_j \in \{0,1\}, \quad \forall j=1,\ldots,n. 
\end{eqnarray}
Note that the FS-SVM formulation presents a feature selection constraint \eqref{numv} that limits the number of selected features in order to construct the separating hyperplane. Additionally, constraints \eqref{waconst} and \eqref{wbconst} are two sets of big M constraints.
 
 A preliminary computational study to check how difficult it is to solve the aforementioned mixed integer linear formulation of FS-SVM with very conservative big M values ($u_j$ and $l_j$ for all $j=1,\ldots,n$) shows that  formulation's performance is not very good (see the columns FS-SVM of Table \ref{tab_prep} for the different data sets). This encourages us to check whether
a strengthening the big M values might improve these computational results. In the following sections, we will analyze the influence of tightening bounds of the $w$-variables in this formulation, $l$ and $u$, to solve the model.  In this sense, we will develop different methodologies to obtain better $w$-variable bounds.

In addition, we have also studied an alternative formulation of FS-SVM by substituing constraints \eqref{waconst} and \eqref{wbconst} by conditional constraints and implementing them by CPLEX command IloIfThen, however we have omitted it because we obtained very bad computational times.
Moreover, from this preliminar computational analysis we have checked that the solution times for solving MILP1 and FS-SVM are similar. Only in the cases in which the optimal value of MILP1 is 0, this model is much faster, but in those cases, MILP1 is useless because the data are separable for the chosen features and many separating hyperplanes can be equally valid.
 
\section{Strategies for obtaining tightened values of the $w_j$ bounds}\label{strategies}

As mentioned above and in terms of developing good solution approaches to our problem, it would be useful to provide tightened values of the upper/lower bounds of $w_j$ for $j=1,\ldots,n$. It should be noted that the literature contains various methods related to bound reduction of $w$ variables in SVM. In particular, two methods are developed in \cite{BELOTTI2016}. One of them was the origin of a CPLEX parameter and the other is based on an iterative process that solves auxiliary MIPs to strengthen big M values associated with certain variables. In our preliminary computational analysis we checked this parameter and it did not improve our computational results. In addition, the second approach in \cite{BELOTTI2016} consists in an iterative process that solves a sequence of MIPs (two for each $w_j,$ $j=1,\ldots,n$). They applied this approach to data with $m=100$ and $n=2$, for this reason, they solve four MILPs in each iteration. However for the datasets analyzed in this paper with a large number of features, this approach does not make sense.

For FS-SVM, we develop two strategies to compute the bounds of $w_j$ for any $j=1,\ldots,n$.  The first strategy proposed is based on solving the maximization of linear problems that report the lower/upper bounds of the variables and the second one uses the Lagrangian relaxation to tighten the bounds.
Note that, in what follows, we will denote the linear relaxation of FS-SVM as LP-FS-SVM.

\subsection{Strategy I}
Given a subset $K\subset\{1,\ldots,n\}$, we will denote the  restricted problem below, which is derived from the original FS-SVM, as FS-SVM($K$):

\begin{eqnarray*}
(\textrm{FS-SVM}(K))&\displaystyle \min_{v,w^+,w^-,b,\xi} & 
 \sum_{j\in K} (w_j^++w_j^-)+ C\sum_{i=1}^m \xi_i  \\
& s.t.&  y_i \Big( \sum_{j\in K} (w_j^+-w_j^-) x_{ij} +b\Big) \ge 1-\xi_i, \quad \forall i=1,\ldots,m,\nonumber \\
& & w_j^+ \le u_j v_j, \quad \forall j\in K , \\
& & w_j^- \le -l_j v_j, \quad \forall j\in K, \\
& & \sum_{j \in K} v_j \le B, \\
& &  v_j \in\{0,1\} , w_j^+,w_j^-\ge 0 \quad \forall j \in K, \\
& &
	 \xi_i \ge 0, \quad \forall i=1,\ldots,m. 
\end{eqnarray*}
Note that in this problem only a subset of variables $v$, $w^+$ and $w^-$ are considered. This is equivalent to considering the FS-SVM where
$v_j=w_j^+=w_j^-=0$ for $j\in\{1,\ldots,n\}\setminus K$. Consequently, the solution to this problem is feasible for the original problem and its objective value, called UB, is an upper bound of our model. Solving FS-SVM$(K)$ will be necessary in the application of Strategy I and
it will also be used in the heuristic approach, as we will see in Section \ref{heuristic}. 
The process given by Strategy I is described in Algorithm \ref{algorithm01}.
\\ \mbox{}\\ 
\scalebox{0.85}{
\begin{algorithm}[H]
 \KwData{Training sample composed by a set of $m$ elements with $n$ features.}
 \KwResult{ Updated values of upper bounds parameters $l_j,u_j$ for $j=1,\ldots,n$.}
\tcc{ Step 1}
For $j=1,\ldots,n$, let $\overline{w}_j^+$ and $\overline{w}_j^-$ be an optimal solution for LP-FS-SVM and set $K_0:=\{j\: : \: \overline{w}_j^++\overline{w}_j^->0\}$. Solve the restricted problem FS-SVM$(K_0)$ to obtain UB.\\
\tcc{ Step 2}
 \For{$k=1$ \KwTo $k=n$}{
 Solve the following linear programming  problems for $k=1,\ldots,n$.
\begin{eqnarray*}
&\mbox{(LP)}\displaystyle \max_{v,w^+,w^-,b,\xi} & w_k^++w_k^- \\ 
& s.t. & \textrm{ \eqref{planes}-\eqref{limxi}},\\
& & \sum_{j=1}^n (w_j^++w_j^-) + C\sum_{i=1}^m \xi_i \le \textrm{UB}, \\
& & 0\le v_j \le1, \quad \forall j=1,\ldots,n. 
\end{eqnarray*}
Let $\overline{u}_k$ be the optimal value of the above problem.\\
\If {{$\overline{u}_k<\max\{-l_k,u_k\}$}}{$ u_k:=\min\{u_k,\overline{u}_k\}$, $-l_k:=\min\{-l_k,\overline{u}_k\}$.}
 }
\label{algorithm01}
 \caption{Strategy I}
\end{algorithm}
}
{
	\begin{rem}
		Note that the bounds obtained by Algorithm \ref{algorithm01} can be  improved by substituting constraint \eqref{numv} in the (LP) for any $k=1,\ldots,n$ for $\sum_{j \neq k} v_j \le B-1$. 
		However, the computational analysis  addressed in a preliminary study showed that the improvement in quality of the bounds is very 
		small and the running times increased when using this modification. For this reason, we decided to keep constraint \eqref{numv} and not use this modification.
	\end{rem}
}
\subsection{Strategy II}
Unlike the previous strategy in which bounds for $w_j^++w_j^-$ have been computed, this strategy will provide us with bounds for $w_j^+$ and $w_j^-$ independently. 
In this case, the strategy is based on the results below.
{
	\begin{teor}
		Let $(\overline{v}, \overline{w}^{+}, \overline{w}^{-}, \overline{b}, \overline{\xi})$ be an optimal solution of LP-FS-SVM; $z_{LB}^{LP}$  
		its objective value; $\overline{\alpha}$ a vector of optimal values for the
		dual variables associated with the constraints \eqref{planes}; and $\overline{w}_{j_0}^++\overline{w}_{j_0}^-=0$ for some $j_0\in\{1,\ldots,n\}$.
		\begin{itemize}
			\item[i)] If   $(\tilde{v}, \tilde{w}^{+}, \tilde{w}^-, \tilde{b}, \tilde{\xi})$ is an optimal solution of LP-FS-SVM restricting
			$w_{j_0}^+=\tilde{w}_{j_0}^+$ where $\tilde{w}_{j_0}^+$ is a positive constant, $z_{\tilde{w}_{j_0}^+}$ its objective value and 
			$\sum_{j=1}^n\overline{v}_j+\frac{\tilde{w}_{j_0}^+}{u_{j_0}}\leq B$, then
			\begin{equation*}
			z_{LB}^{LP}+\tilde{w}_{j_0}^+(1-\sum_{i=1}^{m} \overline{\alpha}_i y_ix_{ij_0})\leq z_{\tilde{w}_{j_0}^+}.
			\end{equation*}
			\item[ii)] 	If $(\tilde{v}, \tilde{w}^{+}, \tilde{w}^-, \tilde{b}, \tilde{\xi})$ is an optimal solution of LP-FS-SVM restricting
			$w_{j_0}^-=\tilde{w}_{j_0}^-$ with ${\tilde w}_{j_0}^-$ a positive constant, $z_{\tilde{w}_{j_0}^-}$  its objective value and 
			$\sum_{j=1}^n\overline{v}_j+\frac{\tilde{w}_{j_0}^-}{-l_{j_0}}\leq B$, then 
				\begin{equation*}z_{LB}^{LP}+\tilde{w}_{j_0}^-(1+\sum_{i=1}^{m}{ \overline{\alpha}_i}y_ix_{ij_0})\leq z_{\tilde{w}_{j_0}^-}.\end{equation*} 
	       \end{itemize} 
		\label{teo1}\end{teor}}
	\dem We are only addressing statement i) here because statement ii) would be proved in a similar manner.
 Since $\overline{\alpha}$ is a vector of optimal values for the dual variables associated with the family of constraints \eqref{planes}, it holds that 
\begin{equation*}
z_{LB}^{LP}=
\sum_{j=1}^n(\overline{w}_j^++\overline{w}_j^-)+C\sum_{i=1}^m\overline{\xi}_i+ 
\sum_{i=1}^m\overline{\alpha}_i\left(1-\overline{\xi}_i-y_i\sum_{j=1}^n(\overline{w}_j^+-\overline{w}_j^-)
x_{ij}-y_i\overline{b}\right).
\end{equation*}
In addition, since $\overline{w}_{j_0}^++\overline{w}_{j_0}^-=0$, we have that 
\begin{equation}
 z_{LB}^{LP}=
\sum_{j=1,j\neq j_0}^n(\overline{w}_j^++\overline{w}_j^-)+ C\sum_{i=1}^m\overline{\xi}_i+
\sum_{i=1}^m\overline{\alpha}_i\left(1-\overline{\xi}_i-y_i\sum_{j=1,j\neq j_0}^n(\overline{w}_j^+-\overline{w}_j^-)
x_{ij}-y_i\overline{b}\right).\label{aux}
\end{equation}
On the other hand, consider the LP-FS-SVM with the additional constraints $w_{j_0}^+=\tilde{w}_{j_0}^+$, $v_{j_0}=\frac{\tilde{w}_{j_0}^+}{u_{j_0}}$ and where the family of constraints \eqref{planes} has been dualized, i.e., 
\begin{eqnarray}
\min_{ v,w^+,w^-,b,\xi} &&
 \sum_{j=1}^n(w_j^++w_j^-)+ C\sum_{i=1}^m \xi_i+
\sum_{i=1}^m\alpha_i \left(1-\xi_i-y_i\sum_{j=1}^n(w_j^+-w_j^-)
x_{ij}-y_i b\right)\nonumber\\
s.t.&& \textrm{\eqref{numv}-\eqref{limxi}},\nonumber\\
&& w_{j_0}^+- w_{j_0}^-=\tilde{w}_{j_0}^+,\nonumber\\
&& 0\leq v_j \leq 1,\;j=1,\ldots,n\nonumber,
\end{eqnarray}
where $\alpha_i\geq 0$. Hence, this problem can be rewritten as follows,
\begin{eqnarray}
\textrm{(Lag-FS-SVM)}\quad\min_{ v,w^+,w^-,b,\xi} &&
\sum_{j=1,j\neq j_0}^n(w_j^++w_j^-)+ C\sum_{i=1}^m \xi_i+
\sum_{i=1}^m\alpha_i ( 1-\xi_i- \nonumber\\
&&-y_i\sum_{j=1,j\neq j_0}^n(w_j^+-w_j^-)
x_{ij}-\sum_{i=1}^m y_i b )+\tilde{w}_{j_0}^+(1-\sum_{i=1}^m \alpha_iy_ix_{ij_0})\nonumber\\
{ s.t}.&&{\sum_{j=1,j\neq j_0}^n v_j\leq B-\frac{\tilde{w}_{j_0}^+}{{ u_{j_0}}},}\nonumber\\
&&w_j^+\leq u_j v_j,\,\,j=1,\ldots,n,\,j\neq j_0,\nonumber\\
&&w_j^-\leq { -l_j} v_j,\,\,j=1,\ldots,n,\,j\neq j_0,\nonumber\\
&&w_j^+,w_j^-\geq 0,\,\,j=1,\ldots,n,\,j\neq j_0,\nonumber\\
&&0\leq v_j\leq 1,\,\,j=1,\ldots,n,\,j\neq j_0,\nonumber\\
&&{ \xi_i\geq 0,\,\,i=1,\ldots,m.}\nonumber
\end{eqnarray}
 Note that $(\overline{v},\overline{w}^+,\overline{w}^-,\overline{b},\overline{\xi})$, an optimal solution of LP-FS-SVM, is feasible for the problem above if $\tilde{w}_{j_0}\leq u_{j_0}(B-\sum_{j=1}^n \bar{v}_j)$. In addition, any feasible solution of the problem Lag-FS-SVM taking $w_{j_0}^+=w_{j_0}^-=v_{j_0}=0$ is feasible for the LP-FS-SVM where family of constraints \eqref{planes} has been dualized. 

Hence, for $\alpha=\overline{\alpha}$, using \eqref{aux}, the optimal objective value of the above problem is 
$z_{LB}^{LP}+\tilde{w}_{j_0}^+(1-\sum_{i=1}^m { \overline{\alpha}}_iy_ix_{ij_0})$, which is the lower bound
of the optimal value of the LP-FS-SVM with the additional constraint of $w_{j_0}^+=\tilde{w}_{j_0}^+$.
\fin
{
\begin{coro}
Under the hypothesis of Theorem \ref{teo1}, if we have an upper bound UB of FS-SVM, then it holds that
\begin{itemize}
	\item[i)] $w_{j_0}^+\leq \min\left\{\frac{\textrm{UB}-z_{LB}^{LP}}{1-\sum_{i=1}^m \overline{\alpha}_i y_i x_{ij_0}},{ u_{j_0}(B-\sum_{j=1}^n\overline{v}_j)}\right\}.$
	\item[ii)] $w_{j_0}^-\leq\min\left\{\frac{\textrm{UB}-z_{LB}^{LP}}{1 +\sum_{i=1}^m \overline{\alpha}_i y_i { x_{ij_0}}}, { -l_{j_0}(B-\sum_{j=1}^n\overline{v}_j)}\right\}.$
\end{itemize} 
\end{coro}}
A detailed description of this second strategy can be found in Algorithm \ref{algorithm02}.

\scalebox{0.85}{
		\begin{algorithm}[H]
			\KwData{Training data ($m$ elements $\times$ $n$ features).}
			\KwResult{Tightened bounds of $u$ and $l$.  }
			\tcc{ Step 1} 
			Solve the { LP-FS-SVM} and obtain  the dual variables (denoted by $\alpha_i$) associated with the family of constraints \eqref{planes}. Let  { $z_{LB}^{LP}$} be its  optimal value and    $(\overline{v}, \overline{w}^{a}, \overline{w}^-, \overline{b}, \overline{\xi}) \in  \mathbb{R}^n_+\times \mathbb{R}^n_+\times \mathbb{R}^n_+\times \mathbb{R} \times \mathbb{R}^m_+$
			an  optimal solution.\\
			\tcc{ Step 2}
			Let UB be  an upper bound of our original model (recall that an upper bound was computed in { Strategy I}). If $\overline{w}_{j_0}^++\overline{w}_{j_0}^-=0$, 
				
				{ Set} $\overline{u}_{j_0}^+ := \min\left\{\frac{\textrm{UB}-{ z_{LB}^{LP}}}{1-\sum_{i=1}^m \alpha_i y_i x_{ij_0}},{ u_{j_0}(B-\sum_{j=1}^n\overline{v}_j)}\right\}$, \\
				\If{$\overline{u}_{j_0}^+ < u_{j_0}$ }{
					$u_{j_0}:=\overline{u}_{j_0}^+$.
				}
				{ Set} $\overline{u}_{j_0}^-{:=}\min\left\{\frac{\textrm{UB}-{ z_{LB}^{LP}}}{1 +\sum_{i=1}^m \alpha_i y_i { x_{ij_0}}}, { -l_{j_0}(B-\sum_{j=1}^n\overline{v}_j)}\right\}$,\\
					\If{$\overline{u}_{j_0}^- < -l_{j_0}$ }{
					${ -l_{j_0}:=\overline{u}_{j_0}^-}$.
				}
			\label{algorithm02}
			\caption{Strategy II}
		\end{algorithm}
	}

\section{Heuristic Solution Approach: Kernel Search}\label{heuristic}
	Among the characteristics of the presented model, we must point out that each data feature ($j$) has an associated binary variable ($v_j$) which indicates 
	whether or not feature $j$ is selected to construct the classifier. Therefore, the size of the problem, and consequently the time required for
	solving it,  grows with the number of features. 
	SVM usually works with real data using quite a large number of features. Hence, a heuristic approach that is suited to the model will help us to very quickly find appropriate, 
	good solutions for those cases where exact methods cannot provide solutions within an acceptable time.
	
	Specifically, we adapt the Kernel Search (KS) proposed by \cite{Ange10}. The basic idea of this heuristic approach is to solve a sequence of restricted MILPs derived from the original problem, thus obtaining a progressively better bound on the solution. The KS has been successfully applied to different kinds of problems such as portfolio optimization 
	(\cite{Ange12}) and location problems (\cite{Gua12}). Even though it was originally applied to pure binary formulations, it has been also used in problems with several continuous or integer variables associated with each binary variable. 
	
	Regarding our model, we observed that the continuous variables $w_j^+$ and $w_j^-$ are related to the binary variables $v_j$ by constraints \eqref{waconst} and \eqref{wbconst}. By applying this heuristic approach to our problem, we will solve a sequence of MILP problems with the same structure as the original one but only considering a subset of variables $v$ and the corresponding subset of continuous variables $\omega^+$ and $\omega^-$ which are associated with it. Since restricted MILPs only take into account a subset of $v$ variables, i.e. the remaining $v$ variables are fixed to $0 $, they will hopefully provide upper bounds in acceptable times.
	
	In the KS, each restricted MILP of the sequence considers the variables that are most likely to take a value different from $0$ in the optimal 
	solution of the original problem. These variables are called promising variables and they form the Kernel set of each restricted MILP. Detailed below is the complete KS procedure for our SVM model, including how to select the promising variables at each step and how to modify the Kernel.
	
	\subsection{Initial step\label{initial_heu}}
	First, feature set $\{1,\ldots,n\}$ must be sorted according to how much the corresponding variables are likely to take a value of $1$ in the optimal solution.
	The LP-FS-SVM is solved with this aim in mind, obtaining a solution $(\bar{v},\bar{\omega}^+,\bar{\omega}^-,\bar{b},\bar{\xi})$
	and the reduced costs of variables $\omega_j^+$ and $\omega_j^-$ for each $j=1,\ldots,n$. Then, features are sorted in non-decreasing order with respect to vector $r$, which is defined as:
	\begin{equation}
	r_{j}=\left\{\begin{array}{c c}
	-(\bar{\omega}_j^++\bar{\omega}_j^-),& \textnormal{if }\bar{\omega}_j^++\bar{\omega}_j^->0,\\
	\min\{r^+_j,r^-_j\},&\textnormal{otherwise}.
	\end{array}
	\right.
	\label{cost_red}
	\end{equation}
	Where, $r_j^+$ and $r_j^-$ are the reduced costs of variables $\omega_j^+$ and $\omega_j^-$ in the LP-FS-SVM,  for $j=1,\ldots,n$.

	To obtain the initial Kernel set ($K_0$), the first $k$ features  are chosen, having been sorted into a non-decreasing order with respect to vector $r$. Specifically, we take 
	$$k:=\left|\left\{j=1,\ldots,n:\bar{\omega}_j^++\bar{\omega}_j^->0\right\}\right|,$$ although $k$ is a parameter of the heuristic that can be modified.
	
	Similar to \cite{Gua12}, the remaining features are divided into $N$ subsets denoted as $K_i$ for $i=1,\ldots,N$. In particular,
	we take $N=\left\lceil \frac{n-k}{k}\right\rceil$. Each subset $K_i$, $i=1,\ldots,N-1$ will be composed of $k$ features and $K_{N}$ will contain the remaining features. In fact, we restrict the KS to analyze only $\bar{N}<N$ of the subsets, due to the size of the instances considered. Computational experiments have shown results when exploring $10\%$ of the total number of subsets, i.e., $\bar{N}:=\left\lceil 0.1\cdot N\right\rceil$.  
	
	Given the initial Kernel $K_0$, the upper bound (UB) of the problem is initialized by solving FS-SVM$(K_0)$. 
	Note that, as observed in Section \ref{strategies}, FS-SVM$(K_0)$ is equivalent to solving the original problem setting $v_j=0$ for $j\notin  K_0$. We should point out that any solution of
	FS-SVM($K_0$) is always a feasible solution for FS-SVM, thus solving FS-SVM($K_0$) we obtain an upper bound.
	
	\subsection{Main step}
	In each iteration ($it$), the
	 heuristic considers the set of features  $\mathcal{K}=K\cup K_{it}$, i.e. the combination of the current Kernel and 
	the features in the set $K_{it}$. To update the UB, in each iteration FS-SVM($\mathcal{K}$) is solved plus the following two constraints, as denoted by 
	FS-SVM($\mathcal{K}$)$+$\eqref{cota}$+$\eqref{nueva},
	
	\begin{eqnarray}
	& & \sum_{j\in\mathcal{K}} (w_j^++w_j^-)+C \sum_{i=1}^m \xi_i \le UB, \label{cota}\\
	& & \sum_{j\in K_{it}}v_{j}\ge 1{ .}\label{nueva} 
	\end{eqnarray}
	Constraint \eqref{cota} restricts the objective function to take a value smaller than or equal to the current upper bound and constraint \eqref{nueva} ensures that at least one feature belonging to $K_{it}$ will be chosen. We also impose the restriction that each problem has to be solved within a time limit of $900$ seconds. If no feasible solution can be found within this time limit, the algorithm skips to the next iteration. Note that this problem may potentially be infeasible due to the presence of constraints \eqref{cota} and \eqref{nueva} together in the formulation. Otherwise, if FS-SVM($\mathcal{K}$)$+$\eqref{cota}$+$\eqref{nueva} is feasible, the objective value, at least, will be equal to the previous UB because of constraint \eqref{cota}.
	\subsection{Update step}
	If the problem FS-SVM($\mathcal{K}$)$+$\eqref{cota}$+$\eqref{nueva} (i.e., FS-SVM($\mathcal{K}$) where constraints \eqref{cota}
	and \eqref{nueva} have been added) is feasible, then some features from $K_{it}$ are chosen in the optimal solution of FS-SVM($\mathcal{K}$)$+$\eqref{cota}$+$\eqref{nueva}. They are added to the current Kernel $K$ for the next iteration
	since adding these features obtains  an identical or better upper bound. Conversely, the set of features of $K$ that has not been chosen in the optimal solution in the previous iterations is removed from the Kernel. The removal of some of the features from the Kernel is decisive in that it does not excessively increase the number of binary variables considered in each FS-SVM($\mathcal{K}$)$+$\eqref{cota}$+$\eqref{nueva}. In our case, we remove the features that were not  selected in the previous two iterations. The set of added features is denoted as $K_{it}^+$ and the set of removed features as $K_{it}^-$. The resulting Kernel for the next iteration is $K=K\cup K_{it}^+\setminus K_{it}^-$.
	
	Coversely, if the problem is infeasible, the kernel is not modified and the procedure skips to the next iteration. The KS for the FS-SVM model is also described in Algorithm \ref{algorithm1}.
	
	\scalebox{0.85}{
		\begin{algorithm}[H]
			\KwData{Training data composed by a set of $m$ elements with { $n$ features.\\ Parameter $k$ is initially fixed as described in Subsection \ref{initial_heu}. }}
			\KwResult{A feasible solution of FS-SVM model. }
			\tcc{ Initial Step}
			Solve LP-FS-SVM. Sort the $n$ features in non-decreasing order with respect to vector $r$ defined in \eqref{cost_red}.\\
			Build the initial kernel $K_0$ taking the first $k$ ordered features. Set $K=K_0$.\\
			Divide the remaining $n-k$ sorted features in a sequence of $N$ subsets ($K_1,\ldots,K_{ N}$).\\
			Select the number of subsets to analyze, $\bar{N}$, {($\bar{N}<N$)}.\\
			Solve the FS-SVM($K$) obtaining the initial upper bound (UB).\\
			\For{$it=0$ \KwTo $it=\bar{N}$}{
				\tcc{ Main Step}
				Build $\mathcal{K}=K\cup K_{it}$.\\
				Solve FS-SVM$(\mathcal{K})$$+$\eqref{cota}$+$\eqref{nueva}.\\
				\tcc{ Update Step}
				\If{FS-SVM($\mathcal{K}$) is feasible (with solution $(\bar{v},\bar{\omega}^+,\bar{\omega}^-,\bar{b},\bar{\xi})$
					and optimal value $z$) and the running time is smaller than $900$ s }{
					$UB=z$\\
					Build $K_{it}^+:=\{j\in K_{it}:\bar{v}_{j}=1\}$.\\
			 Build $K_{it}^-:=\{j\in K:j\textrm{ not selected in the solution of the last two iterations} \}$.\\
					Update $K:=K\cup K_{it}^+\setminus K_{it}^{-}$.\\
				}
				Set $it:=it+1$.\\
			}
			\caption{Kernel Search for FS-SVM \label{algorithm1}}
		\end{algorithm}
	}

 \section{Exact procedure} \label{exact}
This section is devoted to the description of a procedure to get an FS-SVM optimal solution by solving a sequence of semi-relaxed problems.
In this exact procedure, each semi-relaxed problem is associated with a subset $K\subset\{1,\ldots,n\}$ of features in such a way that
only the variables $v_j$ with $j\in K$ will be considered as binary and the remaining ones will be relaxed.
Specifically, the semi-relaxed version of the problem for a set of features $K$ is formulated as follows, 
\begin{eqnarray}
\textrm{(SR-FS-SVM($K$))}&\displaystyle \min_{{ v,w_a,w_b,b,\xi}} &
 \sum_{j=1}^n (w_j^++w_j^-)+C\sum_{i=1}^m \xi_i \nonumber\\
& s.t.& \textrm{\eqref{planes}-\eqref{limxi}},\nonumber\\ 
& & v_j \in \{0,1\}, \quad \forall j\in K, \\
& &0\le v_j\le 1,\quad \forall j \in \{1,\ldots,n\}\setminus K.
\end{eqnarray}
 The optimal value of SR-FS-SVM($K$) provides a lower bound of FS-SVM. By adding and removing certain 
features of $K$, a sequence of semi-relaxed problems is created, providing lower bounds on the solutions. As we will see in the following subsections,  strategies I, II  and the Kernel Search} are also used in this procedure.

\subsection{Initial Step}
To obtain initial bounds on the objective value, the exact procedure exploits techniques detailed previously. First, strategies I and
II are used to tighten parameters $u$ and $l$. It should be noted that the use of these strategies provide an initial upper bound (UB)
and an initial lower bound (LB) for the objective value (by solving the linear relaxation).  The Kernel Search is then performed in order to improve the UB given by the strategies.

\subsection{Main Step}
The main step of the exact procedure consists of solving a sequence of semi-relaxed problems to improve the lower bound of the objective value.
To start with, we must select a subset of features $K\subset\{1,\ldots,n\}$ whose associated $v$ variables will be considered as binary variables in the first semi-relaxed problem. 
 The Kernel Search provides a subset of features that allows us to obtain a good bound
on the optimal objective value. Therefore,  the exact procedure will consider the set provided by the heuristic as the initial $K$ 
and it will obtain an initial LB  solving SR-FS-SVM($K$).

Then, the set $K$ is updated by adding and removing some of the features, improving the bound of the objective value. To this end,
two sets (denoted by $K^+$ and $K^-$) are built in each iteration. Set $K^+$ consists of some of the features in $\{1,\ldots,n\}\setminus K$ whose
associated $v$ variables will be considered as binary in the next iteration, i.e. features of $K^+$ will be added to $K$. Similarly, $K^-$
consists of features in $K$ that will not be considered as binary in the next iteration. In addition and if possible, we will update the UB in the main step.

A general outline of the exact procedure is shown in  Algorithm \ref{alg2}. Since the set $K$ can be modified using different rules to improve the lower bounds, we have provided three different update variants of this procedure in Algorithms  
 \ref{alg2}.\ref{alg3}, \ref{alg2}.\ref{alg4} and \ref{alg2}.\ref{alg6}. In Variant I (Algorithm \ref{alg2}.\ref{alg3}), the set $K$ is updated by adding the features using vector $r$, sorted in non-decreasing order as described in the previous section. 
 This ordered sequence of $v$-variables is based on the idea that the features with the biggest reduced costs of variables $w_j^+$ and $w_j^-$ are less likely to be different from $0$ in the optimal solution of FS-SVM and those features with a positive value in the LP are the most likely to be different from $0$.
 
 \scalebox{0.85}{
 	\begin{algorithm}[H]
 		\KwData{Training data composed by a set of $m$ elements with $n$ features.}
 		\KwResult{Optimal objective value or accurate upper and lower bounds. }
 		\tcc{ Initial Step}
 		Run strategies I and II to tighten $u$ and $l$ and to obtain initial LB and UB.\\
 		 Run the Kernel Search to obtain UB$_{new}$.\\ 
 		\If {UB$_{new}<$UB}{ UB$:=$UB$_{new}$ and run strategies I and II again.}
 		\tcc{Main Step}
 		Let $K$ be the final set obtained using the Kernel Search.
 		Take $it:=0$.\\
 		\While{$\frac{UB-LB}{UB}\cdot 100\le 0.01$}{
 			Solve SR-FS-SVM($K$). Let $(\overline{v},\overline{w}^{a},\overline{w}^{b},\overline{b},\overline{\xi})$ be its solution, $z_{LB}$ its objective value.\\
 			\If {running time of FS-SVM($K$) $>1800$s}{break}
 			Let $\tilde{v}_j$ be the optimal values of variables $v_j$ with $j\in K$ of SR-FS-SVM($K$). Solve the FS-SVM fixing $v_j=0\; \forall j\in \{1,\ldots,n\}\setminus K$ and $v_j=\tilde{v}_j\;\forall j\in K$. Let $z_{UB}$ be its objective value.\\
 			\If{$z_{UB}<UB$}{
 				$UB:=z_{UB}$.
 			}
 			\If{$z_{LB}> LB$}{
 				{ $LB:=z_{LB}$.}
 			}
 			\tcc{ Update Step}
 			Build $K^-$, composed by the features in $K$ that will be relaxed in the next iteration. \\
 			Build $K^+$, composed by the features in $\{1,\ldots,n\}\setminus K$ that will be added to $K$ in the next iteration.\\
 			Update $K:=K\cup K^+\setminus K^-$.\\
 			$it :=it+1$\\
 		}
 		\caption{Exact Procedure}
 		\label{alg2}
 	\end{algorithm}
 }
 
\vspace{0.25cm}
 In contrast, 
in Variant II (Algorithm \ref{alg2}.\ref{alg4}) the set $K$ is increased in each iteration by adding features $j\in \{1,\ldots,n\}\setminus K$ that take a value bigger than $0$ in the SR-FS-SVM($K$). 
 Lastly, in Variant III (Algorithm \ref{alg2}.\ref{alg6}) the set $K$ is modified based on the reduced costs of the resulting linear programming problem after fixing  the binary variables of SR-FS-SVM($K$) to their optimal values. 
We thus obtain the reduced costs of variables $w^+$ and $w^-$, and a vector similar to $r$ is created. In this case, we have denoted it as $\overline{r}$ and it is defined as:
\begin{equation}
\overline{r}_j:=\left\{\begin{array}{c c}
-(\tilde{\omega}_j^++\tilde{\omega}_j^-)& \textnormal{if }\tilde{\omega}_j^++\tilde{\omega}_j^->0,\\
\min\{ \overline{r}^+_j,\overline{r}^-_j\}&\textnormal{otherwise}.
\end{array}
\right.
\label{cost_red2}
\end{equation}
where $\tilde{\omega}^+$, $\tilde{\omega}^-$, $\overline{r}^+$ and $\overline{r}^-$ are the solutions and reduced costs of the problem described above. $\overline{r}$ is sorted in non-decreasing order and $K$ is updated by adding the features in this order. A set of $\bar{n}<n$ features are added in each iteration. In particular we take $\bar{n}=20$ since it provides good results. Additionally, variables that are null in two consecutive iterations are relaxed.

Using these variants, we explore different forms to improve the lower bounds and to update the initial set of binary variables. The 
various performances of the described procedures are analyzed in Section \ref{computational}. 

\scalebox{0.85}{
\begin{megaalgorithm}[H]
\tcc{Modified Main Step: Variant I}
Sort the features $j\in \overline{K} :=\{1,\ldots,n\}\setminus K$ according to vector  $r$ defined in \eqref{cost_red}.\\
Divide $\overline{K}$ in a sequence of subsets {($\overline{K_1},\overline{K_2},\ldots$) of a certain size $S$}, considering the order given by $r$.\\
\tcc{ Update Step: Variant I}
 Build $K^+ :=K_{it}$. Update $K:=K\cup K^+$.\\
 \caption{Update Variant I.}
 \label{alg3}
\end{megaalgorithm}
}

\scalebox{0.85}{
\begin{megaalgorithm}[H]
\tcc{ Update Step: Variant II}
$K^+:=\{j\in \overline{K}| \overline{v}_j>0\}$. Update $K:=K\cup K^+$ and $\overline{K}:=\{1,\ldots,n\}\setminus K$.\\
 \caption{Update Variant II}
 \label{alg4}
\end{megaalgorithm}
}

\scalebox{0.85}{
\begin{megaalgorithm}[H]
\tcc{ Update Step: Variant III}
Build $K^- :=\{j\in K| j$ has not been selected in the solution of the last two iterations$\}$.\\
Solve the LP program resulting of fixing the binary variables of SR-FS-SVM($K$) to its optimal values.\\
 Sort {$\overline{K}:=\{1,\ldots,n\}\setminus K$} in non-decreasing according to the values of $\overline{r}$ vector in \eqref{cost_red2}.\\
Construct the set $K^+$ selecting the first  $\bar{n}<n$ features of the ordered set $\overline{K}$. In particular, we take $\bar{n}:=20$.\\
Update $K:=K\cup K^+\setminus K^-$.\\
 \caption{Update Variant III}
 \label{alg6}
\end{megaalgorithm}
} 

\section{Computational Results}\label{computational}
In this section, we present the results provided by several computational experiments. In particular: i) how the use of Strategies I and II for fixing  the upper and lower bounds of $w$ variables can reduce the computing times for our model; ii) the efficiency of the heuristic approach
(Kernel Search) proposed in this paper; and lastly, iii) the study of the results provided by the different variants of the exact solution approaches.

It should be noted that the various computational experiments were performed using CPLEX 12.6.3 on an Intel(R) Core(TM) i7-4790K CPU 32 GB RAM computer. We should also remark that the CutPass, CutsFactor, EachCutLim, FracCuts, PreInd, RinHeur, EpInt and EpRHS parameters were modified in order to give a clean comparison of the relative performance of the formulations, i.e. by using these parameters we tried to avoid using the CPLEX internal heuristics since they can have a different influence on the previously described solution variants. The computational experiments were carried out on sixteen different datasets. Eight of them can be found in 
the UCI repository (\cite{AsuNew07}), (see Table \ref{small}), where $m$ is the number of elements, $n$ is the number of features and the last column shows the percentage of elements in each class. As can be observed, they contain a small number of features. The other  eight datasets used in the experiments have a larger number of features (see Table \ref{big}).  The Lepiota, Arrythmia, Madelon and MFeat datasets are also in UCI repository. A further description of the remaining datasets in Table \ref{big} can be found in  \cite{Alon}, \cite{car10},
	\cite{Guyon02}, \cite{Maldo14}, \cite{Golub99}, \cite{Shipp} and \cite{Not}.
\begin{table}[!htb]
	{
	\begin{minipage}{.45\linewidth}
		\centering
		\scalebox{0.7}{
				\begin{tabular}{l|ccc}
				\hline
				\multicolumn{4}{c}{Small number of features} \\
				\hline
				Name  & m     & n     &  Class(\%)  \\
				\hline
				BUPA  & 345   & 6     & 42/58 \\
				PIMA  & 768   & 8     & 65/35  \\
				Cleveland & 297   & 13    & 42/58  \\
				Housing & 506   & 13    & 51/49 \\
				Australian & 690   & 14    & 44/56  \\
				GC    & 1000  & 24    & 30/70  \\
				WBC   & 569   & 30    & 37/63  \\
				Ionosphere & 351   & 33    & 64/36  \\
				\hline
			\end{tabular}%
		}
	\subcaption{\label{small}}
	\end{minipage}%
	\begin{minipage}{.55\linewidth}
		\centering
		\scalebox{0.7}{
				\begin{tabular}{l|rrl|l|rrl}
				\hline
				\multicolumn{8}{c}{Big number of features} \\
				\hline
				\multicolumn{4}{c|}{Small sample size} & \multicolumn{4}{c}{Big sample size} \\
				\hline
				Name  & \multicolumn{1}{l}{m} & \multicolumn{1}{l}{n} & Class(\%) & Name  & \multicolumn{1}{l}{m} & \multicolumn{1}{l}{n} & Class(\%) \\
				\hline
				Colon & 62    & 2000  & 35/65 & Lepiota & 1824  & 109   & 52/48 \\
				Leukemia & 72    & 5327  & 47/53 & Arrythmia & 420   & 258   & 57/43 \\
				DLBCL & 77    & 7129  & 75/25 & Madelon & 2000  & 500   & 50/50 \\
				Carcinoma & 36    & 7457  & 53/47 & Mfeat & 2000  & 649   & 10/90 \\
				\hline
			\end{tabular}%
		}
	\subcaption{\label{big}}
	\end{minipage} 
	\caption{Datasets description. \label{datset}}
}
\end{table}

Since the resolution times for the datasets in Table \ref{small} and Lepiota dataset are really good when solving the formulation with CPLEX (a few seconds), in this section we focus our attention of developing alternative solution strategies to the  instances with the largest number of features (different from Lepiota). In Subsection  \ref{smallbig}, we will analyze the instances with a big number of features and a small sample size (Colon, Leukemia, DLBCL, Carcinoma). Lastly, we will apply the best obtained techniques to the instances with a big sample size in Subsection \ref{bigsample}. 

\subsection{Analysis of datasets with small sample size and big number of features}\label{smallbig} In this subsection we will analyze how the use of strategies I and II can affect Colon, Leukemia, DLBCL and Carcinoma datasets. The heuristic and exact procedures applied to these datasets are also studied in this subsection.

\subsubsection{FS-SVM with Strategies I and II \label{Kernel_stra}}  
Section \ref{strategies} described two strategies for obtaining tightened bounds on parameters $u$ and $l$.
Table \ref{tab_prep} reports the computational results of the proposed formulation, both
with and without Strategies I and II. We used small values of $B\in\{10,20,30\}$, $C\in \{2^{-7},2^{-6},2^{-5},2^{-4},2^{-3},2^{-1},1,$ 
$2,2^{2},2^{4},2^{5},2^{6},2^7\}$
and a time limit of two hours (instances exceeding the time limit have been highlighted with their times underlined).  However, since the running times are very short for small values of $C$ and the variables $w$ result to be null in most cases, we have only reported the results for $C\in \{1,2,2^{2},2^{4},2^{5},2^{6},2^7\}$. In this Table, the column labelled ``FS-SVM'' shows the gaps and running times of the proposed model.  The second group of columns for each dataset, titled ``St.+FS-SVM'', shows the results 
associated with the model after strategies I and II have been applied for obtaining tightened lower/upper bounds of $w$-variables. The termination gap (\%) is shown in the ``Gap'' column, whilst the ``$t_{st}$'' column  gives the time required for the two strategies, $t_{solv}$ is the running time  for solving  the formulation once the parameters defining the bounds of $w$ have been fixed and $t_{total}$ is the 
overall process time. Lastly, column $\Delta B=\frac{1}{n}\sum_{j=1}^n (u_j-l_j)$ shows the average difference between the upper and lower bounds after the use of both strategies. Note that initially $u_j=-l_j$ takes a large enough amount. Generally, it can be seen that the use of the strategies provides an average difference between both bounds of less than $5.57$ units. Therefore, Strategies I and II provide tightened bounds.

For the Colon dataset, FS-SVM cannot be solved for  $B=20$, $C\ge 2$ and $B=10$, $C\ge 2^0$ within the time limit. However, 
if Strategies I and II are employed before solving the model, FS-SVM can be solved in less than $30$ minutes for
all cases when $B=20$. For $B=10$ and $C\ge 2$, the model cannot be solved in less than two hours, even if the strategies are performed, but the gaps at termination are smaller.

The second group of columns of Table \ref{tab_prep} shows the results for the Leukemia dataset. In this instance, the model with the strategies solves
the same cases as the model without the strategies. However, most of the cases that cannot be solved in less than two hours present 
smaller gaps if the strategies are used.

Table \ref{tab_prep}, also details the results for the DLBCL dataset. For $B=20$ and $C\ge 2^3$, this instance was not solved within the time limit, but when using the strategies it can be solved in less than approximately $19$ minutes. However, the model in which $B=10$ and $C\ge 1$ cannot be solved, even if $u$ and $l$ are tightened with the strategies, although the gaps are once again better than when the strategies are not utilized.

The last group of columns shows the results of the Carcinoma dataset. In this case, the model cannot be solved in less than two hours for $B=10$. However, if we use both strategies, the solution times improve and only two cases (with parameters $C=2^1,2^2$) remained unsolved after the time limit.

Regarding the reported results, we can conclude that the use of Strategies
I and II leads to a reduction in running times in most cases. Furthermore,
although the model cannot be solved for certain parameters values (even with the use of the strategies), the termination
gaps are better if Strategies I and II are employed. 

\begin{sidewaystable}[htbp]
	\centering
	\scalebox{0.6}{
		\begin{tabular}{r@{/}c||rr|rrrrr||rr|rrrrr||rr|rrrrr||rr|rrrrr}
			\hline
			\multicolumn{2}{c||}{}       & \multicolumn{7}{c||}{Colon m=62 n=2000}                 & \multicolumn{7}{c||}{Leukemia m=72 n=5327}              & \multicolumn{7}{c||}{DLBCL m=77 n=7129}                 & \multicolumn{7}{c}{Carcinoma m=36 n=7457} \\
			\hline
			\multicolumn{2}{c||}{\multirow{2}[4]{*}{B/C}} & \multicolumn{2}{c|}{FS-SVM} & \multicolumn{5}{c||}{{ St. + FS-SVM}} & \multicolumn{2}{c|}{FS-SVM} & \multicolumn{5}{c||}{{ St. + FS-SVM}} & \multicolumn{2}{c|}{FS-SVM} & \multicolumn{5}{c||}{{ St. + FS-SVM}} & \multicolumn{2}{c|}{FS-SVM} & \multicolumn{5}{c}{{ St. + FS-SVM}} \\
			\multicolumn{2}{c||}{}   & Gap   & Time  & Gap   & $t_{st}$ & $t_{solv}$ & $t_{total}$ & $\Delta B$ & Gap   & Time  & Gap   & $T_{st}$ & $t_{solv}$ & $t_{total}$ & $\Delta B$ & Gap   & Time  & Gap   & $t_{st}$ & $t_{solv}$ & $t_{total}$ & $\Delta B$ & Gap   & Time  & Gap   & $t_{st}$ & $t_{solv}$ & $t_{total}$ & $\Delta B$ \\
			\hline
			30    & $2^0$     & \textbf{0.0} & \textbf{1.44} & 0.0   & 34.66 & 2.35  & 37.00 & { 0.02}  & \textbf{0.0} & \textbf{726.84} & 0.0   & 447.50 & 292.67 & 740.17 & { 0.05}  & \textbf{0.0} & \textbf{16.32} & 0.0   & 646.61 & 21.81 & 668.41 & { 0.02}  & \textbf{0.0} & \textbf{0.56} & \textbf{0.0} & 307.27 & 4.32  & 311.60 & { 0.01} \\
			30    & $2^1$      & \textbf{0.0} & \textbf{5.30} & 0.0   & 32.51 & 4.42  & 36.93 & { 0.05}  & 0.0   & 792.12 & \textbf{0.0} & 456.01 & 281.33 & \textbf{737.34} &{ 0.05}   & \textbf{0.0} & \textbf{13.57} & 0.0   & 649.36 & 21.48 & 670.84 & { 0.02}  & \textbf{0.0} & \textbf{0.59} & \textbf{0.0} & 305.81 & 4.31  & 310.11 & { 0.01} \\
			30    & $2^2$    & \textbf{0.0} & \textbf{3.82} & 0.0   & 32.32 & 4.45  & 36.77 & { 0.05}  & \textbf{0.0} & \textbf{814.42} & 0.0   & 457.07 & 598.00 & 1055.06 & { 0.05}   & \textbf{0.0} & \textbf{16.93} & 0.0   & 637.28 & 21.26 & 658.54 &{ 0.02}  & \textbf{0.0} & \textbf{0.61} & \textbf{0.0} & 305.12 & 4.35  & 309.47 & { 0.01} \\
			30    & $2^3$ & \textbf{0.0} & \textbf{3.87} & 0.0   & 32.25 & 4.71  & 36.96 & { 0.05}  & \textbf{0.0} & \textbf{710.88} & 0.0   & 453.51 & 885.04 & 1338.55 & { 0.05}   & \textbf{0.0} & \textbf{12.31} & 0.0   & 644.91 & 21.22 & 666.12 & { 0.02}  & \textbf{0.0} & \textbf{0.58} & \textbf{0.0} & 305.73 & 4.32  & 310.05 & { 0.01} \\
			30    & $2^4$    & \textbf{0.0} & \textbf{4.10} & 0.0   & 32.00 & 5.27  & 37.27 & { 0.05}  & \textbf{0.0} & \textbf{791.81} & 0.0   & 455.26 & 425.10 & 880.36 & { 0.05}   & \textbf{0.0} & \textbf{15.41} & 0.0   & 642.77 & 22.37 & 665.14 & { 0.02}  & \textbf{0.0} & \textbf{0.60} & \textbf{0.0} & 305.54 & 4.32  & 309.86 & { 0.01} \\
			30    & $2^5$    & \textbf{0.0} & \textbf{5.26} & 0.0   & 32.41 & 5.10  & 37.51 & { 0.05} & \textbf{0.0} & \textbf{621.57} & 0.0   & 440.89 & 199.59 & 640.47 & { 0.05}   & \textbf{0.0} & \textbf{16.17} & 0.0   & 640.96 & 24.96 & 665.92 & { 0.02}  & \textbf{0.0} & \textbf{0.59} & \textbf{0.0} & 308.58 & 4.34  & 312.92 & { 0.01} \\
			30    & $2^6$       & \textbf{0.0} & \textbf{5.27} & 0.0   & 32.12 & 5.59  & 37.70 & { 0.05}  & \textbf{0.0} & \textbf{569.28} & 0.0   & 440.59 & 173.85 & 614.44 & { 0.05}   & \textbf{0.0} & \textbf{15.70} & 0.0   & 644.73 & 25.91 & 670.65 & { 0.02}  & \textbf{0.0} & \textbf{0.59} & \textbf{0.0} & 306.58 & 4.32  & 310.90 &{ 0.01} \\
			30    & $2^7$      & \textbf{0.0} & \textbf{5.65} & 0.0   & 32.37 & 5.30  & 37.67 & { 0.05}  & \textbf{0.0} & \textbf{527.17} & 0.0   & 440.78 & 285.78 & 726.56 & { 0.05}   & \textbf{0.0} & \textbf{16.12} & 0.0   & 643.00 & 22.73 & 665.73 & { 0.02}  & \textbf{0.0} & \textbf{0.60} & \textbf{0.0} & 305.90 & 4.33  & 310.23 & { 0.01} \\
			\hline
			20    & $2^0$       & 0.0   & 579.94 & \textbf{0.0} & 31.91 & 93.35 & \textbf{125.26} & { 0.40}  & 1.1   & \underline{7201.03} & \textbf{0.5} & 362.83 & \underline{7205.94} & \underline{7567.88} & { 0.25}  & 0.0   & { 7200.40} & \textbf{0.0} & 536.39 & 511.59 & \textbf{1047.98} &{ 0.25} & \textbf{0.0} & 1457.92 & \textbf{0.0} & 310.88 & 181.54 & \textbf{492.42} & { 0.04} \\
			20    & $2^1$       & 0.3   & \underline{7200.19} & \textbf{0.0} & 35.65 & 1911.11 & \textbf{1946.76} & { 0.79}  & 1.0   & \underline{7200.90} & \textbf{0.5} & 388.25 & \underline{7206.56} & \underline{7593.99} & { 0.25}   & 0.0   & 5759.65 & \textbf{0.0} & 536.98 & 582.46 & \textbf{1119.44} & { 0.25}   & \textbf{0.0} & 1478.06 & \textbf{0.0} & 310.94 & 159.89 & \textbf{470.82} & { 0.04}  \\
			20    & $2^2$       & 0.3   & \underline{7200.10} & \textbf{0.0} & 35.30 & 1634.42 & \textbf{1669.72} &{ 0.79}  & 2.0   & 7201.08 & \textbf{0.6} & 357.72 & \underline{7206.43} & \underline{7563.36} & { 0.25}   & 0.0   & 7168.70 & \textbf{0.0} & 535.46 & 515.52 & \textbf{1050.97} & { 0.25}   & \textbf{0.0} & 1551.47 & \textbf{0.0} & 309.65 & 204.00 & \textbf{513.65} & { 0.04}  \\
			20    & $2^3$       & 0.3   & \underline{7200.10} & \textbf{0.0} & 34.34 & 1652.39 & \textbf{1686.73} & { 0.79}  & 2.0   & \underline{7200.97} & \textbf{0.6} & 362.14 & \underline{7207.29} & \underline{7568.63} &{ 0.25}   & 0.1   & \underline{7200.38} & \textbf{0.0} & 534.11 & 472.57 & \textbf{1006.68} & { 0.25}   & \textbf{0.0} & 1478.82 & \textbf{0.0} & 309.27 & 205.08 & \textbf{514.35} & { 0.04}  \\
			20    & $2^4$      & 0.3   & \underline{7200.10} & \textbf{0.0} & 34.80 & 1549.91 & \textbf{1584.71} & { 0.79}  & 1.9   & 7200.99 & \textbf{0.5} & 385.96 & \underline{7206.79} & \underline{7591.77} & { 0.25}   & 0.0   & \underline{7200.40} & \textbf{0.0} & 535.53 & 513.14 & \textbf{1048.67} & { 0.25}   & \textbf{0.0} & 2016.24 & \textbf{0.0} & 308.58 & 210.23 & \textbf{518.82} & { 0.04} \\
			20    & $2^5$      & 0.3   & \underline{7200.1}0 & \textbf{0.0} & 34.32 & 1380.88 & \textbf{1415.20} & { 0.79} & 1.9   & \underline{7200.82} & \textbf{0.1} & 363.40 & \underline{7205.80} & \underline{7568.42} & { 0.25}   & 0.1   & \underline{7200.38} & \textbf{0.0} & 539.71 & 366.16 & \textbf{905.88} & { 0.25}   & \textbf{0.0} & 2680.07 & \textbf{0.0} & 308.94 & 183.02 & \textbf{491.96} & { 0.04}  \\
			20    & $2^6$    & 0.4   & \underline{7200.10} & \textbf{0.0} & 34.23 & 1098.26 & \textbf{1132.48} &{ 0.79}  & 1.8   & \underline{7201.05} & \textbf{0.2} & 383.67 & \underline{7206.14} & \underline{7588.98} &{ 0.25}   & 0.3   & \underline{7200.42} & \textbf{0.0} & 541.42 & 335.07 & \textbf{876.49} & { 0.25}   & \textbf{0.0} & 3283.47 & \textbf{0.0} & 311.44 & 108.29 & \textbf{419.73} & { 0.04}  \\
			20    & $2^7$      & 0.5   & \underline{7200.11} & \textbf{0.0} & 34.76 & 1009.93 & \textbf{1044.69} & { 0.79}  & 1.8   & \underline{7200.86} & \textbf{0.3} & 375.56 & \underline{7205.54} & \underline{7580.24} & { 0.25}  & 0.3   & \underline{7200.39} & \textbf{0.0} & 540.24 & 373.70 & \textbf{913.94} & { 0.25}   & \textbf{0.0} & 3231.34 & \textbf{0.0} & 309.85 & 108.47 & \textbf{418.32} & { 0.04}  \\
			\hline
			10    & $2^0$       & 2.5   & \underline{7203.00} & \textbf{0.0} & 40.03 & 4598.49 & \textbf{4638.52} & { 2.28}  & \textbf{15.1} & \underline{7202.81} & 15.9  & 764.48 & \underline{7204.22} & \underline{7968.70} & { 1.71}  & 11.8  & \underline{7203.61} & \textbf{8.1} & 1013.27 & \underline{7209.40} & \underline{8222.67} & { 1.62} & 2.4   & \underline{7203.84} & \textbf{0.0} & 412.96 & 6789.65 & \textbf{7202.61} & { 0.37} \\
			10    & $2^1$      & 12.4  & \underline{7204.81} & \textbf{11.2} & 46.16 & \underline{7203.86} & \underline{7250.02} & { 5.19}   & \textbf{17.5} & \underline{7200.71} & 19.4  & 767.57 & \underline{7204.40} & \underline{7971.97} & { 1.71}  & 11.1  & \underline{7200.94} & \textbf{7.2} & 1015.05 & \underline{7207.12} & \underline{8222.17} &{ 1.63}  & 2.3   & \underline{7204.3} & \textbf{0.5} & 413.39 & \underline{7204.63} & \underline{7618.01} &{ 0.37}  \\
			10    &$2^2$      & \textbf{10.5} & \underline{7200.94} & 13.6  & 51.20 & \underline{7206.82} & \underline{7258.02} & { 5.56}  & \textbf{16.9} & \underline{7200.29} & 17.5  & 783.36 & \underline{7204.02} & \underline{7987.37} & { 1.71}  & 13.5  & \underline{7200.40} & \textbf{9.5} & 1004.74 & \underline{7209.84} & \underline{8214.58} & { 1.63}   & 2.3   & \underline{7203.33} & \textbf{0.5} & 415.55 & \underline{7205.41} & \underline{7620.96} & { 0.37}  \\
			10    & $2^3$     & 15.2  & \underline{7201.77} & \textbf{12.3} & 46.35 & \underline{7206.01} & \underline{7252.36} & { 5.56}  & 16.8  & \underline{7203.84} & \textbf{14.9} & 779.71 & \underline{7204.02} & \underline{7983.72} & { 1.71} & 10.8  & \underline{7202.02} & \textbf{10.6} & 1008.17 & \underline{7206.94} & \underline{8215.11} & { 1.63}   & 2.3   & \underline{7203.56} & \textbf{0.0} & 419.50 & 7078.98 & \textbf{7498.48} & { 0.37} \\
			10    & $2^4$       & 16.9  & \underline{7203.00} & \textbf{13.5} & 47.03 & \underline{7206.16} & \underline{7253.19} & { 5.56}  & \textbf{12.9} & \underline{7191.54} & 15.7  & 781.86 & \underline{7207.05} & \underline{7988.91} & { 1.71}  & 13.5  & \underline{7200.41} & \textbf{9.2} & 1006.23 & \underline{7210.54} & \underline{8216.77} & { 1.63}   & 2.3   & \underline{7202.94} & \textbf{0.0} & 420.87 & 4460.00 & \textbf{4880.87} & { 0.37}  \\
			10    & $2^5$    & 17.1  & \underline{7206.07} & \textbf{13.0} & 49.17 & \underline{7206.29} & \underline{7255.46} &{ 5.56}  & 18.8  & \underline{7200.29} & \textbf{15.2} & 784.12 & \underline{7204.11} & \underline{7988.23} & { 1.71}  & 13.5  & \underline{7200.40} & \textbf{10.4} & 1018.63 & \underline{7210.67} & \underline{8229.30} & { 1.63}   & 2.9   & \underline{7206.39} & \textbf{0.0} & 415.71 & 2551.15 & \textbf{2966.86} & { 0.37} \\
			10    & $2^6$    & 16.6  & \underline{7204.87} & \textbf{12.3} & 46.80 & \underline{7205.97} & \underline{7252.77} & { 5.56} & 18.7  & \underline{7203.20} & \textbf{17.9} & 754.17 & \underline{7206.49} & \underline{7960.66} & { 1.71}  & 13.5  & \underline{7200.40} & \textbf{10.1} & 998.42 & \underline{7197.19} & \underline{8195.61} & { 1.63}   & 2.7   & \underline{7205.12} & \textbf{0.0} & 418.61 & 2744.28 & \textbf{3162.89} & { 0.37}  \\
			10    & $2^7$     & 17.0  & \underline{7201.27} & \textbf{10.8} & 46.39 & \underline{7206.23} & \underline{7252.62} & { 5.56}  & \textbf{15.2} & \underline{7203.81} & 16.5  & 775.13 & \underline{7207.06} & \underline{7982.19} & { 1.71}  & 13.5  & \underline{7200.38} & \textbf{9.9} & 1004.78 & \underline{7210.10} & 8214.88 & { 1.63}   & 2.6   & \underline{7206.23} & \textbf{0.0} & 422.03 & 2506.86 & \textbf{2928.89} & { 0.37}  \\
			\hline
		\end{tabular}%
	}
	\caption{Time performance comparison of FS-SVM with and without strategies in datasets with a big number of features. \label{tab_prep}}
\end{sidewaystable}%

\subsubsection{Heuristic procedure}\label{heur}

Table \ref{tab_Kernel} shows the results of applying the Kernel Search to the proposed model for the four datasets: Colon, Leukemia, DLBCL and Carcinoma. The ``Gap'' column details the gap (\%) between the solution provided by the KS and the optimal solution obtained by the FS-SVM formulation (note that for $B=20,\,30$ the formulation was run until the optimal solution was found and for $B=10$ we ran the formulation until a gap of less than 10\% was achieved). For cases where the optimal solutions were not found, this column shows, both the gap between the best lower bound and the KS solution, as well as the gap between the best feasible solution and the KS solution. These two gaps are separated by ``/''.
Furthermore, the columns labelled ``$t_{KS}$'' represent the running time of the KS. Lastly, in ``$t_{Best}$'' columns we can observe the best time for solving the problems.  In those instances that were not solved within two hours, we reported the total time of the solution method (FS-SVM, St+FS-SVM or St+FS-SVM$^{*}$) that provides the best gap (recall that St+FS-SVM$^{*}$ corresponds to the results of FS-SVM using Strategies I and II, establishing a time limit of two hours and using the Kernel Search solution as the initial solution for the model, see Table \ref{kernel_col}).
Shown in superscript is
the gap obtained after the time limit was reached.

In general, Table \ref{tab_Kernel} demonstrates that the Kernel Search provides, in much less time, the optimal solution for all those instances that we were able to solve exactly in less than two hours. In the remaining cases, we provided better upper bounds than the exact approaches, again in much less time.

Due to these good results in terms of gaps, this heuristic may also be useful in reducing the computational times of the formulation. Table \ref{kernel_col} shows a comparison of the gaps and times between the model that uses the KS solution as an initial feasible solution of the problem (column St.+FS-SVM$^*$) and the results reported in Table \ref{tab_prep}. In general, we can observe a reduction of computational times and an improvement of the gaps in the  problems that could not be solved within $7200$ seconds. In each case, the best times and gaps are shown in bold. Additionally, figures \ref{fig:COL1}, \ref{fig:DL1}, \ref{fig:car1} and \ref{fig:car2} show the best improvements that are achieved using the KS solution.

\begin{table}[htbp]
	\centering
	\scalebox{0.7}{
		\begin{tabular}{r@{/}c|rrr|rrr|rrr|rrr}
			\hline
			\multicolumn{2}{c|}{}      & \multicolumn{3}{c|}{Colon} & \multicolumn{3}{c|}{Leukemia} & \multicolumn{3}{c|}{DLBCL} & \multicolumn{3}{c}{Carcinoma} \\
			\hline
			B     & C     & Gap   & $t_{KS}$ & $t_{Best}$ & Gap   & $t_{KS}$ & $t_{Best}$ & Gap   & $t_{KS}$ & $t_{Best}$ & Gap   & $t_{KS}$ & $t_{Best}$ \\
			\hline
			\hline
			30    & $1$ & 0.0   & 0.48  & 1.4   & 0.0   & 9.22  & 661.64 & 0.0   & 2.37  & 16.32 & 0.0   & 1.39  & 0.56 \\
			30    & $2$ & 0.0   & 0.60  & 5.30  & 0.0   & 8.53  & 703.14 & 0.0   & 2.26  & 13.57 & 0.0   & 1.35  & 0.59 \\
			30    & $2^{2}$ & 0.0   & 0.53  & 3.82  & 0.0   & 6.02  & 710.35 & 0.0   & 2.29  & 16.93 & 0.0   & 1.34  & 0.61 \\
			30    & $2^{3}$ & 0.0   & 0.62  & 3.87  & 0.0   & 6.23  & 710.88 & 0.0   & 2.23  & 12.31 & 0.0   & 1.42  & 0.58 \\
			30    & $2^{4}$ & 0.0   & 0.60  & 4.10  & 0.0   & 8.14  & 791.81 & 0.0   & 2.54  & 15.41 & 0.0   & 1.33  & 0.60 \\
			30    & $2^{5}$ & 0.0   & 0.59  & 5.26  & 0.0   & 9.16  & 586.30 & 0.0   & 2.36  & 16.17 & 0.0   & 1.37  & 0.59 \\
			30    & $2^{6}$ & 0.0   & 0.63  & 5.27  & 0.0   & 13.25 & 569.28 & 0.0   & 2.40  & 15.70 & 0.0   & 1.41  & 0.59 \\
			30    & $2^{7}$ & 0.0   & 0.63  & 5.65  & 0.0   & 10.53 & 527.17 & 0.0   & 2.56  & 16.12 & 0.0   & 1.30  & 0.60 \\
			\hline
			20    & $1$ & 0.0   & 5.77  & 125.26 & 0.1   & 426.77 & 7567.88$^{(0.5)}$ & 0.0   & 9.13  & 911.76 & 0.0   & 3.31  & 492.42 \\
			20    & $2$ & 0.0   & 20.87 & 1388.57 & 0.1   & 524.43 & 7593.99$^{(0.5)}$ & 0.0   & 11.98 & 908.02 & 0.0   & 4.66  & 470.82 \\
			20    & $2^{2}$ & 0.0   & 12.90 & 1564.57 & 0.1   & 475.04 & 7563.36$^{(0.6)}$ & 0.0   & 9.86  & 1004.88 & 0.0   & 3.67  & 501.49 \\
			20    & $2^{3}$ & 0.0   & 10.85 & 1528.63 & 0.1   & 473.22 & 7605.20$^{(0.5)}$ & 0.0   & 10.72 & 1006.68 & 0.0   & 3.85  & 514.35 \\
			20    & $2^{4}$ & 0.0   & 12.96 & 1580.48 & 0.1   & 372.09 & 7591.77$^{(0.5)}$ & 0.0   & 10.75 & 869.59 & 0.0   & 4.33  & 493.65 \\
			20    & $2^{5}$ & 0.0   & 11.82 & 1415.20 & 0.1   & 243.00 & 7568.42$^{(0.1)}$ & 0.0   & 10.86 & 905.16 & 0.0   & 5.66  & 476.89 \\
			20    & $2^{6}$ & 0.0   & 11.15 & 1087.20 & 0.1   & 477.97 & 7588.98$^{(0.2)}$ & 0.0   & 11.81 & 801.61 & 0.0   & 4.62  & 419.73 \\
			20    & $2^{7}$ & 0.0   & 20.51 & 712.30 & 0.1   & 419.63 & 7620.99$^{(0.2)}$ & 0.0   & 12.53 & 903.79 & 0.0   & 8.39  & 418.32 \\
			\hline
			10    & $1$ & 0.0   & 18.12 & 4447.71 & 0.1/9.7 & 1082.17 & 8980.67$^{(11.4)}$ & 0.0/5.5 & 210.48 & 8392.61$^{(6.8)}$ & 0.0   & 32.77 & 6931.78 \\
			10    & $2$ & 0.0/6.3 & 317.79 & 7396.14$^{(8.5)}$ & 0.0/9.8 & 1307.56 & 9079.71$^{(11.4)}$ & 0.0/6.5 & 259.54 & 8390.48$^{(6.6)}$ & 0.0   & 36.10 & 7331.22 \\
			10    & $2^{2}$ & 0.0/7.5 & 453.18 & 7442.52$^{(9.3)}$ & 0.9/10.5 & 1683.48 & 9212.35$^{(11.3)}$ & 0.0/6.6 & 227.28 & 8382.97$^{(6.6)}$ & 0.0   & 40.94 & 6994.41 \\
			10    & $2^{3}$ & 0.0/7.5 & 581.85 & 7470.23$^{(9.1)}$  & 0.9/10.6 & 1581.34 & 9214.53$^{(11.3)}$ & 0.0/6.6 & 200.51 & 8404.06$^{(6.6)}$ & 0.0   & 39.27 & 6625.83 \\
			10    & $2^{4}$ & 0.0/7.8 & 489.36 & 7481.98$^{(9.2)}$  & 0.0/9.4 & 1302.88 & 9236.43$^{(11.1)}$ & 0.0/6.5 & 213.77 & 8410.13$^{(6.5)}$ & 0.0   & 30.85 & 4880.87 \\
			10    & $2^{5}$ & 0.0/8.0 & 508.48 & 7434.46$^{(9.1)}$  & 0.0/9.2 & 1376.44 & 9243.22$^{(10.9)}$ & 0.0/6.2 & 205.48 & 8405.00$^{(6.2)}$ & 0.0   & 23.70 & 2966.86 \\
			10    & $2^{6}$ & 0.0/7.3 & 532.52 & 7485.48$^{(9.1)}$  & 0.0/9.9 & 1372.27 & 9284.75$^{(11.0)}$ & 0.0/6.1 & 218.51 & 8406.76$^{(6.1)}$ & 0.0   & 29.38 & 3162.89 \\
			10    & $2^{7}$ & 0.0/7.4 & 499.11 & 7437.58$^{(9.1)}$  & 0.0/9.4 & 1463.05 & 9261.78$^{(10.9)}$ & 0.0/6.2 & 234.97 & 8407.38$^{(6.2)}$ & 0.0   & 26.53 & 2928.89 \\
			\hline
		\end{tabular}%
	}
	\caption{Kernel Search gaps and running times for large datasets.}
	\label{tab_Kernel}%
\end{table}%

\begin{sidewaystable}[htbp]
	\centering
	\scalebox{0.7}{
		\begin{tabular}{r@{/}c||rr|rr|rr||rr|rr|rr||rr|rr|rr||rr|rr|rr}
			\hline
			\multicolumn{2}{c||}{}     & \multicolumn{6}{c||}{Colon m=62,n=2000}         & \multicolumn{6}{c||}{Leukemia m=72, n=5327}     & \multicolumn{6}{c||}{DLBCL m=77, n=7129}        & \multicolumn{6}{c}{Carcinoma m=36, n=7457} \\
			\hline
			\multicolumn{2}{c||}{\multirow{2}[4]{*}{B/C}}  & \multicolumn{2}{c|}{FS-SVM} & \multicolumn{2}{c|}{St.+FS-SVM} & \multicolumn{2}{c||}{St.+FS-SVM$^*$} & \multicolumn{2}{c|}{FS-SVM} & \multicolumn{2}{c|}{St.+FS-SVM} & \multicolumn{2}{c||}{St.+FS-SVM$^*$} & \multicolumn{2}{c|}{FS-SVM} & \multicolumn{2}{c|}{St+FS-SVM} & \multicolumn{2}{c||}{St.+FS-SVM$^*$} & \multicolumn{2}{c|}{FS-SVM} & \multicolumn{2}{c|}{St.+FS-SVM} & \multicolumn{2}{c}{St.+FS-SVM$^*$} \\
			\multicolumn{2}{c||}{} & Gap   & Time  & Gap   & Time  & Gap   & Time  & Gap   & Time  & Gap   & { Time} & Gap   & Time  & Gap   & Time  & Gap   & Time & Gap   & Time  & Gap   & Time  & Gap   & { Time} & Gap   & Time \\
			\hline
			30    & $2^0$     & \textbf{0.0} & \textbf{1.44} & \textbf{0.0} & 37.00 & \textbf{0.0} & 39.15 & \textbf{0.0} & 726.84 & \textbf{0.0} & 740.17 & \textbf{0.0} & \textbf{661.64} & \textbf{0.0} & \textbf{16.32} & \textbf{0.0} & 668.41 & \textbf{0.0} & 659.59 & \textbf{0.0} & \textbf{0.56} & \textbf{0.0} & 311.60 & \textbf{0.0} & 313.48 \\
			30    & $2^1$    & \textbf{0.0} & \textbf{5.30} & \textbf{0.0} & 36.93 & \textbf{0.0} & 36.20 & \textbf{0.0} & 792.12 & \textbf{0.0} & 737.34 & \textbf{0.0} & \textbf{703.14} & \textbf{0.0} & \textbf{13.57} & \textbf{0.0} & 670.84 & \textbf{0.0} & 662.58 & \textbf{0.0} & \textbf{0.59} & \textbf{0.0} & 310.11 & \textbf{0.0} & 312.27 \\
			30    & $2^2$  & \textbf{0.0} & \textbf{3.82} & \textbf{0.0} & 36.77 & \textbf{0.0} & 35.80 & \textbf{0.0} & 814.42 & \textbf{0.0} & 1055.06 & \textbf{0.0} & \textbf{710.35} & \textbf{0.0} & \textbf{16.93} & \textbf{0.0} & 658.54 & \textbf{0.0} & 651.88 & \textbf{0.0} & \textbf{0.61} & \textbf{0.0} & 309.47 & \textbf{0.0} & 310.43 \\
			30    & $2^3$  & \textbf{0.0} & \textbf{3.87} & \textbf{0.0} & 36.96 & \textbf{0.0} & 35.22 & \textbf{0.0} & \textbf{710.88} & \textbf{0.0} & 1338.55 & \textbf{0.0} & 800.44 & \textbf{0.0} & \textbf{12.31} & \textbf{0.0} & 666.12 & \textbf{0.0} & 657.00 & \textbf{0.0} & \textbf{0.58} & \textbf{0.0} & 310.05 & \textbf{0.0} & 310.52 \\
			30    & $2^4$  & \textbf{0.0} & \textbf{4.10} & \textbf{0.0} & 37.27 & \textbf{0.0} & 35.24 & \textbf{0.0} & \textbf{791.81} & \textbf{0.0} & 880.36 & \textbf{0.0} & 1058.32 & \textbf{0.0} & \textbf{15.41} & \textbf{0.0} & 665.14 & \textbf{0.0} & 656.43 & \textbf{0.0} & \textbf{0.60} & \textbf{0.0} & 309.86 & \textbf{0.0} & 310.58 \\
			30    & $2^5$  & \textbf{0.0} & \textbf{5.26} & \textbf{0.0} & 37.51 & \textbf{0.0} & 35.63 & \textbf{0.0} & 621.57 & \textbf{0.0} & 640.47 & \textbf{0.0} & \textbf{586.30} & \textbf{0.0} & \textbf{16.17} & \textbf{0.0} & 665.92 & \textbf{0.0} & 654.55 & \textbf{0.0} & \textbf{0.59} & \textbf{0.0} & 312.92 & \textbf{0.0} & 310.41 \\
			30    & $2^6$  & \textbf{0.0} & \textbf{5.27} & \textbf{0.0} & 37.70 & \textbf{0.0} & 36.21 & \textbf{0.0} & \textbf{569.28} & \textbf{0.0} & 614.44 & \textbf{0.0} & 586.83 & \textbf{0.0} & \textbf{15.70} & \textbf{0.0} & 670.65 & \textbf{0.0} & 658.33 & \textbf{0.0} & \textbf{0.59} & \textbf{0.0} & 310.90 & \textbf{0.0} & 310.44 \\
			30    & $2^7$  & \textbf{0.0} & \textbf{5.65} & \textbf{0.0} & 37.67 & \textbf{0.0} & 35.33 & \textbf{0.0} & \textbf{527.17} & \textbf{0.0} & 726.56 & \textbf{0.0} & 617.71 & \textbf{0.0} & \textbf{16.12} & \textbf{0.0} & 665.73 & \textbf{0.0} & 655.20 & \textbf{0.0} & \textbf{0.60} & \textbf{0.0} & 310.23 & \textbf{0.0} & 310.96 \\
			\hline
			20    & $2^0$     & \textbf{0.0} & 579.94 & \textbf{0.0} & \textbf{125.26} & \textbf{0.0} & 141.74 & 1.1   & \underline{7201.03} & \textbf{0.5} & \underline{7567.88} & 0.6   & \underline{7642.56} & \textbf{0.0} & { 7200.40} & \textbf{0.0} & 1047.98 & \textbf{0.0} & \textbf{911.76} & \textbf{0.0} & 1457.92 & \textbf{0.0} & \textbf{492.42} & \textbf{0.0} & 511.68 \\
			20    & $2^1$   & 0.3   & \underline{7200.19} & \textbf{0.0} & 1946.76 & \textbf{0.0} & \textbf{1388.57} & 1.0   & \underline{7200.90} & \textbf{0.5} & \underline{7593.99} & 0.6   & \underline{7618.48} & \textbf{0.0} & 5759.65 & \textbf{0.0} & 1119.44 & \textbf{0.0} & \textbf{908.02} & \textbf{0.0} & 1478.06 & \textbf{0.0} & \textbf{470.82} & \textbf{0.0} & 533.47 \\
			20    & $2^2$  & 0.3   & \underline{7200.10} & \textbf{0.0} & 1669.72 & \textbf{0.0} & \textbf{1564.57} & 2.0   & \underline{7201.08} & \textbf{0.6} & \underline{7563.36} & \textbf{0.6} & \underline{7608.74}& \textbf{0.0} & 7168.70 & \textbf{0.0} & 1050.97 & \textbf{0.0} & \textbf{1004.88} & \textbf{0.0} & 1551.47 & \textbf{0.0} & 513.65 & \textbf{0.0} & \textbf{501.49} \\
			20    & $2^3$  & 0.3   & \underline{7200.10} & \textbf{0.0} & 1686.73 & \textbf{0.0} & \textbf{1528.63} & 2.0   & \underline{7200.97} & 0.6   & \underline{7568}.63 & \textbf{0.5} & \underline{7605.20} & 0.1   & \underline{7200.38} & \textbf{0.0} & \textbf{1006.68} & \textbf{0.0} & 1053.16 & \textbf{0.0} & 1478.82 & \textbf{0.0} & \textbf{514.35} & \textbf{0.0} & 618.31 \\
			20    & $2^4$  & 0.3   & \underline{7200.10} & \textbf{0.0} & 1584.71 & \textbf{0.0} & \textbf{1580.48} & 1.9   & \underline{7200.99} & \textbf{0.5} & \underline{7591.77} & \textbf{0.5} & \underline{7614.72} & \textbf{0.0} & \underline{7200.40} & \textbf{0.0} & 1048.67 & \textbf{0.0} & \textbf{869.59} & \textbf{0.0} & 2016.24 & \textbf{0.0} & 518.82 & \textbf{0.0} & \textbf{493.65} \\
			20    & $2^5$  & 0.3   & \underline{7200.10} & \textbf{0.0} & \textbf{1415.20} & \textbf{0.0} & 1580.30 & 1.9   & \underline{7200.82} & \textbf{0.1} & \underline{7568.42} & 0.2   & \underline{7653.38} & 0.1   & \underline{7200.38} & \textbf{0.0} & 905.88 & \textbf{0.0} & \textbf{905.16} & \textbf{0.0} & 2680.07 & \textbf{0.0} & 491.96 & \textbf{0.0} & \textbf{476.89} \\
			20    & $2^6$  & 0.4   & \underline{7200.10} & \textbf{0.0} & 1132.48 & \textbf{0.0} & \textbf{1087.20} & 1.8   & \underline{7201.05} & \textbf{0.2} & \underline{7588.98} & \textbf{0.2} & \underline{7621.43} & 0.3   & \underline{7200.42} & \textbf{0.0} & 876.49 & \textbf{0.0} & \textbf{801.61} & \textbf{0.0} & 3283.47 & \textbf{0.0} & \textbf{419.73} & \textbf{0.0} & 456.29 \\
			20    & $2^7$  & 0.5   & \underline{7200.11} & \textbf{0.0} & 1044.69 & \textbf{0.0} & \textbf{712.30} & 1.8   & \underline{7200.86} & 0.3   & \underline{7580.24} & \textbf{0.2} & \underline{7620.99} & 0.3   & \underline{7200.39} & \textbf{0.0} & 913.94 & \textbf{0.0} & \textbf{903.79} & \textbf{0.0} & 3231.34 & \textbf{0.0} & \textbf{418.32} & \textbf{0.0} & 469.33 \\
			\hline
			10    & $2^0$     & 2.5   & \underline{7203.00} & \textbf{0.0} & 4638.52 & \textbf{0.0} & \textbf{4447.71} & 15.1  & \underline{7202.81} & 15.9  & \underline{7968.70} & \textbf{11.4} & \underline{8980.67} & 11.8  & \underline{7203.61} & 8.1   & \underline{8222.67} & \textbf{6.8} & \underline{8392.61} & 2.4   & \underline{7203.84} & \textbf{0.0} & \underline{7202.61} & \textbf{0.0} & \textbf{6931.78} \\
			10    & $2^1$     & 12.4  & \underline{7204.81} & 11.2  & \underline{7250.02} & { \textbf{8.5}}   & \underline{7396.14} & 17.5  & \underline{7200.71} & 19.4  & \underline{7971.97} & \textbf{11.4} & \underline{9079.71} & 11.1  & \underline{7200.94} & 7.2   & \underline{8222.17} & \textbf{6.6} & \underline{8390.48} & 2.3   & \underline{7204.30} & 0.5   & \underline{7618.01} & \textbf{0.0} & \textbf{7331.22} \\
			10    & $2^2$  & 10.5  & \underline{7200.94} & 13.6  & \underline{7258.02} & { \textbf{9.3}}   & \underline{7442.52} & 16.9  & \underline{7200.29}& 17.5  & \underline{7987.37} & \textbf{11.3} & \underline{9212.35} & 13.5  & \underline{7200.40} & 9.5   & \underline{8214.58} & \textbf{6.6} & \underline{8382.97} & 2.3   & \underline{7203.33} & 0.5   & \underline{7620.96} & \textbf{0.0} & \textbf{6994.41} \\
			10    & $2^3$  & 15.2  & \underline{7201.77} & 12.3  & \underline{7252.36} & { \textbf{9.1}}   & \underline{7470.23} & 16.8  & \underline{7203.84} & 14.9  & \underline{7983.72} & \textbf{11.3} & \underline{9214.53} & 10.8  & \underline{7202.02} & 10.6  & \underline{8215.11} & \textbf{6.6} & \underline{8404.06} & 2.3   & \underline{7203.56} & \textbf{0.0} & \underline{7498.48} & \textbf{0.0} & \textbf{6625.83} \\
			10    & $2^4$  & 16.9  & \underline{7203.00} & 13.5  & \underline{7253.19} & { \textbf{9.2}}   & \underline{7481.98} & 12.9  & \underline{7191.54} & 15.7  & \underline{7988.91} & \textbf{11.1} & \underline{9236.43} & 13.5  & \underline{7200.41} & 9.2   & \underline{8216.77} & \textbf{6.5} & \underline{8410.13} & 2.3   & \underline{7202.94} & \textbf{0.0} & \textbf{4880.87} & \textbf{0.0} & 5016.86 \\
			10    & $2^5$  & 17.1  & \underline{7206.07} & 13.0  & \underline{7255.46} & { \textbf{9.1}}   & \underline{7434.46} & 18.8  & \underline{7200.29} & 15.2  & \underline{7988.23} & \textbf{10.9} & \underline{9243.22} & 13.5  & \underline{7200.40} & 10.4  & \underline{8229.30} & \textbf{6.2} & \underline{8405.00} & 2.9   & \underline{7206.39} & \textbf{0.0} & \textbf{2966.86} & \textbf{0.0} & 3569.74 \\
			10    & $2^6$  & 16.6  & \underline{7204.87} & 12.3  & \underline{7252.77} & {\textbf{9.1}}   & \underline{7485.48} & 18.7  & \underline{7203.20} & 17.9  & \underline{7960.66} & \textbf{11.0} & \underline{9284.75} & 13.5  & \underline{7200.40} & 10.1  & \underline{8195.61} & \textbf{6.1} & \underline{8406.76} & 2.7   & \underline{7205.12} & \textbf{0.0} & \textbf{3162.89} & \textbf{0.0} & 3457.01 \\
			10    & $2^7$  & 17.0  & \underline{7201.27} & 10.8  & \underline{7252.62} & { \textbf{9.1}}   & \underline{7437.58} & 15.2  & \underline{7203.81} & 16.5  & \underline{7982.19} & \textbf{10.9} & \underline{9261.78} & 13.5  & \underline{7200.38} & 9.9   & \underline{8214.88} & \textbf{6.2} & \underline{8407.38} & 2.6   & \underline{7206.23} & \textbf{0.0} & \textbf{2928.89} & \textbf{0.0} & 3327.63 \\
			\hline
		\end{tabular}%
	}
	\caption{Computational times and gaps using KS as initial solution of the formulation for Colon, Leukemia and DLBCL datasets.	\label{kernel_col}}%
\end{sidewaystable}%

\begin{figure}[!htb]
	\begin{minipage}{.5 \linewidth}
		\centering
		\includegraphics[width=0.9\textwidth]{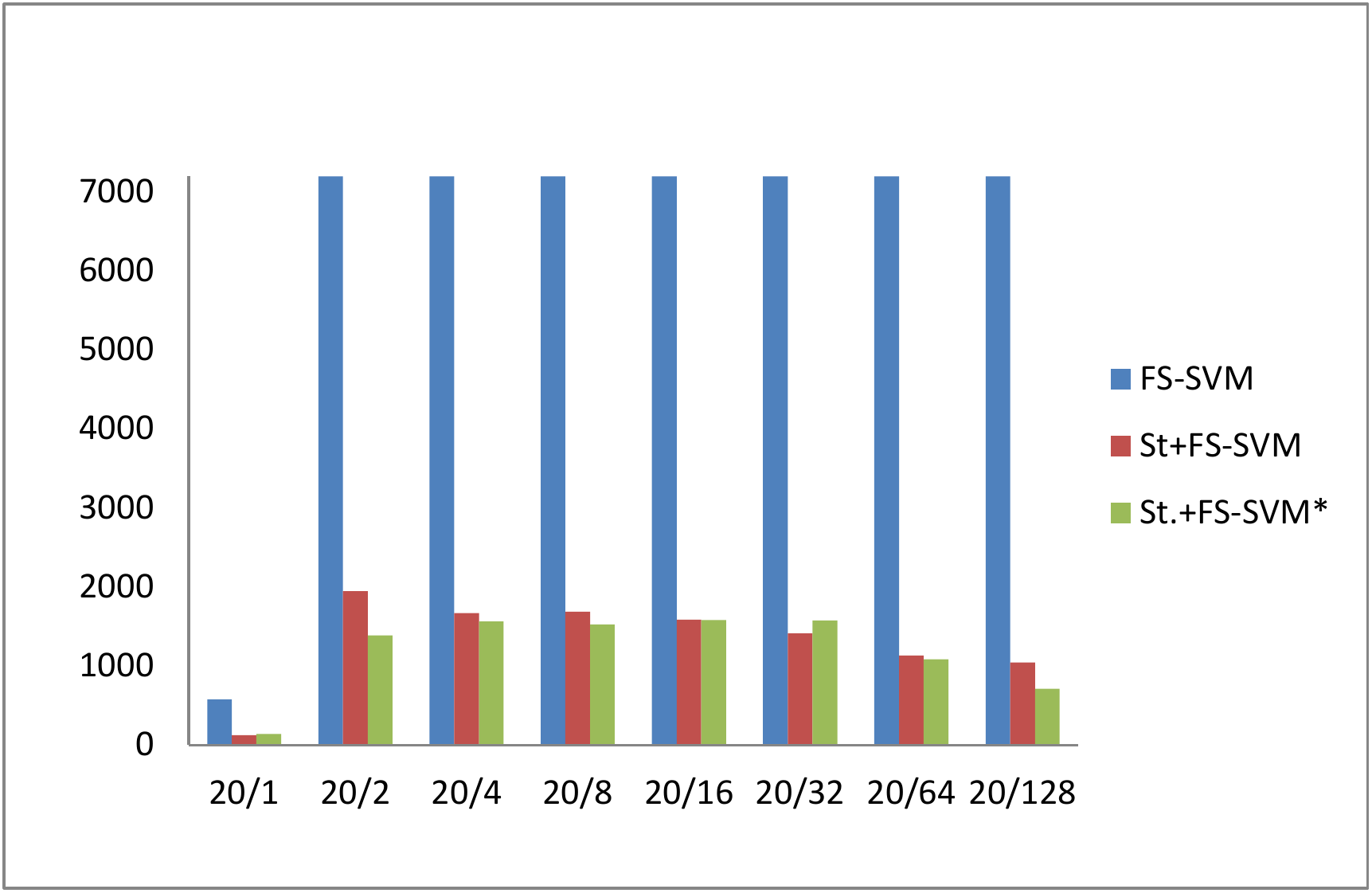}
	\end{minipage}%
	\begin{minipage}{.5 \linewidth}
		\centering
		\includegraphics[width=0.9\textwidth]{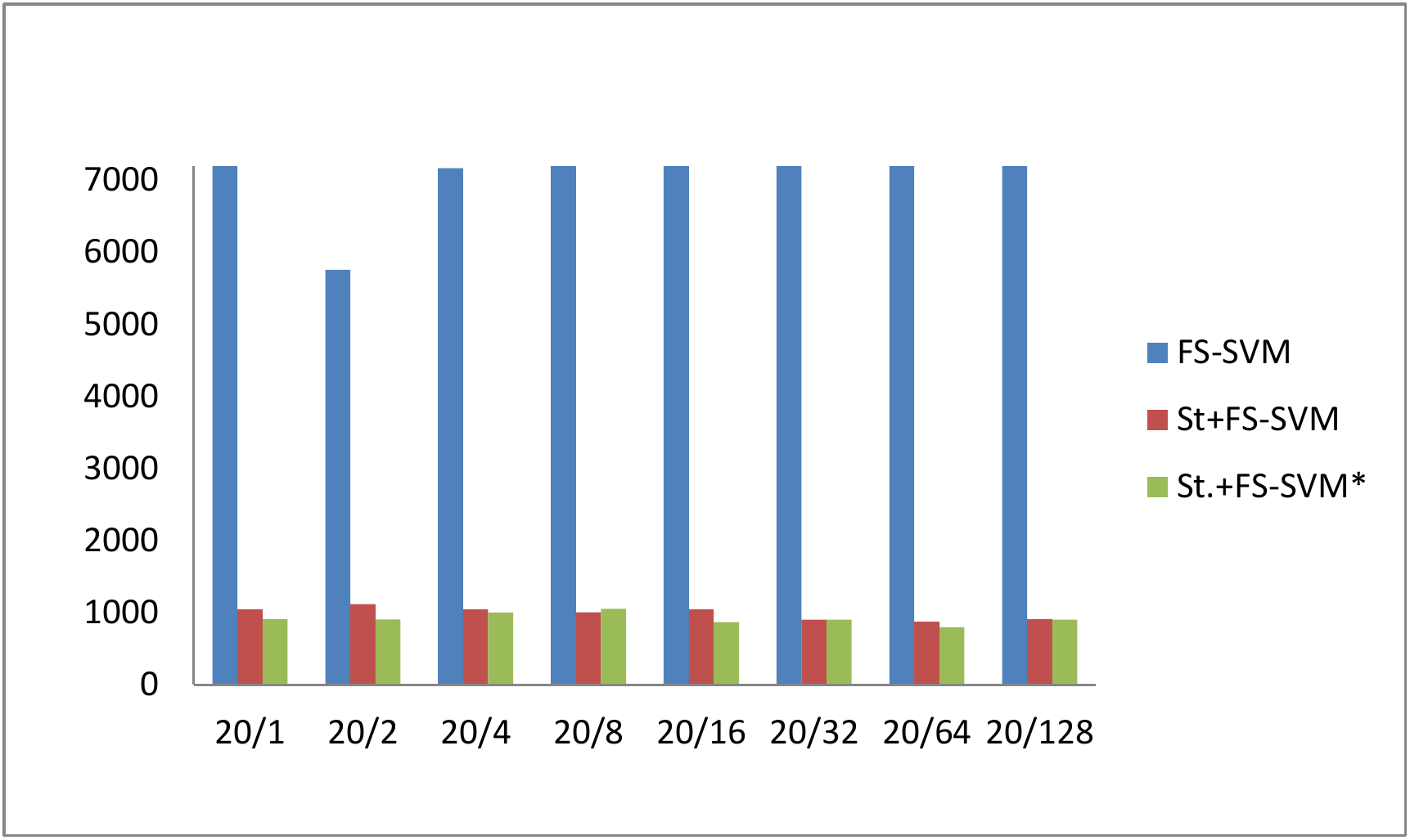}
	\end{minipage} 
	\begin{minipage}[t]{.48\linewidth}
		\caption{\footnotesize{Computational times using KS solution as initial solution of the formulation for Colon dataset with $B=20$.} 		\label{fig:COL1}}
	\end{minipage}%
	\hfill%
	\begin{minipage}[t]{.48\linewidth}
		\caption{\footnotesize{Computational times using KS solution as initial solution of the formulation for DLBCL dataset with $B=20$.} \label{fig:DL1}}
	\end{minipage}%
\end{figure}

\begin{figure}[!htb]
	\begin{minipage}{.5 \linewidth}
		\centering
		\includegraphics[width=0.9\textwidth]{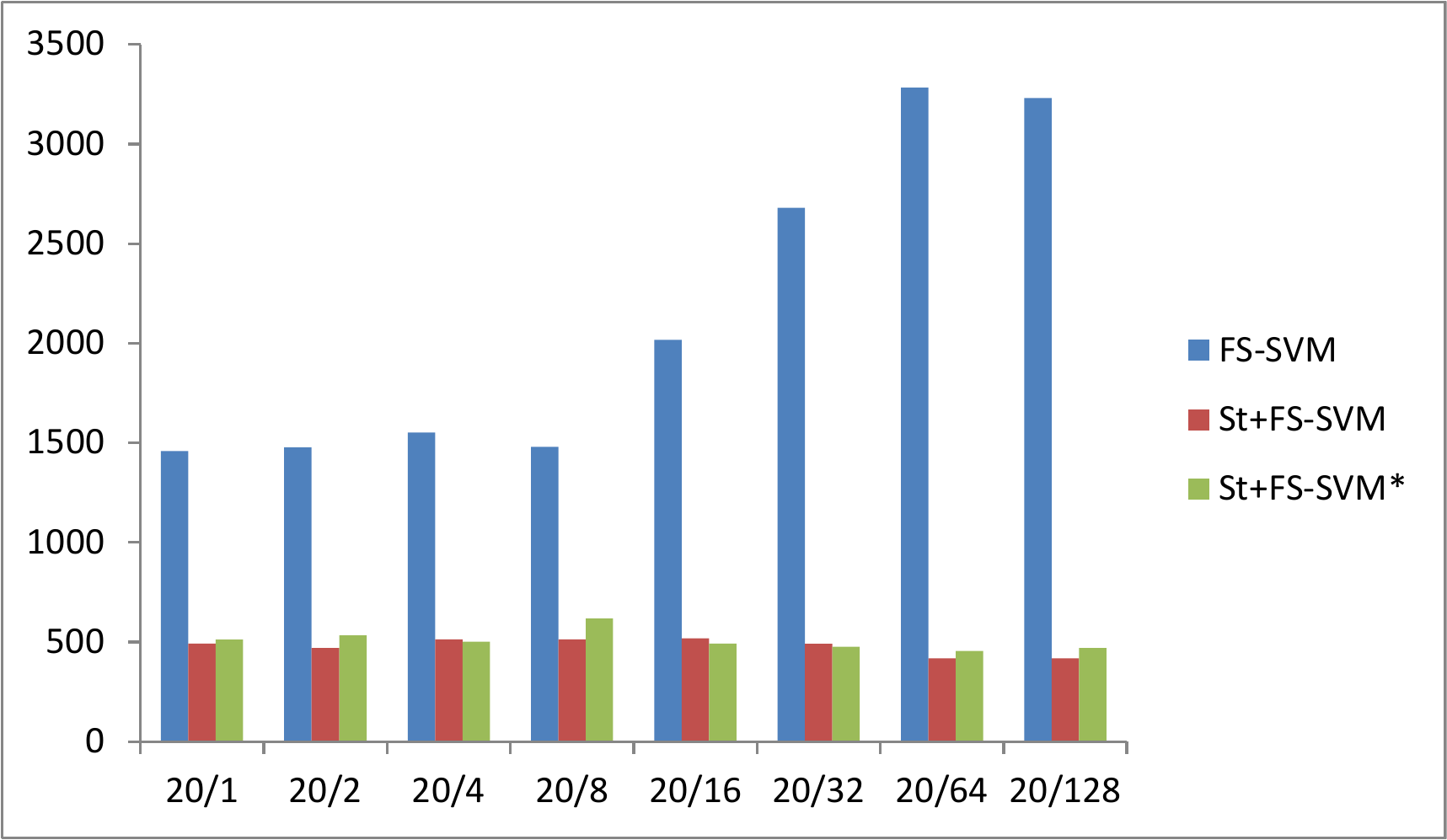}
	\end{minipage}%
	\begin{minipage}{.5 \linewidth}
		\centering
		\includegraphics[width=0.9\textwidth]{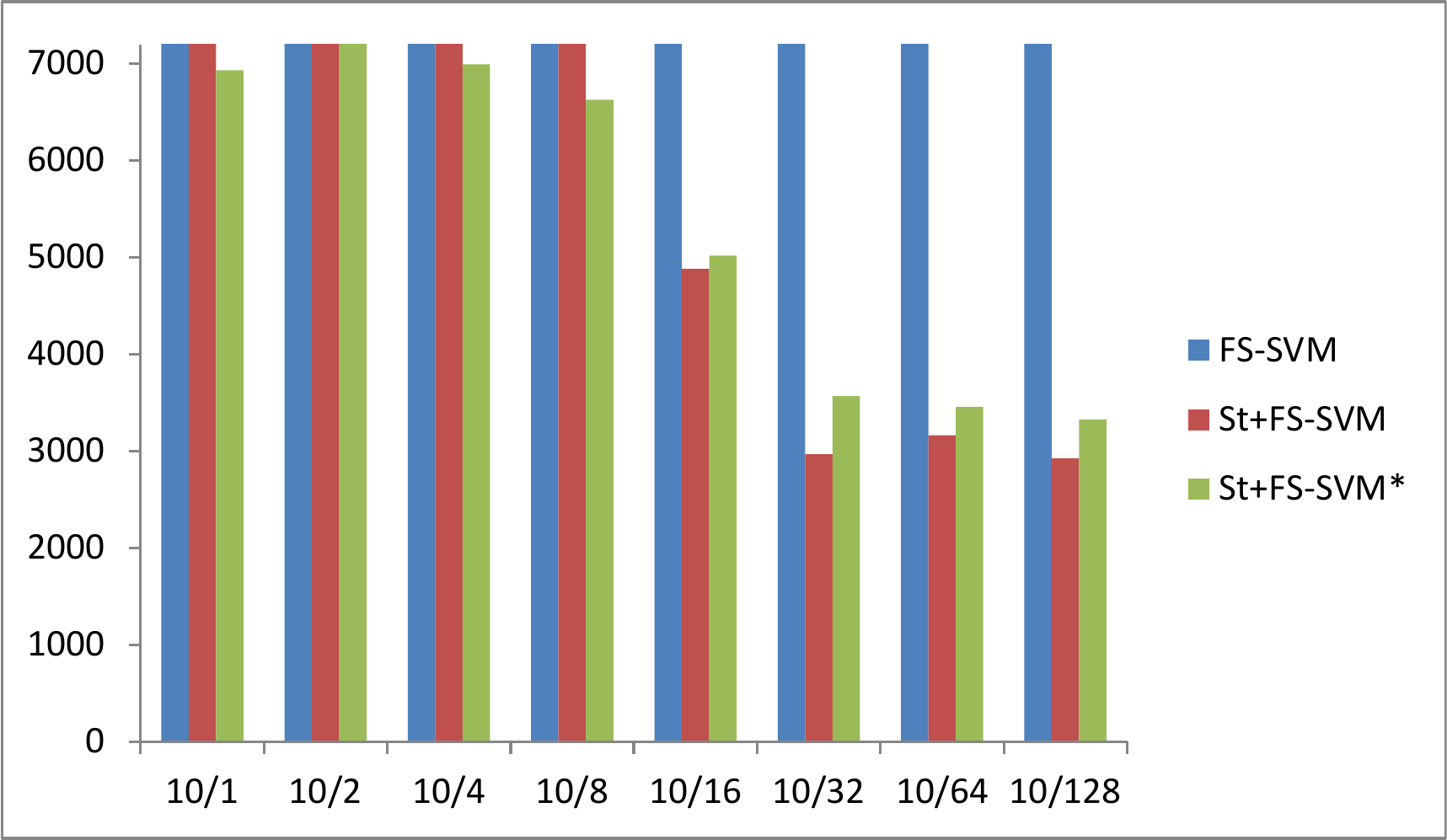}
	\end{minipage} 
	\begin{minipage}[t]{.48\linewidth}
		\caption{\footnotesize{Computational times using KS solution as initial solution of the formulation for Carcinoma dataset with $B=20$.} 		\label{fig:car1}}
	\end{minipage}%
	\hfill%
	\begin{minipage}[t]{.48\linewidth}
		\caption{\footnotesize{Computational times using KS solution as initial solution of the formulation for Carcinoma dataset with $B=10$.} \label{fig:car2}}
	\end{minipage}%
\end{figure}

\vspace{-0.5cm}
\subsubsection{Exact procedure }\label{proc}
Table \ref{tab_proc} reports the termination gaps and running times of the different procedure variants. 
The first  two-column block, St+FS-SVM$^*$, is taken from Table \ref{kernel_col}. The following three blocks of two columns show the gaps and times of the exact procedure using Variant I with the different values of parameter $S$ described in Algorithm \ref{alg2}.\ref{alg3}. 
Columns with an asterisk symbol ``*'' indicate that, for each iteration, we added a constraint to restrict the objective value to greater than or equal to the best lower bound known so far and it also means that the KS solution was used as the initial solution for SR-FS-SVM($K$).
To finish, the last two columns show the results of Variants II and III as described in Algorithms \ref{alg2}.\ref{alg4} and \ref{alg2}.\ref{alg6}, respectively.
 
The stopping rule for the exact procedure, was for the gap between the upper and lower bounds to be smaller than $0.01\%$ or if
the resolution of the semi-relaxed problem took more than $1800$ seconds. Moreover, if the bound did not improve in 5 iterations, the procedure was terminated.
 We used $B=10$ and
$C\in\{1,2,2^{2},2^{3},2^{4},2^{5},2^{6},2^{7}\}$ since the formulation could not be solved within the time limit using these parameter values. 

In the case of the Colon dataset, Table \ref{tab_proc} shows that Variant I with $S=20$ reports the best rates between the upper and lower bounds of the objective values. However, the best running times are provided by Variant III. We  can observe that the fastest procedure for the Leukemia dataset is also Variant III. Additionally, Variant II is the one that provides the best gaps for said dataset. In the case of the DLBCL dataset, Table \ref{tab_proc} shows that Variant I (taking $S=20$ and using the Kernel result as the initial solution) together with Variant II results in the best gaps. The best running times for the DLBCL dataset are also the ones obtained using Variant III. Lastly, in the Carcinoma dataset, we can observe that the exact procedure does not provide better results than the ones obtained by the formulation. In spite of that, we observed that when using Variant III, we get gaps smaller than $3\%$ in less than 2300 seconds. 

In general, these results show that exact procedures are useful for the cases
in which the formulation (together with the two strategies for fixing the bounds of $w$ variables) is not able to find the optimal
solution within the time limit. In those cases, the gaps reported by almost all variants are smaller than the ones given by the formulation and they take less time.

\begin{table}[htbp]
	\centering
	\scalebox{0.65}{
	\begin{tabular}{r@{/}c|rr|rr|rr|rr|rr|rr}
		\hline
		\multicolumn{14}{c}{Colon m=62 n=2000} \\
		\hline
	\multicolumn{2}{c|}{\multirow{2}[4]{*}{B/C}} & \multicolumn{2}{c|}{St.+FS-SVM$^*$} & \multicolumn{2}{c|}{V. I (S=20)} & \multicolumn{2}{c|}{V. I (S=20)$^*$} & \multicolumn{2}{c|}{V. I (S=40)$^*$} & \multicolumn{2}{c|}{Variant II} & \multicolumn{2}{c}{Variant III} \\
		\multicolumn{2}{c|}{}  & Gap   & Time  & Gap  & Time & Gap   & Time & Gap   & Time & Gap   & Time  & Gap   & Time \\
		\hline
		\hline
		10    & $2^0$  & \textbf{0.0} & 4447.71 & 0.8   & \underline{7281.55} & 1.4   & 5377.66 & 1.4   & 3712.85 & 0.7   & \underline{7953.52} & 3.3   & \textbf{296.96} \\
		10    & $2^1$  & 8.5   &\underline{7396.14} & 7.9   & 2679.46 & 7.1   & 3695.09 & \textbf{6.8} & 3923.84 & 7.7   & 3225.88 & 10.2  & \textbf{401.30} \\
		10    & $2^2$  & 9.3   & \underline{7442.52} & 7.7   & 4166.56 & 7.8   & 3617.22 & \textbf{7.5} & 3925.70 & 8.1   & 3663.01 & 11.0  & \textbf{429.52} \\
		10    & $2^3$  & 9.1   & \underline{7470.23} & \textbf{7.8} & 3803.96 & \textbf{7.8} & 3549.82 & \textbf{7.8} & 2218.81 & 8.1   & 3029.71 & 10.5  & \textbf{424.62} \\
		10    & $2^4$  & 9.2   & \underline{7481.98} & 8.7   & 2519.67 & \textbf{7.8} & 3819.17 & 8.0   & 2167.09 & 8.1   & 3494.11 & 10.7  & \textbf{417.16} \\
		10    & $2^5$  & 9.1   & \underline{7434.46} & \textbf{7.8} & 4216.14 & \textbf{7.8} & 4433.42 & 7.9   & 2197.12 & 8.1   & 3423.91 & 10.6  & \textbf{416.82} \\
		10    & $2^6$   & 9.1   & \underline{7485.48} & \textbf{7.8} & 3872.19 & \textbf{7.8} & 4235.20 & 7.5   & 3880.49 & 8.1   & 3355.99 & 10.9  & \textbf{397.24} \\
		10    &$2^7$   & 9.1   & \underline{7437.58} & \textbf{7.5} & 5150.32 & 7.8   & 4690.11 & 8.8   & 2159.65 & 8.1   & 2852.67 & 11.0  & \textbf{411.92} \\
		\hline
		\multicolumn{14}{c}{Leukemia m=72 n=5327} \\
		\hline
		10    & $2^0$   & 11.4  & \underline{8980.67} & \textbf{9.7} & 6015.89 & \textbf{9.7} & 5556.15 & 11.0  & 4273.43 & 12.6  & 4314.27 & 10.0  & \textbf{3903.36} \\
		10    & $2^1$   & 11.4  & \underline{9079.71} & 11.6  & 4397.74 & 10.7  & 4398.99 & 11.3  & 4395.28 & \textbf{10.4} & 4849.97 & 10.5  & \textbf{3776.61} \\
		10    & $2^2$  & 11.3  & \underline{9212.35} & 11.0  & 4499.99 & 10.6  & 4595.05 & 11.0  & 4491.02 & \textbf{10.4} & 4962.93 & 10.5  & \textbf{4006.07} \\
		10    & $2^3$  & 11.3  & \underline{9214.53} & 11.5  & 4445.23 & 11.1  & 4411.35 & 11.0  & 4455.09 & \textbf{10.4} & 5002.57 & 10.5  & \textbf{3862.61} \\
		10    & $2^4$   & 11.1  & \underline{9236.43} & 11.7  & 4452.47 & 10.5  & 4476.44 & 10.8  & 4408.88 & \textbf{10.4} & 5054.30 & 10.5  & \textbf{3741.09} \\
		10    & $2^5$  & 10.9  & \underline{9243.22} & 11.6  & 4513.21 & \textbf{10.4} & 4528.56 & 10.6  & 4543.10 & \textbf{10.4} & 4978.80 & 10.5  & \textbf{3734.30} \\
		10    & $2^6$ & 11.0  & \underline{9284.75} & 11.4  & 4383.05 & \textbf{10.1} & 4403.61 & 10.3  & 4436.12 & 10.4  & 4829.99 & 10.5  & \textbf{3656.79} \\
		10    & $2^7$ & 10.9  & \underline{9261.78} & 11.0  & 4486.71 & { \textbf{9.9} }  & 4508.91 & 10.5  & 4523.24 & 10.4 & 4899.46 & 10.5  & { \textbf{3670.81}} \\
		\hline
		\multicolumn{14}{c}{DLBCL m=77 n=7129} \\
		\hline
		10    & $2^0$  & 6.8   & \underline{8392.61} & 6.9   & 4502.23 & \textbf{5.4} & 6150.31 & 6.9   & 4467.66 & 6.2   & 4958.08 & 6.9   & \textbf{3924.76} \\
		10    & $2^1$   & 6.6   & \underline{8390.48} & 5.4   & 6134.29 & 5.4   & 6054.37 & 6.6   & 4509.11 & \textbf{5.2} & 6791.88 & 6.9   & \textbf{3940.43} \\
		10    & $2^2$  & 6.6   & \underline{8382.97} & \textbf{5.4} & 6251.76 & \textbf{5.4} & 5906.59 & 6.0   & 4476.80 & 6.2   & 5125.47 & 6.9   & \textbf{3927.96} \\
		10    &$2^3$  & 6.6   & \underline{8404.06} & 6.1   & 4498.68 & \textbf{5.4} & 6285.88 & 6.6   & 4516.28 & 6.2   & 5570.18 & 6.9   & \textbf{3722.99} \\
		10    & $2^4$   & 6.5   & \underline{8410.13} & 6.2   & 4467.94 & 5.4   & 6054.82 & 6.1   & 4479.48 & \textbf{5.2} & 7214.17 & 6.9   & \textbf{3869.21} \\
		10    &$2^5$   & 6.2   & \underline{8405.00} & 6.1   & 4508.95 & 5.4   & 5500.90 & 5.8   & 4537.00 & \textbf{5.2} & 6835.76 & 6.9   & \textbf{3796.52} \\
		10    & $2^6$ & 6.1   & \underline{8406.76} & \textbf{5.4} & 6093.68 & \textbf{5.4} & 5555.90 & 5.7   & 4487.47 & 6.2   & 5039.73 & 6.9   & \textbf{3628.86} \\
		10    & $2^7$  & 6.2   & \underline{8407.38} & 6.3   & 4498.92 & \textbf{5.4} & 5561.95 & 5.5   & 4546.70 & 6.2   & 5112.92 & 6.9   & \textbf{3746.41} \\
		\hline
		\multicolumn{14}{c}{Carcinoma m=36 n=7457} \\
		\hline
		10    &$2^0$  & \textbf{0.0} & 6931.78 & 1.0   & 5793.66 & 1.6   & 4511.14 & 1.6   & 3697.36 & 1.3   & 5963.29 & 2.5   & \textbf{1989.73} \\
		10    & $2^1$  & \textbf{0.0} & 7331.22 & 1.0   & 5473.46 & 1.6   & 5190.58 & 1.6   & 3649.83 & 1.3   & 6682.76 & 2.9   & \textbf{1841.85} \\
		10    &$2^2$  & \textbf{0.0} & 6994.41 & 1.1   & 5900.01 & 1.6   & 5137.24 & 1.6   & 3670.44 & 0.9   & 7563.22 & 2.7   & \textbf{2068.55} \\
		10    & $2^3$  & \textbf{0.0} & 6625.83 & 1.1   & 5155.73 & 2.1   & 3470.51 & 1.6   & 3655.70 & 0.9   & 7654.51 & 2.2   & \textbf{1882.08} \\
		10    & $2^4$  & \textbf{0.0} & 5016.86 & 1.1   & 5494.97 & 1.6   & 4912.55 & 1.6   & 3591.92 & 1.3   & 6222.99 & 2.3   & \textbf{1881.46} \\
		10    & $2^5$  & \textbf{0.0} & 3569.74 & 1.1   & 5368.17 & 1.6   & 4998.79 & 1.6   & 3586.95 & 1.3   & 5641.86 & 2.0   & \textbf{2295.48} \\
		10    & $2^6$ & \textbf{0.0} & 3457.01 & 0.8   & 7084.94 & 1.6   & 4545.17 & 1.6   & 3585.37 & 0.9   & 6623.43 & 2.6   & \textbf{1961.13} \\
		10    & $2^7$  & \textbf{0.0} & 3327.63 & 1.1   & 5799.65 & 1.6   & 5047.74 & 1.6   & 3533.14 & 0.9   & 7392.51 & 2.0   & \textbf{2404.20} \\
		\hline
	\end{tabular}%
       }
	\caption{ Different variants of the exact procedure for the Colon, Leukemia, DLBCL and Carcinoma datasets.}\label{tab_proc}%
\end{table}%

\subsection{Analysis of datasets with big sample sizes and big number of features }\label{bigsample}
In this subsection we will focus on the datasets with big sample size and big number of features (Arrythmia, Madelon and MFeat datasets). 
Table \ref{tab:newdata_proc} reports four blocks of columns with some of the techniques described in the previous subsection. For each dataset, the blocks of columns FS-SVM and St+FS-SVM* are defined similarly to Table \ref{kernel_col},  the Kernel block reports the results of Kernel Search algorithm. Besides, we include $\Delta B$ column like in Table \ref{tab_prep}. Finally, the last block of columns report the results of exact procedure V. I (S=20)* (see Table \ref{tab_proc}). 

If we focus our attention on the Arrythmia dataset, we observe that none of the instances can be solved either using the strategies or not. However, the final gaps when using the strategies are smaller than the ones without the strategies. A similar behaviour can be appreciated in Madelon dataset. However, for the MFeat dataset we observe that the instances can be solved in less than two hours if strategies and Kernel are applied. In $\Delta B$ column of MFeat dataset, we can observe that the bounds of $w$ variables are very tightened when using the strategies.

Regarding the Kernel Search results, we observe that, in almost all the cases, the solution provided by the heuristic gives the same final gap as the formulation in smaller times. For this reason, we can conclude that the upper bound provided by the Kernel Search is quite good because it cannot be improved.

The last block of columns of Table \ref{tab:newdata_proc} reports the results of exact procedure Variant I (S=20)*. We have chosen this variant because it seems that it provides  the best gaps in Table \ref{tab_proc}. For Arrythmia and Madelon datasets, in most of the cases, the gaps provided by this exact procedure are better than the final gaps of the model using the strategies. In the MFeat dataset although the gaps have not improved these are very small. Moreover, some instances provide better  running times.

\begin{table}[htbp]
	\centering
	\scalebox{0.7}{
	\begin{tabular}{r|rr|rrr|rr|rr}
		\hline
		& \multicolumn{9}{c}{Arrythmia m=420 n=258} \\
		\hline
	 & \multicolumn{2}{c|}{FS-SVM} & \multicolumn{3}{c|}{St+FS-SVM*} & \multicolumn{2}{c|}{Kernel} & \multicolumn{2}{c}{V. I (S=20)*} \\
		\hline  	\multicolumn{1}{c|}{B/C}        & \multicolumn{1}{r}{Gap} & \multicolumn{1}{r|}{Time} & \multicolumn{1}{r}{Gap} & \multicolumn{1}{r}{t$_{total}$} & \multicolumn{1}{r|}{$\Delta B$} & Gap   & \multicolumn{1}{r|}{Time} & \multicolumn{1}{r}{Gap} & \multicolumn{1}{r}{Time} \\
		\hline
		30/1  & 8.4   & 7201.24 & 8.3   & 7342.68 & 33.34 & \multicolumn{1}{r}{8.8/0.6} & 128.50 & 7.0   & 1851.22 \\
		30/2  & 19.6  & 7201.14 & 17.4  & 8121.25 & 112.51 & \multicolumn{1}{r}{17.4/0.0} & 909.35 & 17.6  & 1868.85 \\
		30/4  & 35.4  & 7200.80 & 31.1  & 8114.76 & 367.35 & \multicolumn{1}{r}{31.1/0.0} & 902.93 & 30.9  & 2000.83 \\
		30/8  & 52.7  & 7200.80 & 47.2  & 8117.53 & 1079.00 & \multicolumn{1}{r}{47.2/0.0} & 909.76 & 48.0  & 3610.61 \\
		30/16 & 70.1  & 7201.10 & 67.5  & 8111.17 & 3169.37 & \multicolumn{1}{r}{67.5/0.0} & 904.49 & 67.4  & 3304.04 \\
		30/32 & 85.1  & 7201.02 & 84.6  & 7252.30 & 8388.60 & \multicolumn{1}{r}{84.8/1.3} & 44.23 & 83.2  & 3631.16 \\
		30/64 & 92.5  & 7200.82 & 92.3  & 7257.44 & 18027.74 & \multicolumn{1}{r}{92.3/0.0} & 49.87 & 91.6  & 3274.35 \\
		30/128 & 96.3  & 7200.89 & 96.1  & 7244.39 & 37249.68 & \multicolumn{1}{r}{96.1/0.0} & 38.03 & 95.7  & 2849.88 \\
		\hline
		& \multicolumn{9}{c}{Madelon m=2000 n=500} \\
		\hline
		 & \multicolumn{2}{c|}{FS-SVM} & \multicolumn{3}{c|}{St+FS-SVM*} & \multicolumn{2}{c|}{Kernel} & \multicolumn{2}{c}{V. I (S=20)*} \\
		\hline   	\multicolumn{1}{c|}{B/C}       & \multicolumn{1}{r}{Gap} & \multicolumn{1}{r|}{Time} & \multicolumn{1}{r}{Gap} & \multicolumn{1}{r}{t$_{total}$} & \multicolumn{1}{r|}{$\Delta B$} & Gap   & \multicolumn{1}{r|}{Time} & \multicolumn{1}{r}{Gap} & \multicolumn{1}{r}{Time} \\
		\hline
		30/1  & 24.6  & 7200.22 & 23.0  & 10602.11 & 721.48 & \multicolumn{1}{r}{23.0/0.0} & 1800.51 & 23.0  & 10580.60 \\
		30/2  & 26.6  & 7200.27 & 25.5  & 10593.11 & 1600.91 & \multicolumn{1}{r}{25.5/0.0} & 1800.54 & 25.3  & 8823.08 \\
		30/4  & 29.2  & 7200.26 & 27.0  & 10587.78 & 3389.77 & \multicolumn{1}{r}{27.0/0.0} & 1801.07 & 27.2  & 9690.67 \\
		30/8  & 29.7  & 7200.16 & 27.5  & 10554.38 & 6892.47 & \multicolumn{1}{r}{27.5/0.0} & 1801.37 & 27.7  & 9115.10 \\
		30/16 & 30.2  & 7200.26 & 27.9  & 10540.10 & 13938.45 & \multicolumn{1}{r}{27.9/0.0} & 1800.56 & 28.4  & 9734.53 \\
		30/32 & 30.1  & 7200.13 & 28.7  & 10556.32 & 28932.99 & \multicolumn{1}{r}{28.7/0.0} & 1800.95 & 28.2  & 9740.28 \\
		30/64 & 30.5  & 7200.24 & 28.7  & 10609.76 & 57902.00 & \multicolumn{1}{r}{28.7/0.0} & 1800.68 & 28.6  & 8817.60 \\
		30/128 & 30.6  & 7200.21 & 28.4  & 10577.54 & 114302.34 & \multicolumn{1}{r}{28.4/0.0} & 1800.91 & 28.8  & 9732.01 \\
		\hline
		& \multicolumn{9}{c}{Mfeat m=2000 n=649} \\
		\hline
& \multicolumn{2}{c|}{FS-SVM} & \multicolumn{3}{c|}{St+FS-SVM*} & \multicolumn{2}{c|}{Kernel} & \multicolumn{2}{c}{V. I (S=20)*} \\
	\hline       	\multicolumn{1}{c|}{B/C}   & \multicolumn{1}{r}{Gap} & \multicolumn{1}{r|}{Time} & \multicolumn{1}{r}{Gap} & \multicolumn{1}{r}{t$_{total}$} & \multicolumn{1}{r|}{$\Delta B$} & Gap   & \multicolumn{1}{r|}{Time} & \multicolumn{1}{r}{Gap} & \multicolumn{1}{r}{Time} \\
		\hline
		30/1  & 0.1   & 7201.11 & 0.0   & 4095.02 & 0.10  & 0.0   & 1120.02 & 0.5   & 4508.91 \\
		30/2  & 0.1   & 7201.19 & 0.0   & 4752.77 & 0.10  & 0.0   & 1226.71 & 0.5   & 3428.11 \\
		30/4  & 0.2   & 7201.19 & 0.0   & 5018.38 & 0.10  & 0.0   & 1039.17 & 0.5   & 4863.50 \\
	30/8  & 0.1   & 7201.07 & 0.0   & 6189.81 & 0.10  & 0.0   & 1116.03 & 0.5   & 4876.60 \\
		30/16 & 0.2   & 7201.19 & 0.0   & 5563.56 & 0.10  & 0.0   & 1103.97 & 0.5   & 4619.32 \\
		30/32 & 0.1   & 7201.15 & 0.0   & 3362.12 & 0.10  & 0.0   & 1040.65 & 0.5   & 4984.04 \\
		30/64 & 0.1   & 7202.74 & 0.0   & 4207.90 & 0.10  & 0.0   & 1151.23 & 0.6   & 4731.86 \\
		30/128 & 0.1   & 7201.15 & 0.0   & 3763.38 & 0.10  & 0.0   & 1219.63 & 0.6   & 4656.22 \\
		\hline
	\end{tabular}%
}
\caption{ Strategies, Kernel and procedure V. I(20)* results for Arrythmia, Madelon and MFeat datasets.}
	\label{tab:newdata_proc}%
\end{table}%

\section{Validation of the model}\label{validation}
\label{comparison}
In this section, the proposed model is analyzed together with other classification approaches. To this end, a 10-fold-cross validation (10-FCV) was performed for the datasets described above. We compare our model not only with classical SVM such as $\ell_1$-SVM (\cite{BraMa98}), 
$\ell_2$-SVM (\cite{Vap98}) and  {LP}-SVM (\cite{Zhou02}),  but also 
with recent SVM that include feature selection constraints, such as MILP1 described in Section \ref{svm} and MILP2 (an extension of LP-SVM developed by \cite{Zhou02} with feature selection through a budget constraint), see \cite{Maldo14}. Furthermore, the following classification techniques are also analyzed and compared with our model: FSV (\cite{BraMa98}),
 RFE-SVM (\cite{Guyon02}) and Fisher Criterion Score (\cite{Guyon06}), referred to as Fisher-SVM, 

We compute two different predictive performance measures: the accuracy (ACC) and the area under the curve (AUC). The ACC is given by the percentage of good classified elements of the test set. 
However, the AUC is the mean of the percentage of good classified positive elements and the percentage of good classified negative elements. Hence,
$$ACC=\frac{TP+TN}{TP+TN+FP+FN},$$
$$AUC=\frac{\frac{TP}{TP+FN}+\frac{TN}{TN+FP}}{2},$$
where TP are true positives, TN are true negatives, FP false positives and FN false negatives.

The different models are compared using these two measures. We should point out that each model has different parameters and we chose the parameters with the best performance in each model. In fact, the penalty parameter $C$, which  is necessary for the MILP2 and FS-SVM models, is varied in the set
$\{2^{-7},2^{-6},2^{-5},2^{-4},2^{-3},2^{-2},2^{-1},2^0,2^{1},$
$2^{2},2^{3},2^{4},2^{5},2^{6},2^{7}\}$. Moreover, the upper bounds on the variables related  to the MILP1
and MILP2 models are fixed to a sufficiently large amount. The budget parameter (B) is varied,  as shall be explained in the following subsections.

To select the best parameters for a model, 10-FCV is performed with each possible combination of parameter values. For each parameter value and each fold of 10-FCV, we obtain  an ACC and AUC measure value. After the application of 10-FCV to this parameter combination, we obtain an average of the measures associated with each fold. We repeat this process for the remaining parameter combinations and then compare the average performance measures. Lastly, we indicate the average measures associated with the parameters that provide best results for each instance.

 In Subsection \ref{SmallInstances}, the eight datasets with the smallest number of features will be analyzed. Subsection \ref{BigInstances} will address the Colon, Leukemia, DLBCL and Carcinoma datasets. Finally, Subsection \ref{SuperBigInstances} will be focused on the remaining datasets that present big sample sizes and  big number of features.
\subsection{Instances with a small number of features}\label{SmallInstances}

This subsection concerns the comparative analysis of different classification techniques and our model for the following datasets: BUPA, PIMA, Cleveland, Housing, Australian Credit, German Credit, WBC
and Ionosphere. The number of features in all these datasets is below 50. Therefore, the budget parameter (B) varies with all the possible numbers of features, i.e. $B=1,\ldots,n$. 

The first column of Tables \ref{tab_bupa}-\ref{tab_iono} specifies the classification method used.
For each method, we have provided the best average ACC achieved in the 10-FCV (column ACC) and the parameters used to obtain this ACC (columns $B$ and $C$). Using these  parameters, column AUC reports the average AUC value. The fourth column gives the average number of selected features in the 10-FCV and lastly, the last column, labelled ``Time'', details the average time that is required to run  a single fold of the 10-FCV.
 It should be noted that the last two rows of each table correspond to our formulation (FS-SVM) and the Kernel Search (KS FS-SVM).
 The model with best performance for each dataset is shown in bold.  

 Furthermore, figures \ref{fig:bupa_acc}-\ref{fig:iono_acc} show the ACC performance of different models that include feature selection in terms of parameter $B$. We should emphasize that parameter $B$ represents the maximum number of features that can be selected in the model. Note that the results shown in these graphs are the average ACC of 10-FCV corresponding to the $C$ values reported in the corresponding Tables \ref{tab_bupa}-\ref{tab_iono}.

For the BUPA dataset (Table \ref{tab_bupa}), the best performance is achieved using FS-SVM together with $\ell_1$-SVM.  Figure
\ref{fig:bupa_acc} reveals that MILP2 and FS-SVM are the models that provide the best performances for the different $B$ values.
In the case of the PIMA dataset (Table \ref{tab_pima}), the best performance in terms of ACC is also obtained using our model together with $\ell_1$-SVM, while the best performance in terms of AUC is achieved by RFE-SVM.
Figure \ref{fig:pima_acc} shows that the best ACC performance in general terms is obtained using the proposed FS-SVM model.

In Table \ref{tab_clev}, the best performance for the Cleveland dataset is reported by the RFE-SVM method. Additionally, Figure \ref{fig:clev_acc} shows that MILP1, MILP2 and FS-SVM present the best average ACC when varying the $B$ parameter. In the case of the Housing dataset (Table \ref{tab_hous}), we can observe that RFE-SVM provides the best average ACC and AUC. FS-SVM also shows similar values but with half the number of features. In Figure \ref{fig:hous_acc}, we can conclude that our approach for $B\leq 7$ is among the best ones, together with MILP2.

Table \ref{tab_aus} shows the results for the Australian Credit dataset. In this case, all the models present a similar ACC and AUC performance. Although the results in bold are those related with the
 Fisher-SVM  technique, Figure \ref{fig:australian_acc} however, shows that Fisher-SVM performs worse than other models when the budget parameter has a lower value ($B=1,\ldots,4$). Therefore, there is no significant difference between the models.  

For the German Credit dataset, the best ACC performance is given by MILP1 and the best AUC performance is produced by Fisher-SVM, as shown in Table \ref{tab_gc}. Figure \ref{fig:gc_acc} shows that  the best ACC is achieved with models MILP1, MILP2 and FS-SVM. 
The best performance for the WBC dataset is obtained with the $\ell_2$-SVM model. However, $30$ features are required to achieve these results. Our model also provides good 
results when using only $4$ of the $30$ features. In Figure \ref{fig:wbc_acc}, the overall best results are obtained with the MILP2 and FS-SVM models, using
different $B$ values.

The dataset with the biggest number of features among the smallest instances (the Ionosphere dataset) performs best when using FS-SVM, as can be seen in Table \ref{tab_iono}. 
 Figure \ref{fig:iono_acc} also shows that all the models provide an ACC performance between $82\%$ and $88\%$. In general, FS-SVM always has the best performance or else is among the best performers.
 	
 In terms of these datasets, we can conclude that the proposed model (FS-SVM) produces, in general terms, a good and stable performance, improving, in some cases, the existing approaches in the literature.

\begin{table}[!htb]
	\begin{minipage}{.5\linewidth}
		\centering
		\scalebox{0.65}{
			\begin{tabular}{rrrrrrr}
						\hline
						\multicolumn{7}{c}{BUPA m=345 n=6} \\
						\hline
						Form. & ACC   & AUC   & Av. F.    & B   & C  & Time \\
						\hline
						{ $\ell_1$-SVM} & \textbf{69.62\%} & \textbf{66.75\%} & 6      &   -    & 2 &0.01 \\
						{ $\ell_2$-SVM} & 69.33\% & 66.58\% & 6        &   -  &1  & 0.01 \\
						{  LP-SVM} & 69.33\% & 66.58\% & 6          &  - & 1   & 0.01 \\
						MILP1  & 69.33\% & 66.58\% & 6          & 6  &  - & 0.03 \\
						MILP2 & 69.33\% & 66.58\% & 6         & 6   & 1 & 0.01 \\
						FSV   & 49.48\% & 54.92\% & 3          &  -  &- & 0.01 \\
						{ Fisher-SVM} & 66.33\% & 66.50\% & 6         & 6  & 2  & 0.05 \\
						{ RFE-SVM}   & 68.52\% & 66.75\% & 4        & 4 & $2^{-1}$   & 0.07 \\
						FS-SVM   & \textbf{69.62\%} & \textbf{66.75\%} & 2& 6         & 6     & 0.01 \\
						KS-FS-SVM & \textbf{69.62\%} & \textbf{66.75\%}& 2 & 6          & 6     & 0.01 \\
						\hline         
			\end{tabular}%
		}
	\end{minipage}%
	\begin{minipage}{.5\linewidth}
		\centering
		\scalebox{0.65}{
			\begin{tabular}{rrrrrrr}
					    \hline
	\multicolumn{7}{c}{PIMA m=768 n=8}  \\
		\hline
		Form.	      &Av. ACC	&Av. AUC	 &Av. F. &	B& C&	Time\\
		\hline
		{ $\ell_1$-SVM}	    &{\bf 77.75\%}&	72.79\%&	8&	-&	2&0.01\\
		{ $\ell_2$-SVM}	    &77.49\%&	72.50\%&	8&	-&2&	0.01\\
		{LP-SVM}	    &77.10\%&	72.11\%&	8&	-& 1& 0.02\\
		MILP1 	    &77.49\%&	72.59\%&	7&	7&- &	0.07\\
		MILP2	      &77.62\%& 72.69\%&	7&	7&$2^{-4}$	&0.06\\
		FSV	        &75.01\%& 70.97\%&	4&	-&	-&0.02\\
		{ Fisher-SVM}	&76.59\%&	74.49\%&	5&	5&	$2^2$&0.47\\
		{ RFE-SVM}    &76.98\%&	{\bf 74.87\%}&	5&	5&$2^7$	&0.44\\
		FS-SVM	  &{\bf 77.75\%}&	72.79\%&	8&	8&2	&0.05\\
		KS FS-SVM 	&{\bf 77.75\%}&	72.79\%&	8&	8&2	&0.06\\	
		\hline		
			\end{tabular}%
		}
	\end{minipage} 
	\begin{minipage}[t]{.45\linewidth}
		\caption{\footnotesize{Best average ACC and AUC for BUPA dataset. 	\label{tab_bupa}}}
	\end{minipage}%
	\hfill%
	\begin{minipage}[t]{.45\linewidth}
		\caption{\footnotesize{Best average ACC and AUC for PIMA dataset.	\label{tab_pima}}}
	\end{minipage}%
\end{table}

\begin{table}[!htb]
    \begin{minipage}{.5\linewidth}
      \centering
			\scalebox{0.65}{
        \begin{tabular}{rrrrrrr}	
        	\hline	
			\multicolumn{6}{c}{Cleveland m=297 n=13} \\
			\hline
			Formulación & ACC   & AUC   & Av. F. & B     &C& Time \\
			\hline
			{ $\ell_1$-SVM} & 84.35\% & 83.67\% & 12.9  &   -    & 4&0.83 \\
			{ $\ell_2$-SVM} & 84.69\% & 84.08\% & 13    &   -    &2 &0.02 \\
			{ LP-SVM} & 84.01\% & 83.28\% & 13    &   -    & 1 &0.00 \\
			MILP1  & 85.02\% & 84.51\% & 10    & 10    &- &0.05 \\
			MILP2 & 84.69\% & 84.05\% & 10    & 10    &1 &0.06 \\
			FSV   & 68.87\% & 66.01\% & 2   &       & -&0.01 \\
			{ Fisher-SVM}& 84.69\% & 84.06\% & 10    & 10    &8 &0.03 \\
			{ RFE-SVM}  & \textbf{85.25\%} & \textbf{84.53\%} & 11    & 11& $2^7$   & 0.48 \\
			FS-SVM & 84.35\% & 83.75\% & 10    & 10    &$2^7$ &0.04 \\
			KS-FS-SVM & 84.35\% & 83.75\% & 10    & 10    &$2^7$ &0.08 \\
      \hline     
    \end{tabular}%
		}
    \end{minipage}%
    \begin{minipage}{.5\linewidth}
     \centering
				\scalebox{0.65}{
				\begin{tabular}{rrrrrrr}
		\hline
		\multicolumn{7}{c}{Housing m=506 n=13} \\
		\hline
		Formu. & ACC   & AUC   & Av. F. & B   &  C& Time \\
		{ $\ell_1$-SVM} & 85.34\% & 85.34\% & 13    & -  & $2^7$   & 0.01 \\
		{ $\ell_2$-SVM} & 85.91\% & 85.90\% & 13    & -  & $2^{-2}$   & 0.04 \\
		{ LP-SVM} & 85.34\% & 85.34\% & 13    & -  & $2^7$   & 0.01 \\
		MILP1  & 86.12\% & 86.14\% & 6     & 6  & -  & 0.07 \\
		MILP2 & 86.12\% & 86.14\% & 6     & 6   &  $2^7$  & 0.08 \\
		FSV   & 59.45\% & 58.84\% & 3.1   & -   &  - & 0.01 \\
		{ Fisher-SVM} & 86.12\% & 86.11\% & 8     & 8 &  $2^6$    & 0.12 \\
		{ RFE-SVM}  & \textbf{86.50\%} & \textbf{86.45\%} & 12    & 12& $2^6$   & 0.94 \\
		FS-SVM & 86.32\% & 86.34\% & 6     & 6 &  $2^2$  & 0.11 \\
		KS-FS-SVM & 86.32\% & 86.34\% & 6     & 6& $2^2$     & 0.15 \\
\hline
    \end{tabular}%
		}
    \end{minipage} 
		\begin{minipage}[t]{.45\linewidth}
        \caption{\footnotesize{Best average ACC and AUC for Cleveland dataset. 	\label{tab_clev}}}
    \end{minipage}%
		\hfill%
    \begin{minipage}[t]{.45\linewidth}
        \caption{\footnotesize{Best average ACC and AUC for Housing dataset.	\label{tab_hous}}}
    \end{minipage}%
\end{table}

\begin{table}[!htb]
	\begin{minipage}{.5\linewidth}
		\centering
		\scalebox{0.65}{
			\begin{tabular}{rrrrrrr}	
				\hline	
					\multicolumn{7}{c}{Australian Credit m=690 n=14}        \\
				\hline
				Form.	    & Av. ACC& 	Av. AUC& 	Av. F.& 	B&C	&  Time\\
				\hline
				{ $\ell_1$-SVM}   &	85.51\%&	86.21\%&	1		  &  -  & $2^{-1}$& 0.01\\
				{ $\ell_2$-SVM}   &	85.51\%&	86.21\%&	14	  &  -	& 2 &0.13\\
				{ LP-SVM}   &	85.22\%&	85.89\%&	12	  &  -	& $2^7$ &0.01\\
				MILP1     &	85.51\%&	86.21\%&	1	    &  1	& - &0.06\\
				MILP2	    & 85.51\%& 	86.21\%& 	1	    &  1	& $2^7$ &0.04\\
				FSV	      & 85.51\%& 	86.21\%& 	1	    &  -  & - &0.02\\
				{ Fisher-SVM}&	{\bf 85.65\%}&{\bf	86.34\%}&	10    &  10 & 2	&0.06\\
				RFE-SVM	  & 85.51\%& 	86.21\%& 	1	    &  1	& $2^{-1}$ &0.18\\
				FS-SVM&	85.51\%&	86.21\%&	1	    &  1	& 1 &0.06\\
				KS FS-SVM& 85.51\%& 	86.21\%& 	1		  &  1	& 1 &0.11\\
				\hline
			\end{tabular}%
		}
	\end{minipage}%
	\begin{minipage}{.5\linewidth}
		\centering
		\scalebox{0.65}{
			\begin{tabular}{rrrrrrr}
			\hline
			\multicolumn{7}{c}{German Credit m=1000 n=24}      \\
			\hline
			Form. & Av. ACC & Av. AUC & Av. F. & B  &C   &Time \\
			\hline
			{ $\ell_1$-SVM} & 76.20\% & 68.71\% & 24    &-       & $2^4$&0.04 \\
			{ $\ell_2$-SVM} & 76.10\% & 68.55\% & 24    &  -     &$2^4$ &0.18 \\
			{ LP-SVM}& 76.30\% & 68.79\% & 24    & -      &$2^{-3}$ &0.06 \\
			MILP1  & \textbf{76.90\%} & 68.93\% & 17    & 17 & -  & 0.54 \\
			MILP2 & 76.70\% & 69.26\% & 22    & 22    & $2^{-7}$&0.18 \\
			FSV   & 63.70\% & 68.36\% &  1     & -      &- &0.09 \\
			{ Fisher-SVM} & 75.10\% & \textbf{72.21\%} & 17    & 17 &$2^{-3}$   & 0.07 \\
			RFE-SVM	& 74.30\% & 71.93\% & 18    & 18  & $2^{-3}$ & 0.06 \\
			FS-SVM & 76.40\% & 68.95\% & 21    & 21   &$2^6$ & 0.32 \\
			KS FS-SVM& 76.50\% & 69.12\% & 21    & 21  &  &$2^5$ 0.36 \\
			\hline
			\end{tabular}%
		}
	\end{minipage} 
	\begin{minipage}[t]{.45\linewidth}
		\caption{\footnotesize{Best average ACC and AUC for Australian dataset. 	\label{tab_aus}}}
	\end{minipage}%
	\hfill%
	\begin{minipage}[t]{.45\linewidth}
		\caption{\footnotesize{Best average ACC and AUC for German credit dataset.	\label{tab_gc}}}
	\end{minipage}%
\end{table}

\begin{table}[!htb]
    \begin{minipage}{.5 \linewidth}
      \centering
			\scalebox{0.65}{
			\begin{tabular}{rrrrrrr} 
    			\hline
    \multicolumn{7}{c}{WBC m=569 n=30}         \\
    \hline
    Form. & Av. ACC & Av. AUC & Av. F. & B   &C  & Time \\
    \hline
    { $\ell_1$-SVM} & 97.37\% & 96.92\% & 10    &  -     &1 &0.01 \\
    { $\ell_2$-SVM} & \textbf{98.07\%} & \textbf{97.58\%} & 30    &   -  & $2^2$ & 0.09 \\
    {LP-SVM} & 97.89\% & 97.44\% & 30    &   -    &$2^{-5}$ &0.01 \\
    MILP1  & 97.02\% & 96.45\% & 3     & 3     &- &0.57 \\
    MILP2 & 97.89\% & 97.54\% & 23    & 23    &$2^{-5}$ &0.19 \\
    FSV   & 42.35\% & 54.03\% & 20    &  -     &- &0.02 \\
    { Fisher-SVM} & 95.60\% & 96.19\% & 30    & 30    &$2^6$ &0.17 \\
    RFE-SVM & 95.78\% & 96.33\% & 23    & 23    & $2^6$ &0.20 \\
    FS-SVM   & 97.72\% & 97.20\% & 4     & 4     & $2^4$ &1.04 \\
    KS FS-SVM & 97.72\% & 97.20\% & 4     & 4     &$2^4$  &0.90 \\
    \hline
    \end{tabular}%
		}
    \end{minipage}%
    \begin{minipage}{.5 \linewidth}
		 \centering
				\scalebox{0.65}{
    \begin{tabular}{rrrrrrr}
    \hline
    \multicolumn{7}{c}{Ionosphere m=351 n=33}       \\
    \hline
    Form. & Av. ACC   & Av. AUC   & Av. F.    & B& C&Time \\
    \hline
    { $\ell_1$-SVM} & 88.33\% & 85.39\% & 28    &  -     &2 &0.01 \\
    { $\ell_2$-SVM} & 86.76\% & 83.92\% & 33    &  -     &$2^3$ &0.05 \\
    {LP-SVM}& 86.20\% & 83.82\% & 33    &     -  & 1&0.01 \\
    MILP1  & 88.06\% & 85.62\% & 16    & 16    &- &4.98 \\
    MILP2 & 88.52\% & 85.88\% & 6     & 6     &$2^{-7}$  &0.15 \\
    FSV   & 67.22\% & 54.23\% & 5     &    -   &- &0.02 \\
    { Fisher-SVM}  & 88.24\% & 84.05\% & 28    & 28    &$2^{-2}$ &0.05 \\
    RFE-SVM	& \textbf{88.61\%} & 84.44\% & 2     & 2     & $2^{-4}$&0.07 \\
    FS-SVM  & \textbf{88.61\%} & \textbf{86.39\%} & 16    & 16    &$2^5$ &6.39 \\
    KS FS-SVM & \textbf{88.61\%} & \textbf{86.39\%} & 16    & 16    &$2^5$ &5.61 \\
    \hline
    \end{tabular}%
		}
    \end{minipage} 

		\begin{minipage}[t]{.45\linewidth}
        \caption{\footnotesize{\footnotesize{Best average ACC and AUC for WBC dataset.} 	\label{tab_wbc}}}
    \end{minipage}%
		\hfill%
    \begin{minipage}[t]{.45\linewidth}
        \caption{\footnotesize{Best average ACC and AUC for Ionosphere dataset.} \label{tab_iono}}
    \end{minipage}%
\end{table}

\begin{figure}[H]
	\begin{minipage}{.5 \linewidth}
		\centering
		\includegraphics[width=0.9\textwidth]{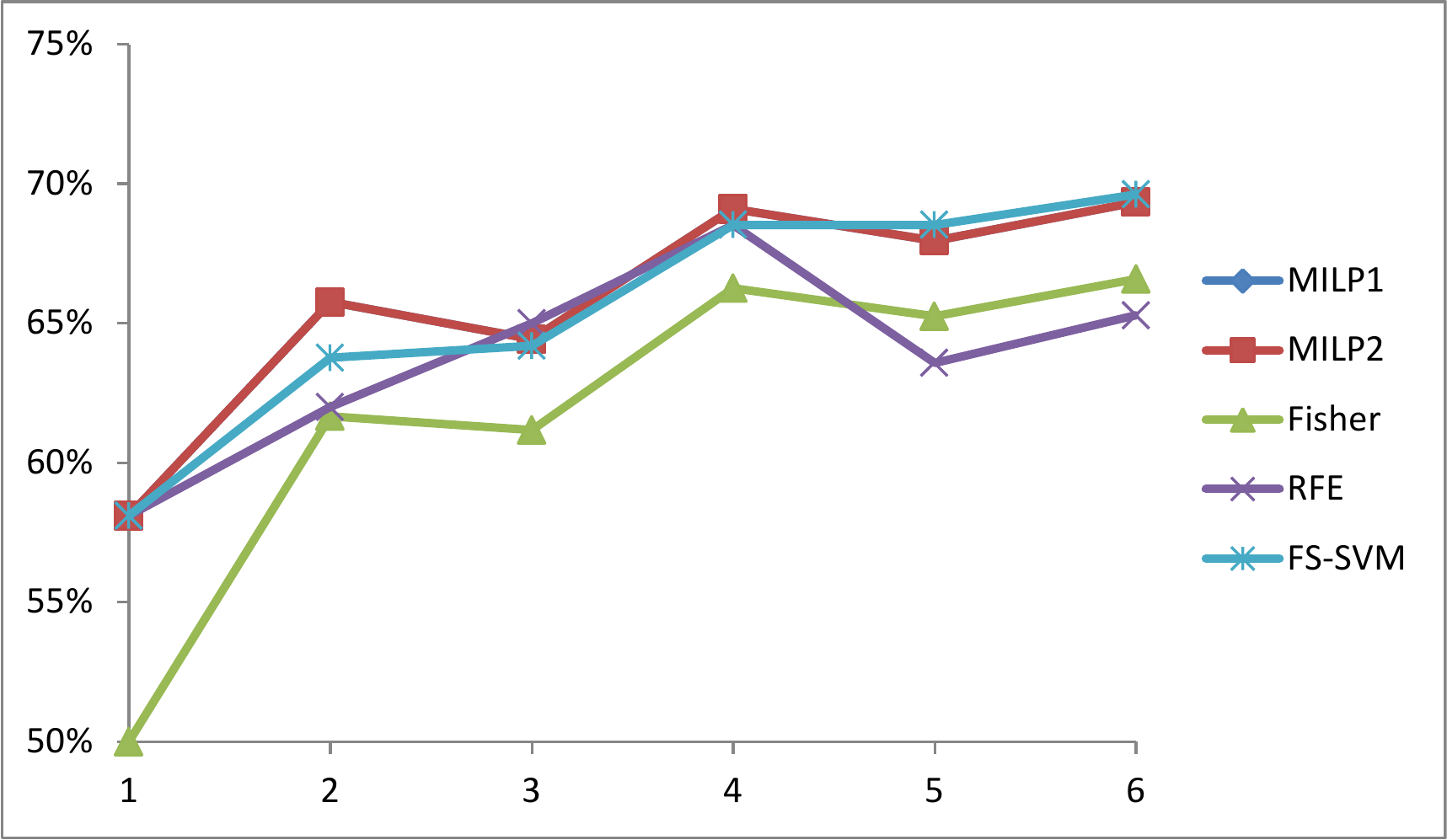}
	\end{minipage}%
	\begin{minipage}{.5 \linewidth}
		\centering
		\includegraphics[width=0.9\textwidth]{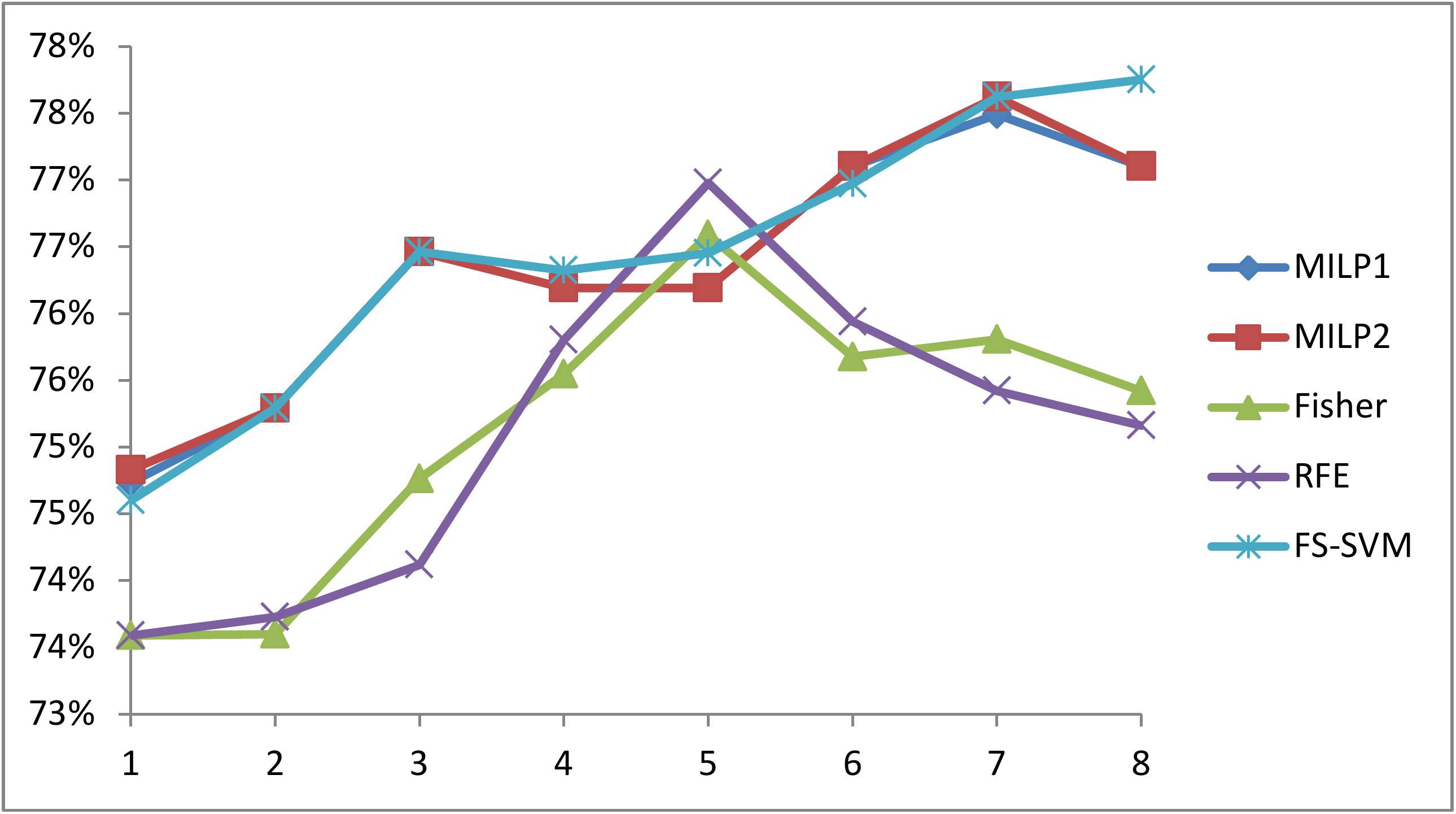}
	\end{minipage} 
	\begin{minipage}[t]{.5\linewidth}
		\caption{\footnotesize{\footnotesize{Average ACC for BUPA dataset.} 		\label{fig:bupa_acc}}}
	\end{minipage}%
	\hfill%
	\begin{minipage}[t]{.5\linewidth}
		\caption{{\footnotesize{Average ACC for PIMA dataset.}\label{fig:pima_acc}}}
	\end{minipage}%
\end{figure}
\

\begin{figure}[H]
	\begin{minipage}{.5 \linewidth}
		\centering
		\includegraphics[width=0.9\textwidth]{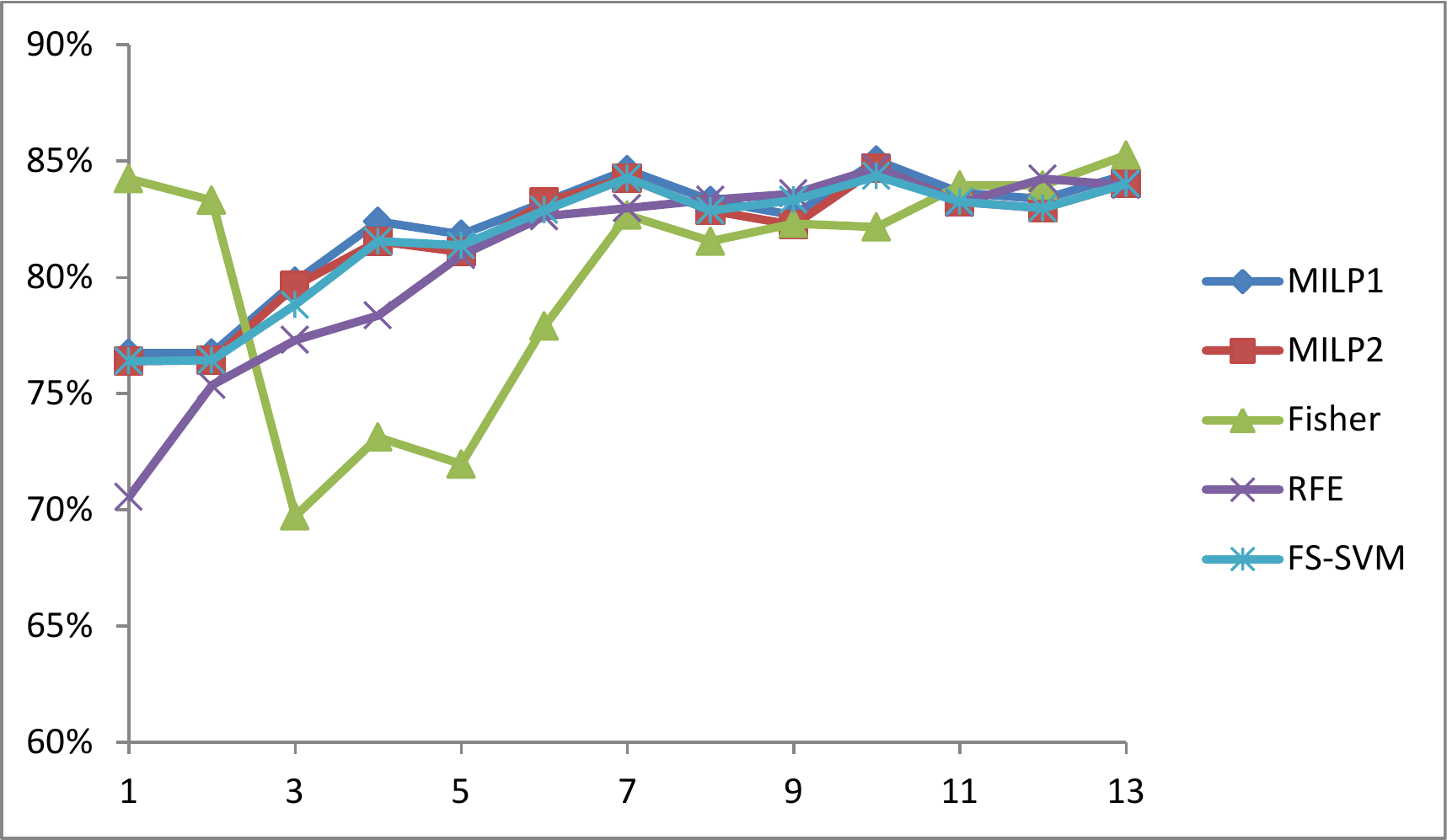}
	\end{minipage}%
	\begin{minipage}{.5 \linewidth}
		\centering
		\includegraphics[width=0.9\textwidth]{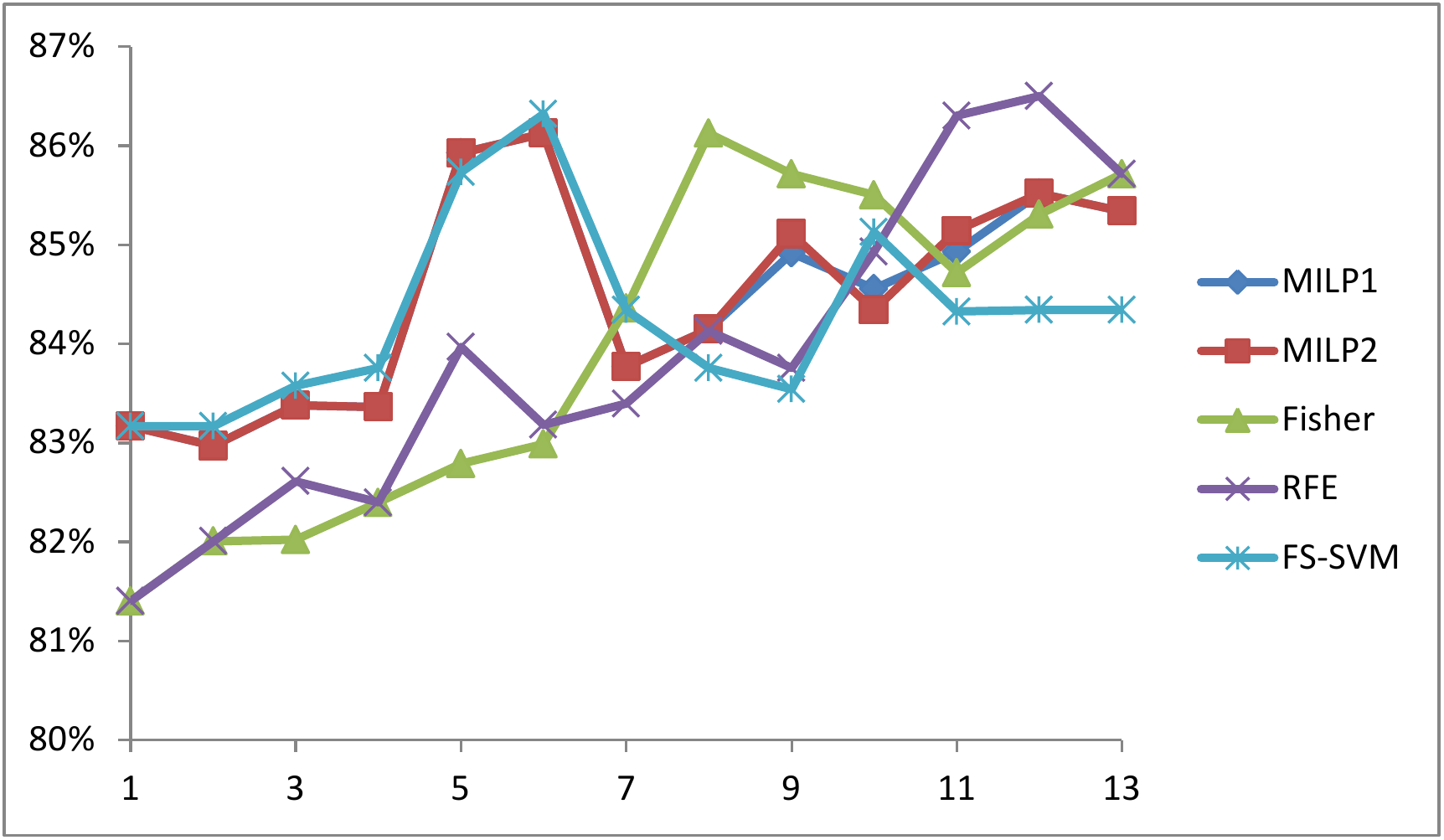}
	\end{minipage} 
	\begin{minipage}[t]{.5\linewidth}
		\caption{\footnotesize{\footnotesize{Average ACC for Cleveland dataset.} 		\label{fig:clev_acc}}}
	\end{minipage}%
	\hfill%
	\begin{minipage}[t]{.5\linewidth}
		\caption{\footnotesize{Average ACC for Housing dataset.}\label{fig:hous_acc}}
	\end{minipage}%
\end{figure}

\begin{figure}[H]
    \begin{minipage}{.5 \linewidth}
				\centering
					\includegraphics[width=0.9\textwidth]{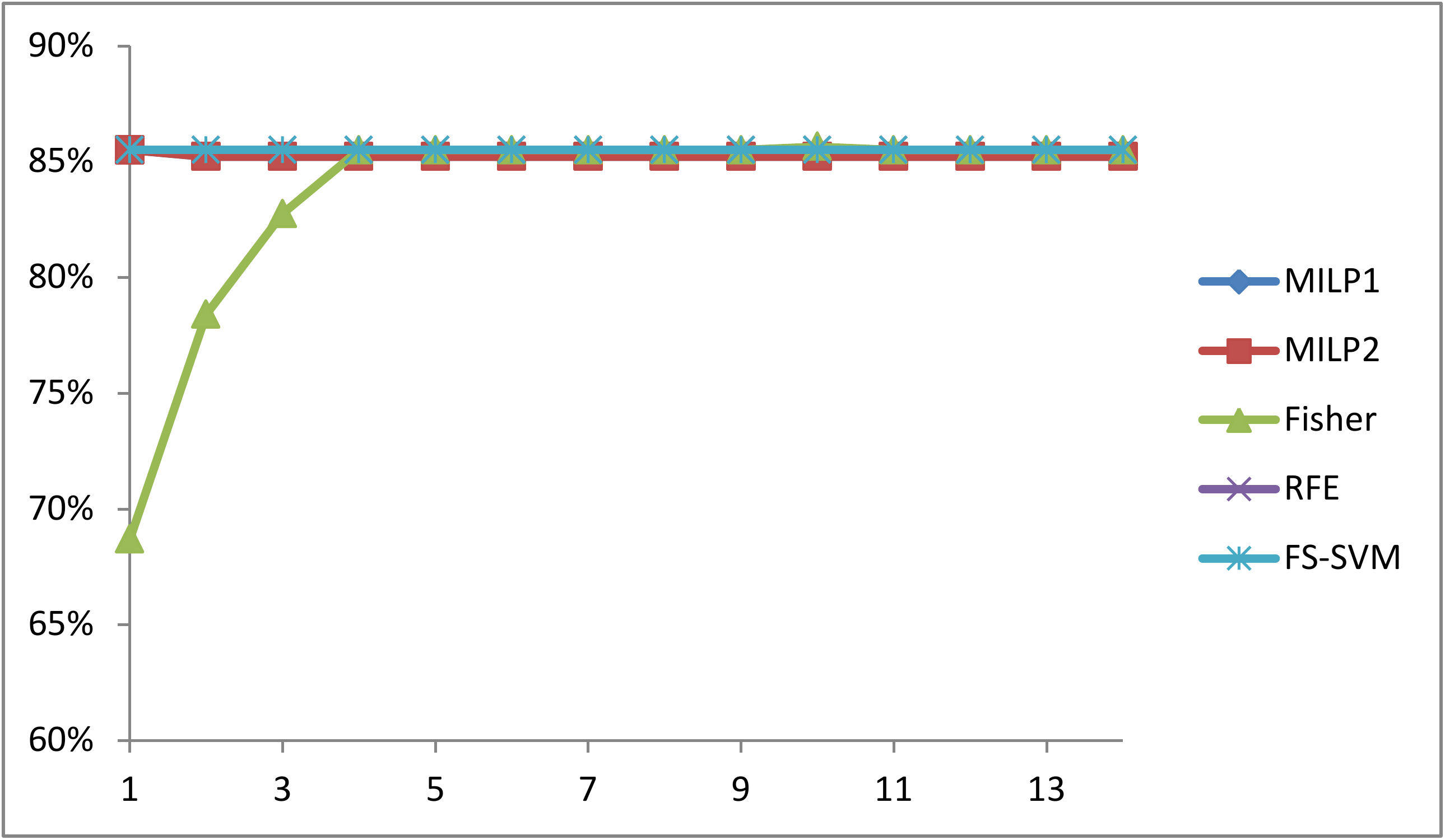}
    \end{minipage}%
    \begin{minipage}{.5 \linewidth}
     \centering
			\includegraphics[width=0.9\textwidth]{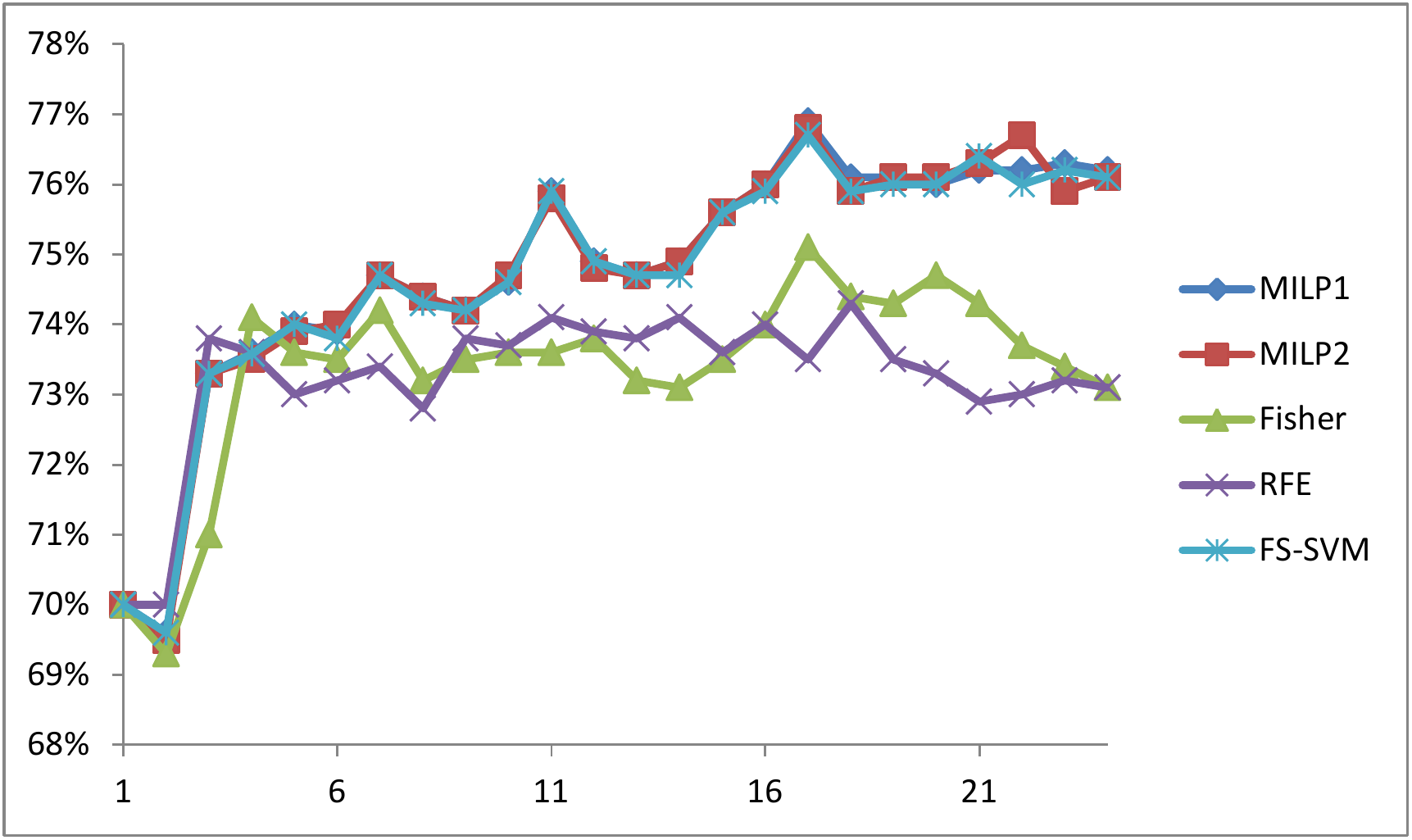}
    \end{minipage} 
		\begin{minipage}[t]{.5\linewidth}
        \caption{\footnotesize{\footnotesize{Average ACC for Australian dataset.} 		\label{fig:australian_acc}}}
    \end{minipage}%
		\hfill%
    \begin{minipage}[t]{.5\linewidth}
        \caption{\footnotesize{Average ACC for G. Credit dataset.} \label{fig:gc_acc}}
    \end{minipage}%
\end{figure}

\begin{figure}[H]
    \begin{minipage}{.5 \linewidth}
				\centering
					\includegraphics[width=0.9\textwidth]{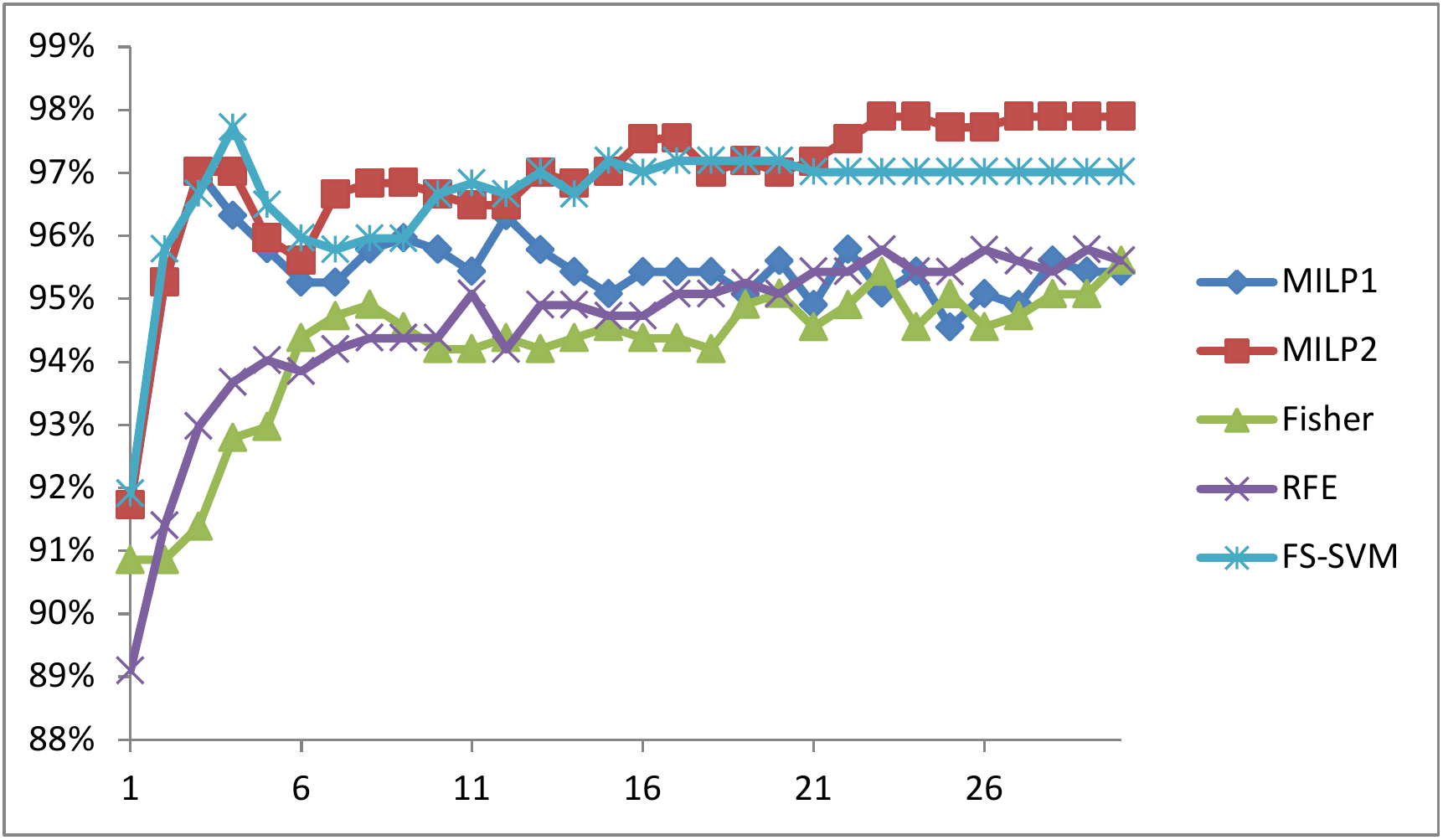}
    \end{minipage}%
    \begin{minipage}{.5 \linewidth}
		 \centering
			\includegraphics[width=0.9\textwidth]{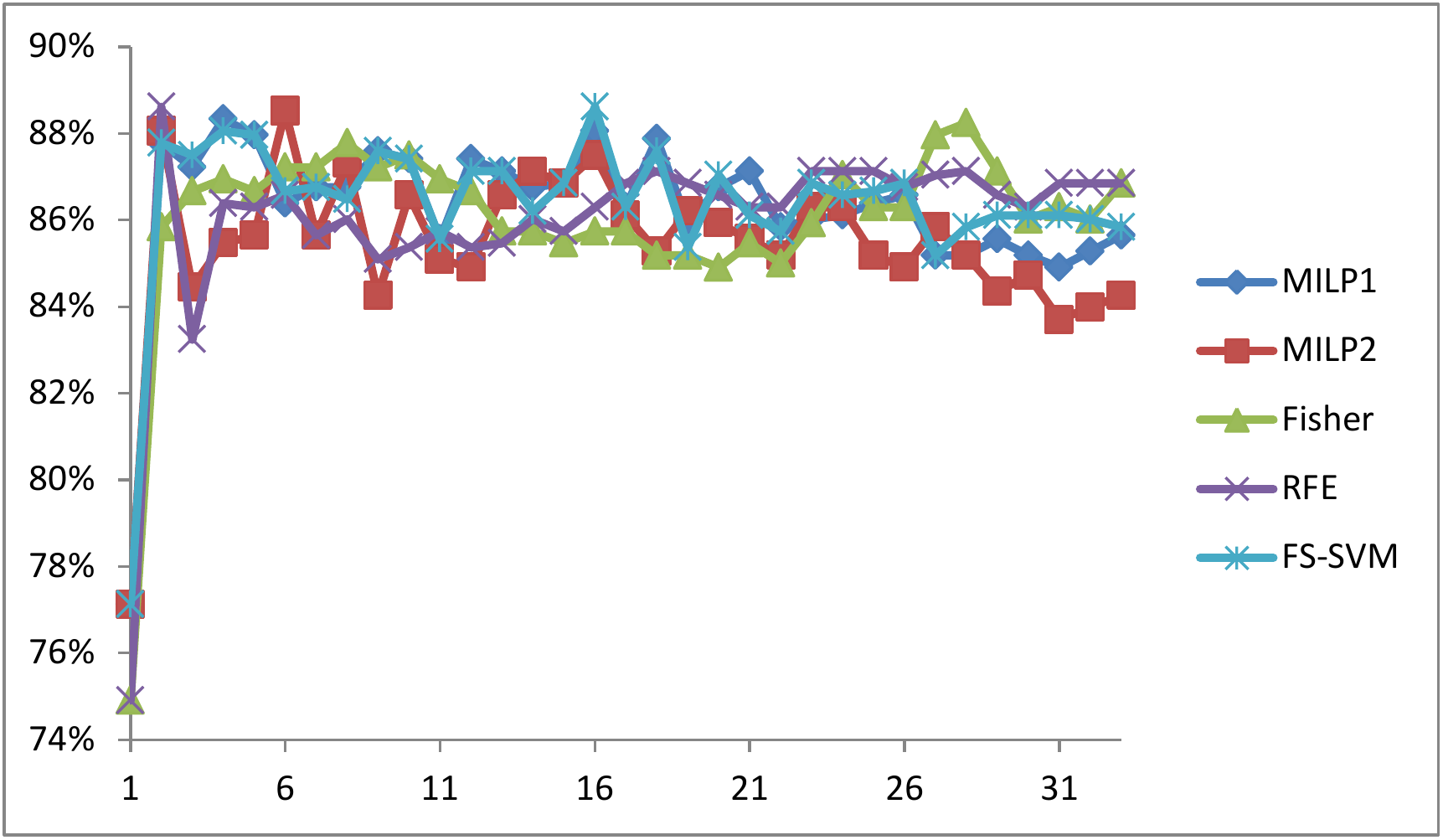}
    \end{minipage} 
		\begin{minipage}[t]{.5\linewidth}
        \caption{\footnotesize{\footnotesize{Average ACC for WBC dataset.}} \label{fig:wbc_acc}}
    \end{minipage}%
		\hfill%
    \begin{minipage}[t]{.5\linewidth}
        \caption{\footnotesize{Average ACC for IONO dataset.}\label{fig:iono_acc}}
    \end{minipage}%
\end{figure}
\subsection{Instances with small sample size and big number of features}\label{BigInstances}
 This subsection is focused on Colon, Leukemia, DLBCL and Carcinoma datasets. Tables \ref{tab_colon}-\ref{tab_car} do not include the results of MILP2 due to the huge times to solve this model for big instances. Since the KS provides good approximations of the solutions for our model within short time, we will only use the KS to analyze FS-SVM performance for most of the datasets. In addition, the performance of ACC is illustrated in Figures \ref{colon_acc}-\ref{car_acc} in which $B$ varies betweeen $10$ and $100$ for the $C$ values shown in Tables \ref{tab_colon}-\ref{tab_car}.

The results for the Colon dataset are indicated in Table \ref{tab_colon}. We observe that the best performance are given by RFE-SVM and FS-SVM. Moreover, Figure \ref{colon_acc} shows that FS-SVM
achieves stable results with all the possible values of $B$. In fact, RFE-SVM, Fisher-SVM and FS-SVM present similar performances.

\begin{table}[!htb]
	\begin{minipage}{.5 \linewidth}
		\centering
		\scalebox{0.6}{
			\begin{tabular}{rrrrrrr}
				\hline
				\multicolumn{7}{c}{Colon m=62 n=2000}       \\
				\hline
				Form. & Av. ACC   & Av. AUC   & Av. F.    & B&&Time \\
				\hline
				{ $\ell_1$-SVM} & 90.42\% & 88.75\% & 10   &   -    &$2^{-2}$ &0.38 \\
				{ $\ell_2$-SVM} & 88.75\% & 88.75\% & 2000  &  -     &$2^{-6}$ &0.05 \\
				{ LP-SVM} & 87.08\% & 86.25\% & 2000  &  -     & -&0.11 \\
				MILP1  & 85.83\% & 85.00\% & 10    & 10    &- &41.16 \\
				FSV   & 62.92\% & 60.00\% &  10     &  -     &- &1.76 \\
				{ Fisher-SVM} & 90.42\% & 90.00\% & 100   & 100   &$2^{-1}$ &0.00 \\
				RFE-SVM   & \textbf{92.08\%} & \textbf{91.25\%} & 20    & 20    &$2^{-4}$ &0.00 \\
				FS-SVM  & \textbf{92.08\%} & \textbf{ 91.25\%} & 20    & 20    & $2^{7}$&678.78 \\
				KS FS-SVM &  90.42\% &  88.75\% & 20    & 20    &$2^{7}$ &9.51 \\
				\hline
			\end{tabular}%
		}
	\end{minipage}%
	\begin{minipage}{.5 \linewidth}
		\centering
		\scalebox{0.6}{
			\begin{tabular}{rrrrrrr}
				\hline
				\multicolumn{7}{c}{Leukemia m=72 n=5327}      \\
				\hline
				Form. & Av. ACC & Av. AUC & Av. F. & B     & C&Time \\
				\hline
				{ $\ell_1$-SVM} & { 97.14\%} & { 97.50\%} & 38    &   -    &$2^{-1}$ &1.72 \\
				{ $\ell_2$-SVM} & \textbf{98.57}\% & \textbf{98.75}\% & 5327  &  -     & 1&0.11 \\
				{ LP-SVM}& 96.03\% & 96.37\% & 5327  &  -     & 1&0.28 \\
				MILP1  & 80.63\% & 78.87\% & 10    &  10     &- &374.85 \\
				FSV   & 60.79\% & 65.00\% & 0.2   &    -   &- &10.04 \\
					{ Fisher-SVM} & { 97.14\%} & { 97.50\%} & { 60} & { 60}  & { $2^{-5}$} & { 0.01} \\
				RFE-SVM  & { 97.14\%} & { 97.50\%} & { 30}   & { 30}   & { $2^{-5}$}& { 0.01} \\
				KS FS-SVM & 97.14\% & 97.50\% & 30    & 30    & $2^{-1}$&13.72 \\
				\hline
			\end{tabular}%
		}
	\end{minipage} 
	\begin{minipage}[t]{.45\linewidth}
		\caption{\footnotesize{\footnotesize{Best average ACC and AUC for Colon dataset.}\label{tab_colon}}}
	\end{minipage}%
	\hfill%
	\begin{minipage}[t]{.45\linewidth}
		\caption{\footnotesize{Best average ACC and AUC for Leukemia dataset.}\label{tab_leukemia}}
	\end{minipage}%
\end{table}

\begin{table}[!htb]
	\begin{minipage}{.5 \linewidth}
		\centering
		\scalebox{0.6}{
			\begin{tabular}{rrrrrrr}
				\hline
				\multicolumn{7}{c}{DLBCL m=77 n=7129}         \\
				\hline
				Form. & Av. ACC & Av. AUC & Av. F. & B    &C & Time \\
				\hline
				{ $\ell_1$-SVM} & {\textbf{ 98.75\%}} & { 97.50\%} & 32    &   -  & 1 & 3.15 \\
				{ $\ell_2$-SVM} & 96.25\% & 94.17\% & 7129  &  - & 1   & 0.15 \\
				{ LP-SVM} & 97.50\% & 95.00\% & 7129  &  -   &$2^{-4}$  & 0.50 \\
				MILP1  & 88.75\% & 85.83\% & 10    & 10   &- & 728.08 \\
				FSV   & 49.75\% & 66.67\% & 23      &  -   & - & 15.66 \\
				{ Fisher-SVM} & { 83.50\%} & { 88.17\%} & { 100}   & { 100 }& { $2^{-5}$} & { 0.01} \\
				RFE-SVM  & \textbf{98.75\%} & \textbf{99.17\%} & 40    & 40 & $2^7$  & 0.01 \\
				KS FS-SVM & \textbf{98.75\%} & 97.50\% & 30    & 30    &$2^6$ &1.41 \\
				\hline
			\end{tabular}%
		}
	\end{minipage}%
	\begin{minipage}{.5 \linewidth}
		\centering
		\scalebox{0.6}{
			\begin{tabular}{rrrrrrr}
				\hline
				\multicolumn{7}{c}{Carcinoma m=36 n=7457} \\
				\hline
				Form. & ACC   & AUC   & Av. F. & B     & C     & Time \\
				\hline
				{ $\ell_1$-SVM} & 92.50\% & 87.50\% & 14.6  & -      & $2^{-2}$ & 0.81 \\
				{ $\ell_2$-SVM} & \textbf{95.00\%} & \textbf{90.00\%} & 7457  &   -    & $2^7 $   & 0.06 \\
				{ LP-SVM} & 92.50\% & 87.50\% & 7457  &-       & $2^7$   & 0.17 \\
				MILP1  & 74.17\% & 70.00\% & 90    & 90    &   -    & 278.92 \\
				FSV   & 74.17\% & 67.50\% & 2.9   & -      &    -   & 5.42 \\
				{ Fisher-SVM} & { \textbf{95.00\%}} & { \textbf{90.00\%}} & { 30}   & { 30}   & { $2^{-4}$}   & { 0.01} \\
				{ RFE-SVM}  & \textbf{95.00\%} & \textbf{90.00\%} & 20    & 20    & $2^7$   & 2.28 \\
				KS-FS-SVM & 92.50\% & 87.50\% & 10    & 10    & $2^{-4}$ & 0.17 \\
				\hline			
			\end{tabular}%
		}
	\end{minipage} 
	\begin{minipage}[t]{.45\linewidth}
		\caption{\footnotesize{\footnotesize{Best average ACC and AUC for DLBCL dataset.}\label{tab_DLBCL}}}
	\end{minipage}%
	\hfill%
	\begin{minipage}[t]{.45\linewidth}
		\caption{\footnotesize{Best average ACC and AUC for Carcinoma dataset.}\label{tab_car}}
	\end{minipage}%
\end{table}

\begin{figure}[!htb]
	\begin{minipage}{.5 \linewidth}
		\centering
		\includegraphics[width=0.9\textwidth]{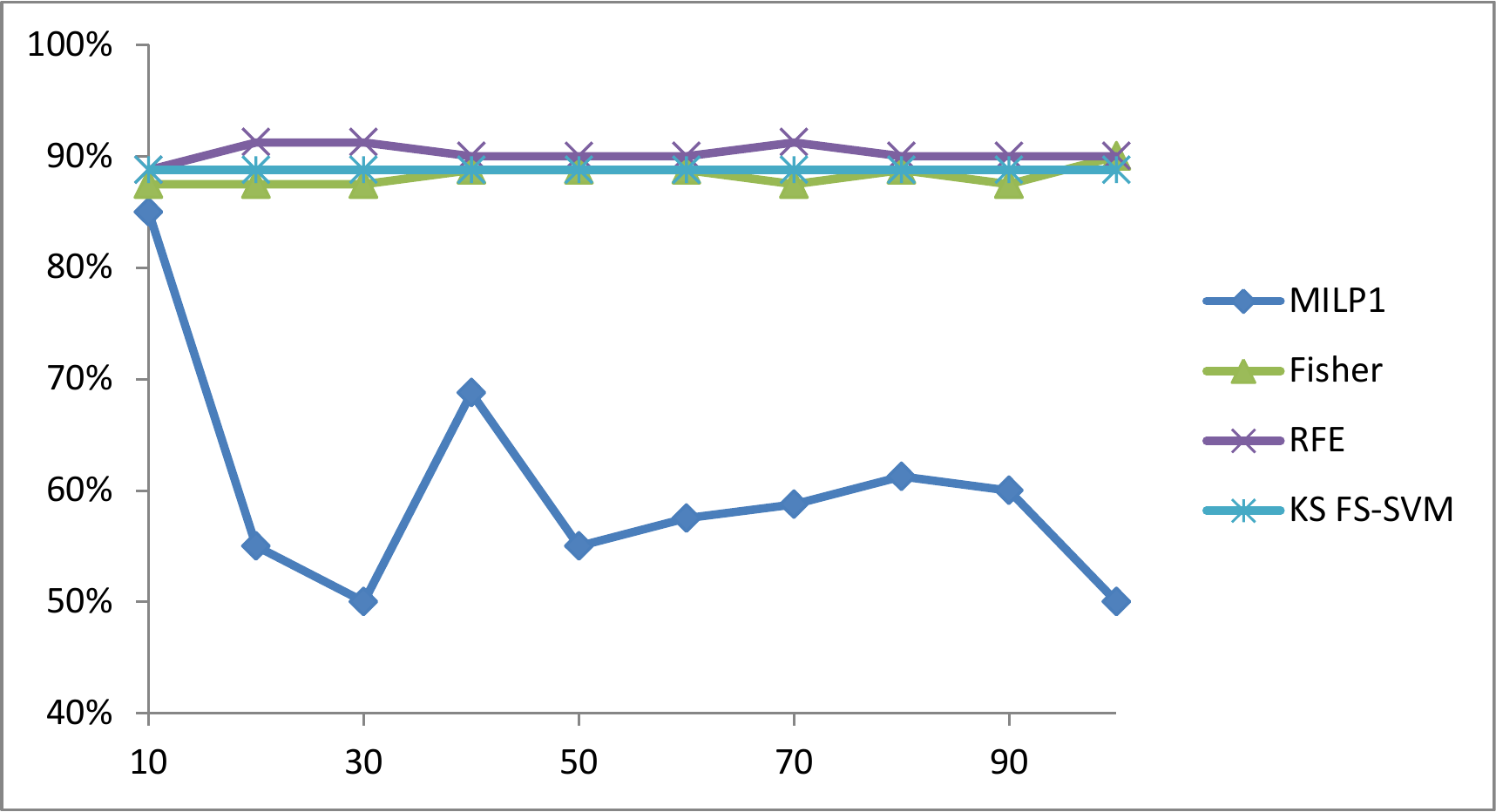}
	\end{minipage}%
	\begin{minipage}{.5 \linewidth}
		\centering
		\includegraphics[width=0.9\textwidth]{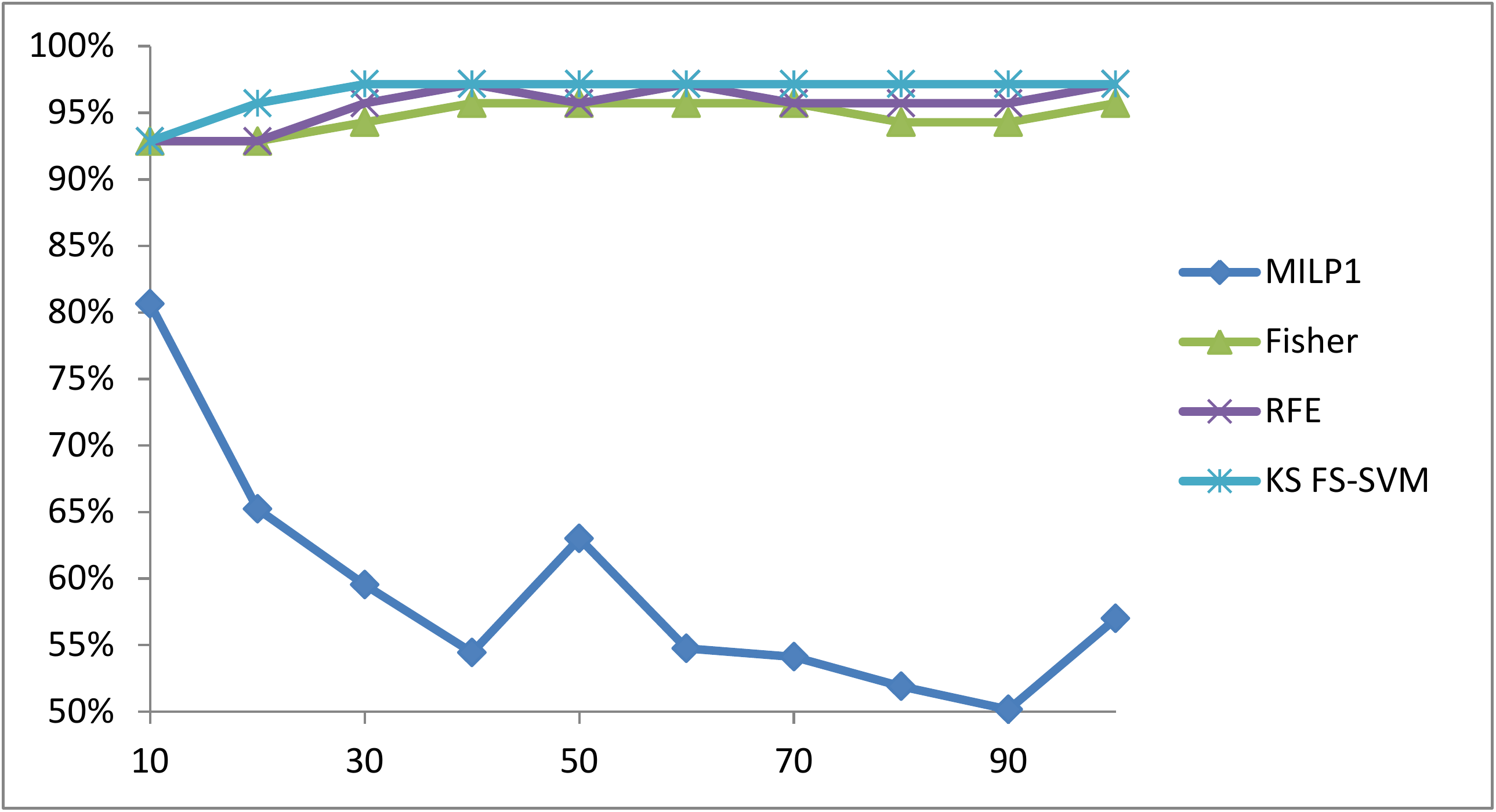}
	\end{minipage} 
	\begin{minipage}[t]{.5\linewidth}
		\caption{\footnotesize{\footnotesize{Average ACC for Colon dataset.}}\label{colon_acc}}
	\end{minipage}%
	\hfill%
	\begin{minipage}[t]{.5\linewidth}
		\caption{\footnotesize{Average ACC for Leukemia dataset.}\label{fig:leukemia_acc}}
	\end{minipage}%
\end{figure}

\begin{figure}[!htb]
	\begin{minipage}{.5 \linewidth}
		\centering
		\includegraphics[width=0.9\textwidth]{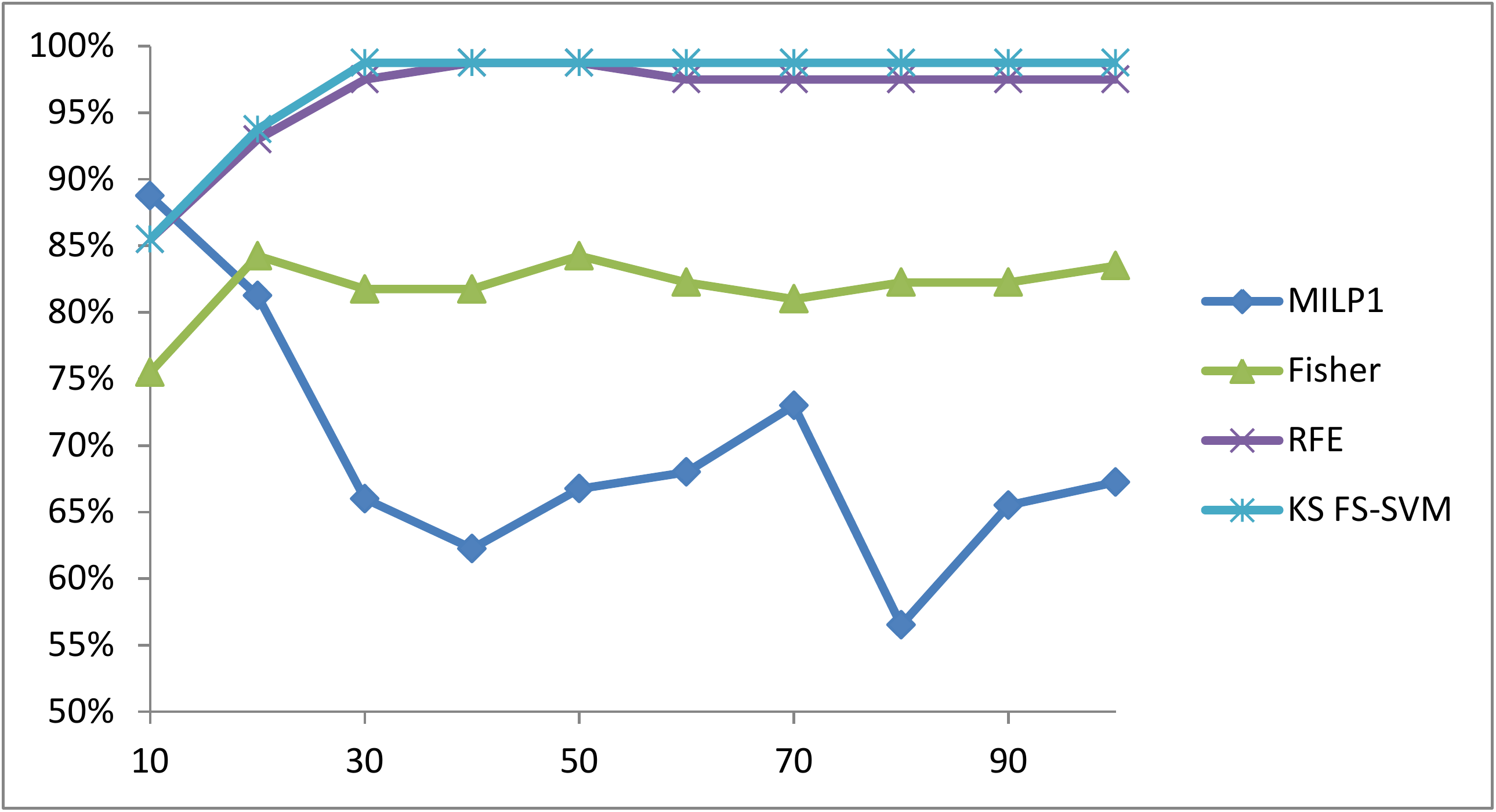}
	\end{minipage}%
	\begin{minipage}{.5 \linewidth}
		\centering
		\includegraphics[width=0.9\textwidth]{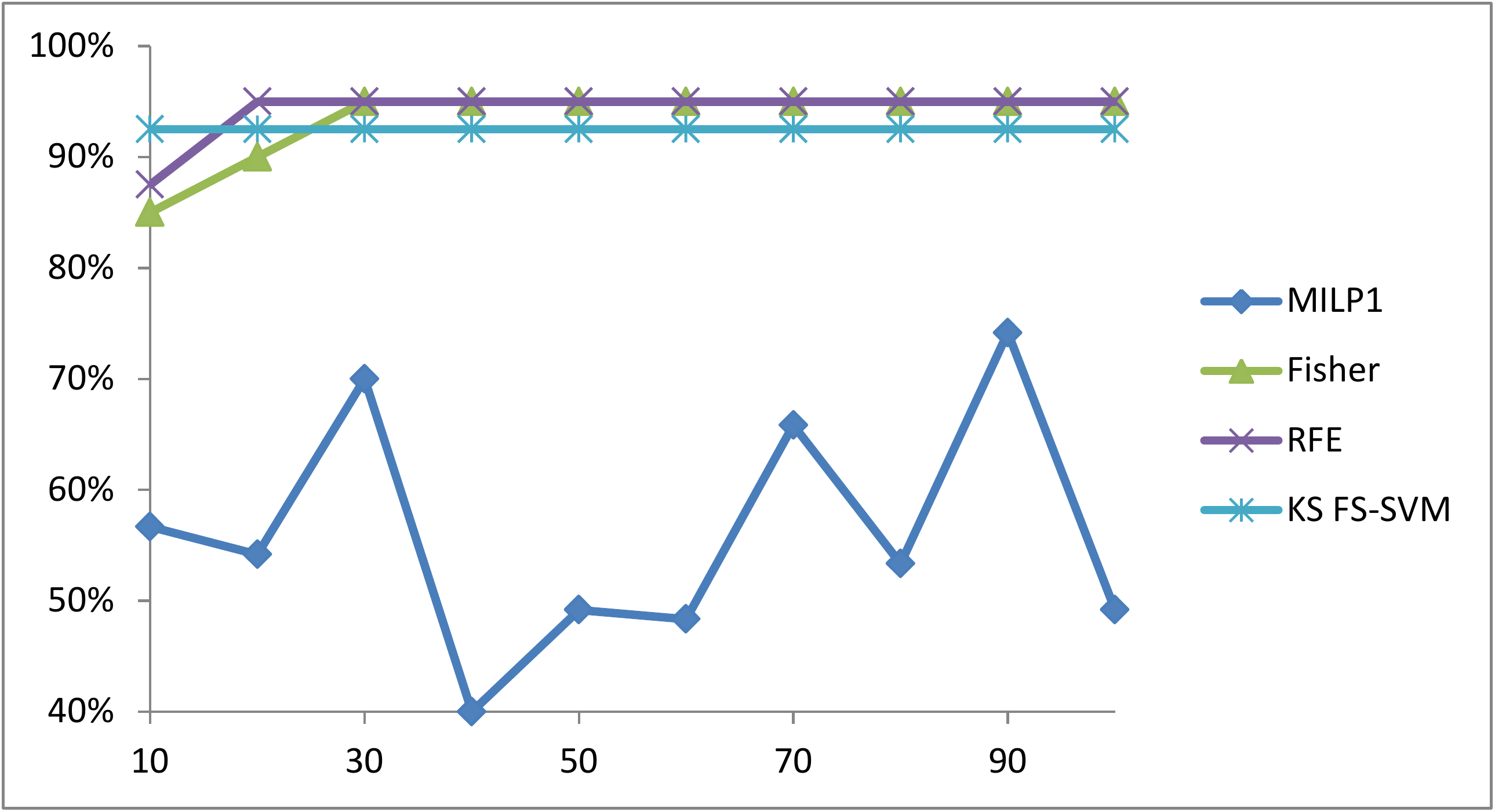}		
	\end{minipage} 
	\begin{minipage}[t]{.5\linewidth}
		\caption{\footnotesize{\footnotesize{Average ACC for DLBCL dataset.}}\label{fig:dlbcl_acc}}
	\end{minipage}%
	\hfill%
	\begin{minipage}[t]{.5\linewidth}
		\caption{\footnotesize{Average ACC for Carcinoma dataset.}\label{car_acc}}
	\end{minipage}%
\end{figure}

For the Leukemia dataset, the best performance is achieved in $\ell_2$-SVM, but FS-SVM also attains a good performance using only 30 features in contrast 
with the 5327 used by $\ell_2$-SVM. Besides, Figure
\ref{fig:leukemia_acc} shows that KS FS-SVM provides the best performance.
 For the DLBCL dataset, we obtain similar results to the Leukemia case and the results are shown in Table \ref{tab_DLBCL} and Figure \ref{fig:dlbcl_acc}.
The best performances in this case are achieved by the $\ell_1$-SVM, the RFE-SVM and KS FS-SVM models.

Lastly, for the Carcinoma dataset, we can see that the best average ACC and AUC are obtained with $\ell_2$-SVM, Fisher-SVM and RFE-SVM. In addition, Figure \ref{car_acc} shows that KS FS-SVM, Fisher-SVM and RFE-SVM provide good results for all studied $B$ values.

To summarize, the performances of the different models in these four datasets show that our model provides results which are better than MILP1  models and that it attains similar results to RFE-SVM. Regarding these results, we can therefore conclude that our model is
a good classifier for big  instances where a small number of features must be selected.

 \subsection{Instances with big sample size and big number of features}\label{SuperBigInstances}
 This subsection analyzes Lepiota, Arrythmia, Madelon and MFeat datasets. As in the previous section, Tables \ref{tab_lepiota}-\ref{tab_mfeat} do not include the results of MILP2 due to the huge performing times of this model. For the same reason, Tables \ref{tab_arrythmia} and \ref{tab_madelon} do not show the results of MILP1. Besides, we will use the KS to analyze FS-SVM performance. The results of ACC is illustrated in Figures \ref{Lepiota_acc}-\ref{MFeat_acc} in which $B$ varies betweeen $10$ and $100$ for the $C$ values shown in Tables \ref{tab_lepiota}-\ref{tab_mfeat}.

 The results of the Lepiota dataset are presented in Table \ref{tab_lepiota}. Note that for this dataset almost all models provide the same good results. It is important to point out that KS FS-SVM is the model which attains the best accuracy using the smallest number of features. In Figure \ref{Lepiota_acc} we can see that for each value of parameter $B$, KS FS-SVM provides the best average accuracy. In the case of the Arrythmia dataset (Table \ref{tab_arrythmia}), the model with the best behaviour is $\ell_2$-SVM. However, we would highlight the fact that  KS FS-SVM obtains a similar result whilst using  many less features than $\ell_2$-SVM.  Figure \ref{arrythmia_acc} shows that KS-FS-SVM provides the best results.
	
Fisher-SVM is the best model  for Madelon dataset, but $\ell_1$-SVM and KS FS-SVM also provide good accuracy and AUC  with a lesser number of features, (see Table \ref{tab_madelon}). In addition to this, Figure \ref{Madelon_acc} illustrates that KS FS-SVM provides the best behaviour for different values of $B$. For MFeat dataset, the best results are achieved by KS FS-SVM, as can be seen in \ref{tab_mfeat}. 
Moreover, it can be seen in Figure \ref{MFeat_acc}, the best accuracy is obtained when using KS FS-SVM for different values of $B$.

\begin{table}[!htb]
		\begin{minipage}{.5 \linewidth}
			\centering
			\scalebox{0.6}{
				\begin{tabular}{rrrrrrr}
					\hline
					\multicolumn{7}{c}{Lepiota m=1824 n=109}       \\
					\hline
					Form.&Av. ACC&Av. AUC & Av. F.& B&C&Time\\
					\hline
					$\ell_1$-SVM & \textbf{100.00\%} & \textbf{100.00\%}
					& 12.7  &   &  2$^6$     & 4.28 \\
					$\ell_2$-SVM &\textbf{ 100.00\%} &\textbf{  100.00\%} & 110   & -  &   2$^6$    & 135.98 \\
					LP-SVM & \textbf{ 100.00\%} & \textbf{ 100.00\%} & 114.8 &  -&  2$^6$      & 1.38 \\
					MILP1  & \textbf{ 100.00\%} &\textbf{  100.00\%} &   10    &   10    & -  & 329.58 \\
					FSV   & 71.70\% & 71.88\% & 2     &   -    &   -    & 20.13 \\
					Fisher-SVM & 99.90\% & 99.90\% & 95    & 100  &  2$^{-7}$ & 7.73 \\
					RFE-SVM   & 99.91\% & 99.91\% & 90    &90   &  2$^{-7}$  & 7.33 \\
					KS FS-SVM & \textbf{ 100.00\%} & \textbf{ 100.00\%} &    10  & 10  & 2$^6$   & 0.67 \\
					\hline
				\end{tabular}%
			}
		\end{minipage}%
		\begin{minipage}{.5 \linewidth}
			\centering
			\scalebox{0.6}{
				\begin{tabular}{rrrrrrr}
					\hline
					\multicolumn{7}{c}{Arrythmia m=420 n=258}      \\
					\hline
					Form. & Av. ACC & Av. AUC & Av. F. & B     & C&Time \\
					\hline
					$\ell_1$-SVM & 77.86\% & 75.93\% & 69.2  &       & 2\^-1 & 0.76 \\
					$\ell_2$-SVM & \textbf{78.81\%} & \textbf{77.20\%} & 258   &       & 2$^{-3}$ & 0.13 \\
					LP-SVM & 70.48\% & 69.21\% & 278   &       & 2$^{-7}$ & 0.16 \\
					FSV   & 64.76\% & 62.08\% & 117.6 &       &       & 0.23 \\
					Fisher-SVM & 72.14\% & 68.22\% & 100   & 100   & 2$^{-4}$ & 0.32\\
					RFE-SVM & 74.76\% & 71.27\% & 80    & 80    & 2$^{-5}$ & 13.29\\
					KS FS-SVM & 77.14\% & 75.47\% & 40    & 40    & 2$^0$  & 1851.09 \\
					\hline
				\end{tabular}%
			}
		\end{minipage} 
		\begin{minipage}[t]{.45\linewidth}
			\caption{\footnotesize{\footnotesize{Best average ACC and AUC for Lepiota dataset.}\label{tab_lepiota}}}
		\end{minipage}%
		\hfill%
		\begin{minipage}[t]{.45\linewidth}
			\caption{ \footnotesize{Best average ACC and AUC for Arrythmia dataset.}\label{tab_arrythmia}}
		\end{minipage}%
\end{table}

\begin{figure}[!htb]
	\begin{minipage}{.5 \linewidth}
		\centering
		\includegraphics[width=0.9\textwidth]{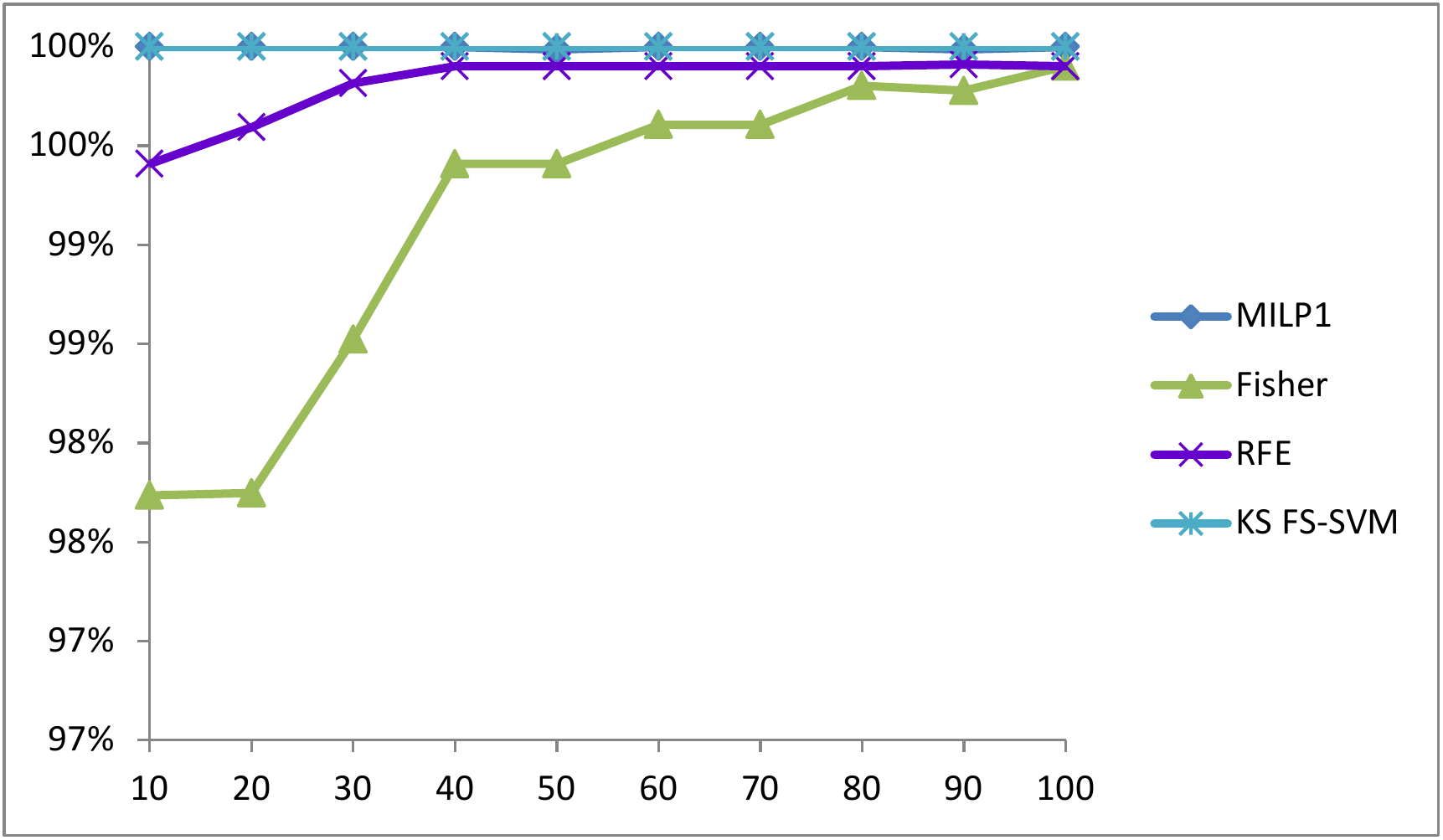}
	\end{minipage}%
	\begin{minipage}{.5 \linewidth}
		\centering
		\includegraphics[width=0.9\textwidth]{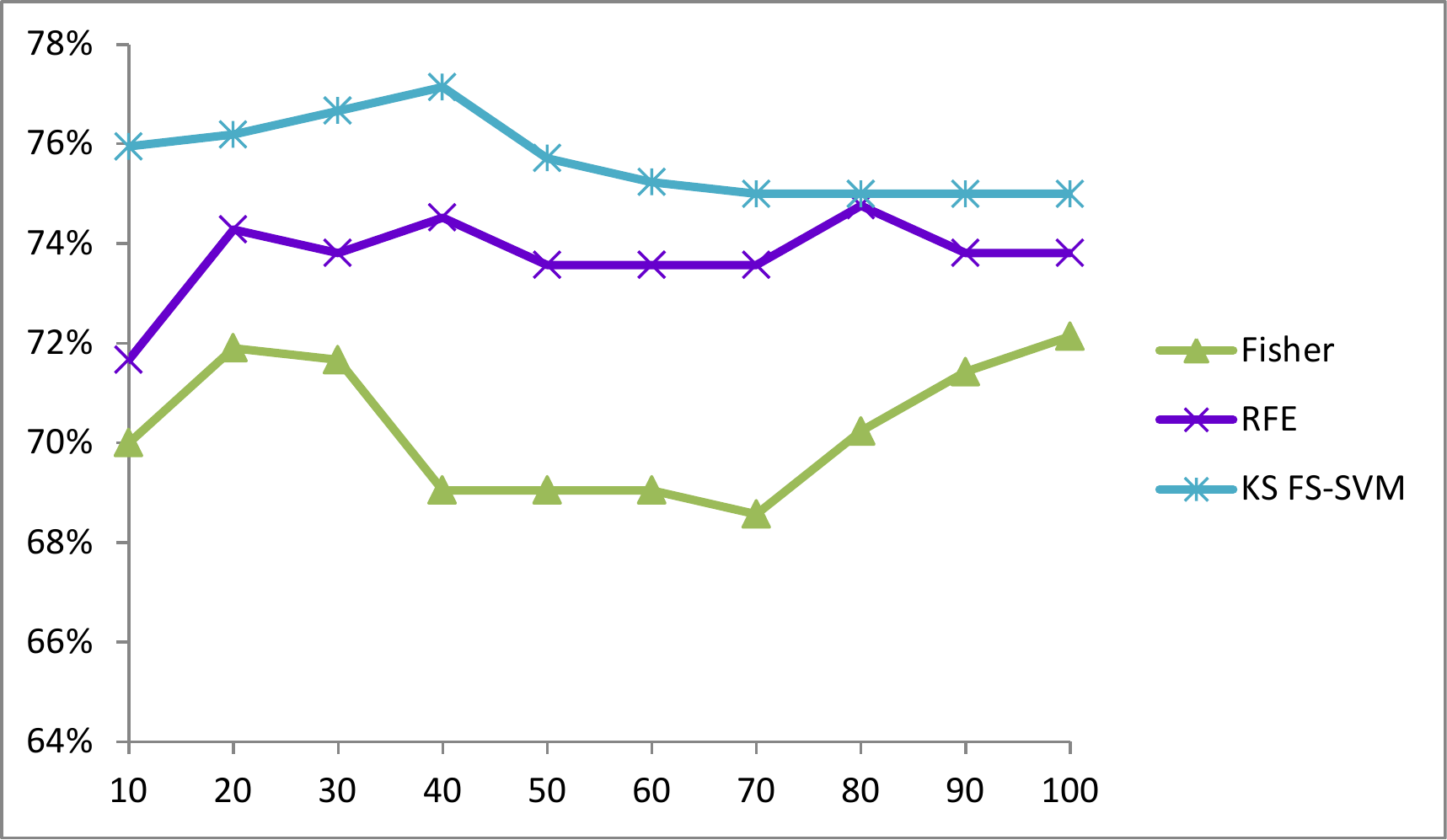}
	\end{minipage} 
	\begin{minipage}[t]{.5\linewidth}
		\caption{\footnotesize{\footnotesize{Average ACC for Lepiota dataset.}}\label{Lepiota_acc}}
	\end{minipage}%
	\hfill%
	\begin{minipage}[t]{.5\linewidth}
		\caption{\footnotesize{Average ACC for Arrythmia dataset.}\label{arrythmia_acc}}
	\end{minipage}%
\end{figure}

\begin{figure}[!htb]
	\begin{minipage}{.5 \linewidth}
		\centering
		\includegraphics[width=0.9\textwidth]{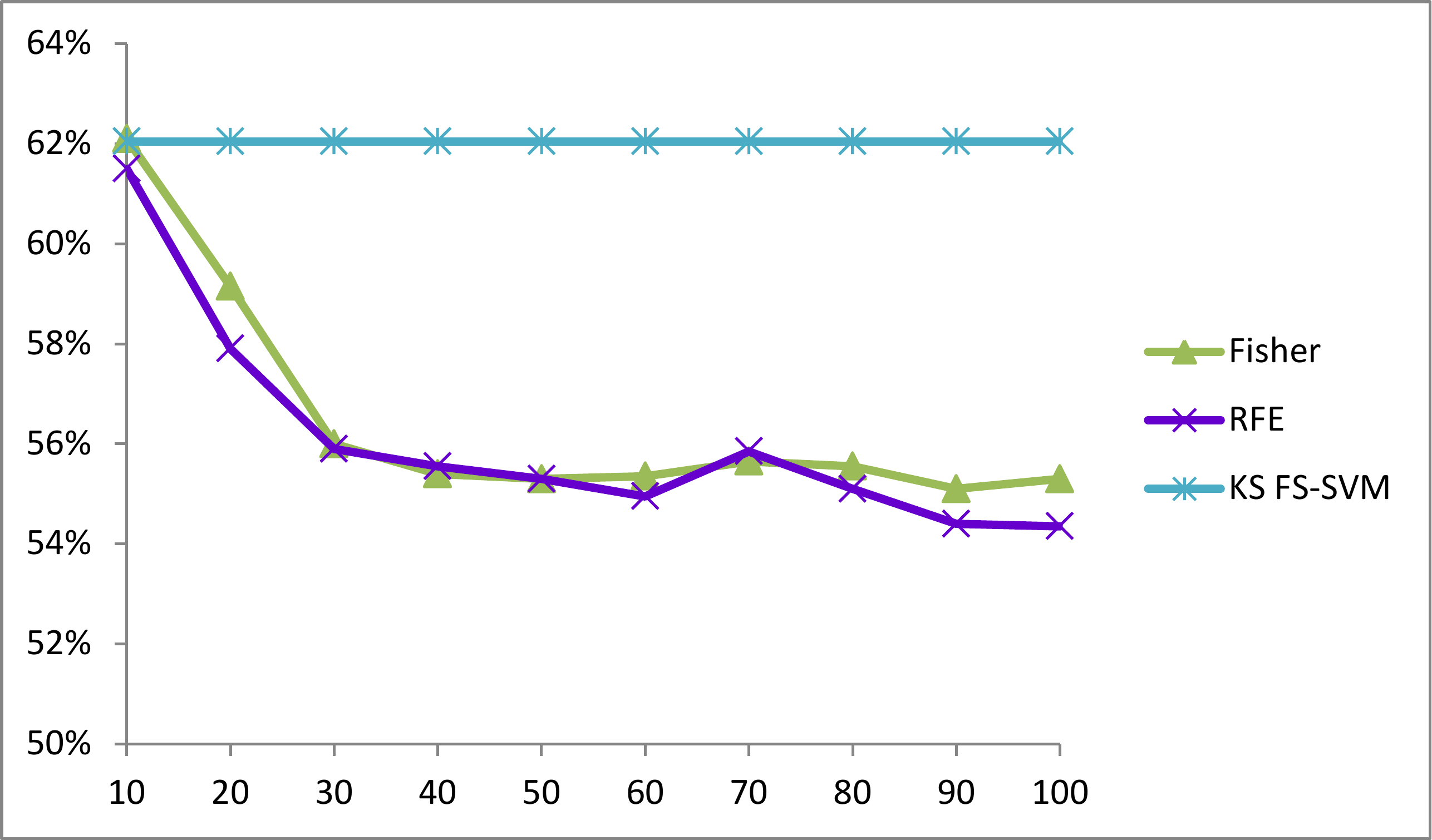}
	\end{minipage}%
	\begin{minipage}{.5 \linewidth}
		\centering
		\includegraphics[width=0.9\textwidth]{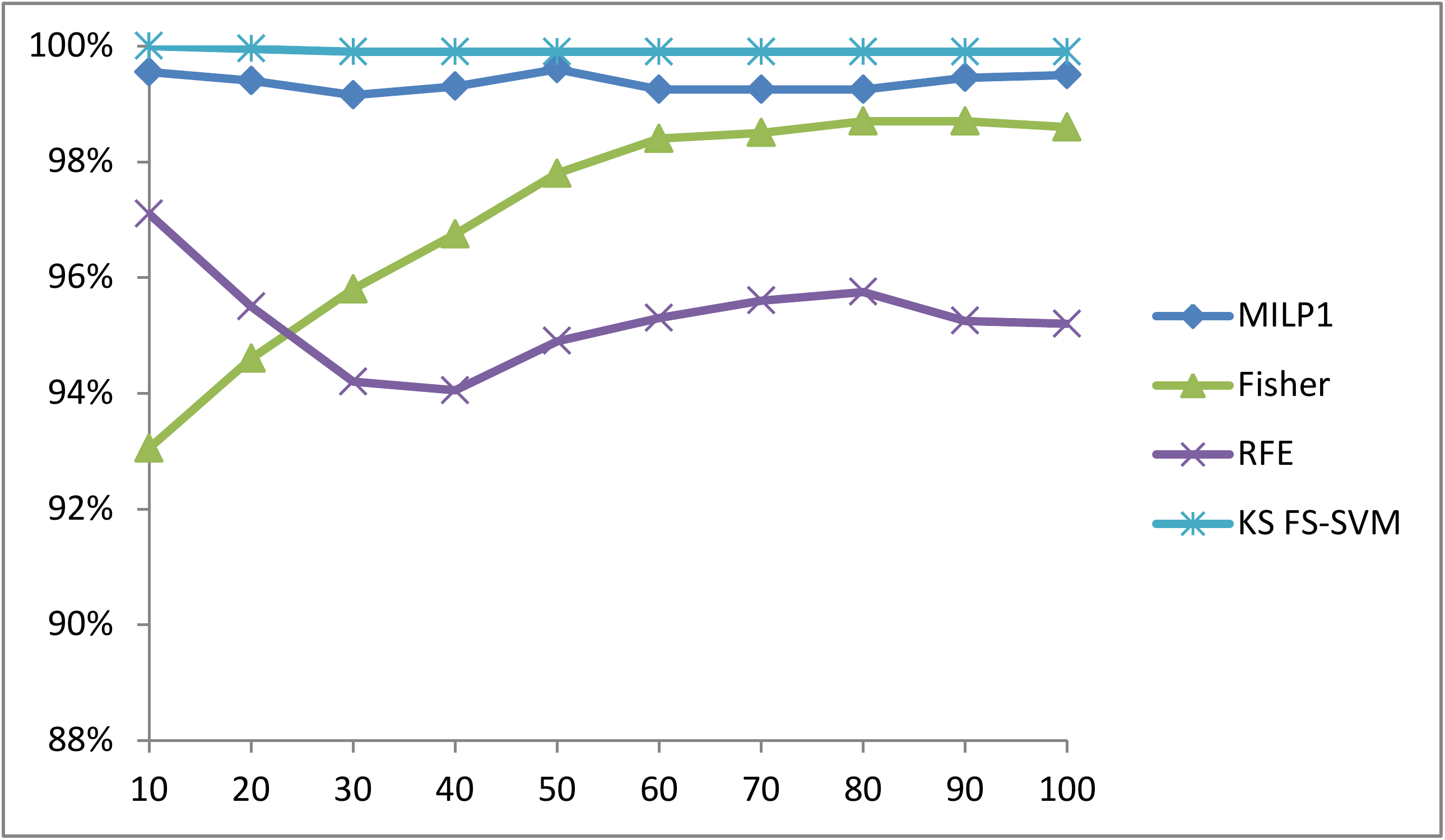}
	\end{minipage} 
	\begin{minipage}[t]{.5\linewidth}
		\caption{\footnotesize{\footnotesize{Average ACC for Madelon dataset.}}\label{Madelon_acc}}
	\end{minipage}%
	\hfill%
	\begin{minipage}[t]{.5\linewidth}
		\caption{\footnotesize{Average ACC for MFeat dataset.}\label{MFeat_acc}}
	\end{minipage}%
\end{figure}

\begin{table}[!htb]
		\begin{minipage}{.5 \linewidth}
			\centering
			\scalebox{0.6}{
				\begin{tabular}{rrrrrrr}
					\hline
					\multicolumn{7}{c}{Madelon m=2000 n=500}       \\
					\hline
					Form.&Av. ACC&Av. AUC & Av. F.& B&C&Time\\
					\hline
					$\ell_1$-SVM & 62.05\% & 62.05\% & 4     &   -    & 2$^{-5}$ & 8.49 \\
					$\ell_2$-SVM & 58.85\% & 58.85\% & 500   &    -   & 2$^{-7}$ & 2.90\\
					LP-SVM & 54.05\% & 54.05\% & 500   &   -    & 2$^6$  & 18.61 \\
					FSV   & 61.78\% & 61.78\% & 2     &  -     &     -  & 20.73 \\
					Fisher-SVM & \textbf{62.10\%} & \textbf{62.10\%} & 10    & 10    & 2$^1$  & 5.85 \\
					RFE-SVM   & 61.50\% & 61.50\% & 10    & 10    & 2$^7$  & 5.35\\
					KS FS-SVM &62.05\%& 62.05\% & 4    & 10    & 2$^{-5}$ & 21.97 \\
					\hline
				\end{tabular}%
			}
		\end{minipage}%
		\begin{minipage}{.5 \linewidth}
			\centering
			\scalebox{0.6}{
				\begin{tabular}{rrrrrrr}
					\hline
					\multicolumn{7}{c}{MFeat m=2000 n=649}      \\
					\hline
					Form. & Av. ACC & Av. AUC & Av. F. & B     & C&Time \\
					\hline
					$\ell_1$-SVM & 99.90\% & 99.50\% & 44.5  & \multicolumn{1}{l}{-} & 2$^7$  & 6.77 \\
					$\ell_2$-SVM & 99.85\% & 99.25\% & 649   & \multicolumn{1}{l}{-} & 2$^{-4}$ & 3.91 \\
					LP-SVM& 99.80\% & 99.22\% & 649   & \multicolumn{1}{l}{-} & 2${^7}$  & 5.02 \\
					MILP1  & 99.70\% & 99.61\% & 144.9 & 163   & -     & 26.31 \\
					FSV   & 50.00\% & 55.22\% & 10    & \multicolumn{1}{l}{-} & -     & 51.74 \\
					Fisher-SVM &98.70\% & 99.28\% & 80    & 80    & 2$^{-4}$ & 2.02 \\
					RFE-SVM   & 97.10\% & 98.39\% & 10    & 10    & 2$^{-7}$ & 0.85 \\
					KS FS-SVM &\textbf{100\%} & \textbf{100\%} & 10    & 10    & 2${^-2}$ & 13.47 \\
					\hline
				\end{tabular}%
			}
		\end{minipage} 
		\begin{minipage}[t]{.45\linewidth}
			\caption{\footnotesize{\footnotesize{Best average ACC and AUC for Madelon dataset.}\label{tab_madelon}}}
		\end{minipage}%
		\hfill%
		\begin{minipage}[t]{.45\linewidth}
			\caption{ \footnotesize{Best average ACC and AUC for MFeat dataset.}\label{tab_mfeat}}
		\end{minipage}%
\end{table}

\section{Concluding remarks \label{concluding}}
We have in this paper provided an MILP model for a classification problem with feature selection based on an SVM approach.  We anlyzed different solution methods for this model, using exact and heuristic algorithms. In addition, in order to validate the model, we have compared the results of the model with respect to others existing in the literature. The main conclusion is that the results are good, stable and, in some cases, they provide an improvement over existing models for small instances and instances with small sample size and big number of features. However, future work should focus on the development of enhancements for the cases with big sample size and big number of features.

\section*{Acknowledgements}
This research  was supported in part by Ministerio de Econom\'ia y Competitividad via grant numbers MTM2013-46962-C02-02 and MTM2016-74983-C2-2-R. Luisa I. Martínez Merino acknowledges that the research reported here was also supported by Universidad
de Cádiz, via PhD grant number UCA/REC02VIT/2014. This work was also partially supported by the Interuniversity Attraction Poles Programme P7/36 \ COMEX initiated by the Belgian Science Policy Office.

\bibliographystyle{abbrvnat}

\end{document}